\documentclass[10pt]{amsart}

\usepackage{amsthm,amsfonts,amssymb}
\usepackage{color}
\usepackage{graphicx}
\usepackage{float}
\usepackage{multirow}
\usepackage{subfigure}
\usepackage[unicode=true,pdfusetitle,
   bookmarks=true, 
   bookmarksnumbered=false, 
   bookmarksopen=false,
   breaklinks=false, 
   pdfborder={0 0 1}, 
   backref=false, 
   colorlinks=false,hidelinks]{hyperref}
\hypersetup{pdfauthor={Name}}

\graphicspath{{figs/}}
\usepackage[ruled,vlined]{algorithm2e}
\usepackage{caption}

\theoremstyle{definition}
\newtheorem{thm}{Theorem}[section]

\newtheorem{exm}[thm]{Example}

\numberwithin{equation}{section}

\newcommand{\abs}[1]{\left\vert#1\right\vert}

\newcommand{\Real}{\mathbb R}
\newcommand\blfootnote[1]{%
  \begingroup
  \renewcommand\thefootnote{}\footnote{#1}%
  \addtocounter{footnote}{-1}%
  \endgroup
}

\newcommand{\zz}[1]{{\color{red}#1}}
\allowdisplaybreaks[1]

\begin{document}
\title[Two-scale networks]{Two-scale Neural Networks for Partial Differential Equations with  Small Parameters}

\author{Qiao Zhuang}
\address[Qiao Zhuang]{Department of Mathematical Sciences, Worcester Polytechnic Institute, Worcester, MA 01609, USA}
\curraddr{School of Science and Engineering, University of Missouri-Kansas City, Kansas City, MO 64110, USA}
\email{qzhuang@umkc.edu}

\author{Chris Ziyi Yao}
\address[Chris Ziyi Yao]{Department of Aeronautics, Imperial College London, London, SW7 2AZ, UK}
\email{chris.yao20@imperial.ac.uk}

\author{Zhongqiang Zhang}
\address[Zhongqiang Zhang]{Department of Mathematical Sciences, Worcester Polytechnic Institute, Worcester, MA 01609, USA}
\email{zzhang7@wpi.edu}

\author{George Em Karniadakis}
\address[George Em Karniadakis]{Division of Applied Mathematics, Brown University, Providence, RI 02912, USA\\
and Pacific Northwest National Laboratory, P.O. Box 999, Richland, 99352, WA, USA}
\email{george\_karniadakis@brown.edu}

\date{\today}
\maketitle

\begin{abstract}
We propose a two-scale neural network method for solving partial differential equations (PDEs) with small parameters using physics-informed neural networks (PINNs). We directly incorporate the small parameters into the architecture of neural networks. 
The proposed method enables solving PDEs with small parameters in a simple fashion, without adding Fourier features or other computationally taxing searches of truncation parameters. Various numerical examples demonstrate reasonable accuracy in capturing features of large derivatives in the solutions caused by small parameters.
\end{abstract}


\section{Introduction}
In this work, we consider physics-informed neural networks (PINNs) for  the following  equation 
\begin{equation}\label{eq:equation-of-interest}
{p}\partial_t u -\epsilon L_2u +  L_0u = f,\quad  \mathbf{x}\in D\subset\Real^d,~{ t\in (0, T ], }
\end{equation}
where some proper boundary conditions and initial conditions are imposed.  Here 
${p}=0$ or $1$, $\epsilon>0$, { $D$ is the spatial domain,} ${ d\in \mathbb{N}^+}$ {and  $T>0$}.  Also, $L_2u$ consists of the leading-order differential operator   and 
$L_0u$ consists of lower-order linear or nonlinear differential operators. 
For example, consider a singular perturbation problem where {${p}=0$}, $L_2 u = {\rm div} ({ a} \nabla u)$ is second-order 
and $L_0u= b\cdot\nabla u + cu$ is first-order. 

Small parameters in the equation often pose extra difficulties for numerical methods, see e.g., \cite{RoosST-b08}. The difficulties come from  one or more sharp transitions  in regions of small volumes, which implies large first-order derivatives
or even large high-order derivatives. When ${p}=0$, $\epsilon>0$ is very small, $L_2u=\Delta u$ and 
$L_0u = \mathbf{1}\cdot \nabla u$, then at least one boundary layer { arises}. 

When  deep feedforward neural networks
are used, they are usually 
trained with stochastic gradient descent methods but do not resolve the issues above as they learn functions with low frequency and small first-order derivatives, see e.g., in  
\cite{basri2020frequency,rahaman2019spectral, xu2019frequency}. 

\subsection{Literature review}
To deal with functions with high-frequency components such as in singular perturbation problems, at least four approaches have been proposed to address this issue:
\begin{itemize}
	\item Adding features in the neural networks:
	Adding random Fourier features is probably the most non-intrusive approach. Adding $\cos(\omega_{i}^\top x), \sin (\omega_{i}^\top x)$'s
	to deep neural networks as approximations of target functions or solutions. 
    With an explicitly specified range for the random frequencies, one can learn a large class of functions. 
	In \cite{WangZC20,LiuCaiXu20,LiXuZhang20,CaiLiLiu20}, frequency ranges are scheduled to represent complicated solutions to partial differential equations. 
	See also \cite{YangGu22,LiXuZhang23,Zhang23helmholtz} for more elliptic type multi-scale PDEs. 	The Fourier feature networks are also used in  \cite{WangWP21,wang2021eigenvector}, and see \cite{wang2023experts} for a  review.
	\item Enhanced loss by adding first-order and higher-order derivatives:      Another approach is to include the gradient information in the loss function, e.g., in 
	\cite{Czarnecki17Sobolev,hoffman2020robust,lammleBCR2023sequential}. 
	This approach has been applied to solve { PDEs}, e.g., in \cite{ son2021sobolev,yuLMK22}.
	\item Adaptive weights: In the loss function, adaptive weights are assigned to have a better balance of each squared term in the least-squares formulation, e.g., 
	self-adaptive \cite{mclenny23self-adaptive}, 
	  	 attention-based weights \cite{anagnostopoulos2023residualbased}, binary weights \cite{Gu2021selfpaced}. See also \cite{xiang2022self,zhang2023dasa,maddu2022inverse}.
	\item Resampling: Sampling points in the loss function can be made adaptive based on the residuals at sampling points during the training process such as in \cite{lu2021deepxde,wu2023comprehensive,gao2022failure,gao2023failure}.
     Also, in   \cite{mohammad2021efficient,tang23das-pinn,tang2023adversarial}, the density function for sampling during the training is computed using the idea of importance sampling.
     The density function for re-sampling is approximated by Gaussian mixtures and formulated by finding the min-max of the loss function via fixed-point iterations.

\end{itemize} 
	Other approaches are not particularly for the high-frequency issue, such as \cite{jagtap2020adaptive} using adaptive activation functions and 
\cite{wang2022respecting} by respecting causality, {continual learning \cite{howard2023stacked} and sequential training \cite{HanLee23hierarchical}.}

%

 \subsection{Related works}
 One relevant approach to accommodate multiple scales using neural networks is the asymptotic-preserving neural network (APNN) \cite{Jin23,LuWangWu23,BerLP22}. It involves two-scale expansions and transforming the underlying equation into several equations that may be scale-independent. For example, in \cite{Jin23, LuWangWu23}, the loss function is designed based on equations from a macro-micro decomposition, while in \cite{BerLP22}, the loss function is rooted in a system in a macroscopic form that corresponds to the governing PDEs.
In \cite{GuanYangCLG23}, dimension-augmentation is used to enhance the accuracy, which includes adding the Fourier features.
Our approach is closely related to APNN as we explicitly build-in a scale as an extra dimension of input for the neural networks.
Another relevant approach is multi-level neural networks \cite{multilvNN23}, addressing the multiple scales in the solution by using a sequence of neural networks with increasing complexity: at each level of the process, it uses a
new neural network, to generate a correction corresponding to the error in the previous approximation.  Two ingredients are essential in mitigating the effects of high-frequency components: adding Fourier features in the networks and the extra level of refining residuals. A similar approach is in \cite{fang2023ensemble} while the scale of residuals is not explicitly incorporated in the solution.

\subsection{Contribution and significance of the work}

We directly incorporate the scale parameter into the neural networks' architecture, establishing a \textit{two-scale} neural network method. In particular, for a small scale parameter $\epsilon$, the networks employ auxiliary variables related to the scale of ${\epsilon}^{\gamma}$ ($\gamma \in \Real$) designed to capture intricate features that involve large derivatives such as boundary layers, inner layers, and oscillations. 
Specifically, 
we use the feedforward neural network of the form  
\begin{equation}\label{eq:nn-with-aux-general}
w(\mathbf{x}) N({ t},\mathbf{x}, \epsilon^\gamma(\mathbf{x} -\mathbf{x}_c), \epsilon^\gamma), \quad\gamma<0,
\end{equation}
where $\mathbf{x}_c\in D$.  Here $w(\mathbf{x})$ is $1$ when boundary layers arise or a simple function satisfying the Dirichlet boundary condition if no boundary layer arises. 
Without an explicit statement, $\mathbf{x}_c$ is the center of the domain $D$. 
To assess the performance of the two-scale neural networks, we will compare their results with those of the following \textit{one-scale} neural networks
\begin{equation}\label{eq:nn-wn-aux}
w( \mathbf{x} )N({ t},\mathbf{x}).
\end{equation}
{ We note that when addressing steady-state problems, $t$ is removed from \eqref{eq:nn-with-aux-general} and \eqref{eq:nn-wn-aux}.}

{
The network \eqref{eq:nn-with-aux-general} is motivated as follows.
First, it can handle large derivatives for training neural networks. By 
taking derivatives in $x$, we readily obtain a factor of $\epsilon^{\gamma}$.
Thus, the networks 
\eqref{eq:nn-with-aux-general} can facilitate  
learning large derivatives or higher frequencies, especially 
 with low-order optimizers.
Second, the network \eqref{eq:nn-with-aux-general} is analogous to the classical scale expansion of $u$ in $\epsilon$ while we do not have explicit expansions.
Suppose  that 
$u=\sum_{k=0}^\infty u_k \epsilon^{k\beta}$, where $\beta>0$ is problem-dependent.  The intuition is that the scale of $u_k$ is about the order of $\epsilon^{-k\beta}$. 
For the network 
\eqref{eq:nn-with-aux-general}, 
 Taylor's expansion in $\mathbf{x}$ leads to a scale expansion in $\epsilon$. %
Consider the one-dimensional function and the first two derivatives of 
\eqref{eq:nn-with-aux-general}. 
The first derivative is of order 
$\epsilon^{\gamma}$. 
%
%
The second derivative is of the order 
$\epsilon^{2\gamma}$.
 For singular perturbation problems in this work, $\beta$ is usually $1/2$ or $1$. We then choose $\gamma=-1/2$ in all cases. In some cases,  $\gamma=-1$  may be a better choice but still can be covered if we use a smooth activation function since $(\epsilon^{-\frac{1}{2}})^2= \epsilon^{-1}$. It may enhance the flexibility of networks by using $\gamma=-1/2$ in all cases, especially when a smooth activation function implies a valid Taylor's expansion in $\epsilon$.  

The point $\mathbf{x}_c$ is chosen as the center of the spatial domain $D$. However,  other choices are possible. For example, we can take $\mathbf{x}_c$ to be an endpoint if there is only a boundary layer at one endpoint. 
However, in most cases, we do not know where large derivatives are. In these cases, we expect that 
the linear inner layer of the feedforward neural network can help find the best linear combinations of 
$\epsilon^{\gamma}(\mathbf{x}
-\mathbf{x}_c)$ and $\epsilon^\gamma$ in training and consequently have locations of large derivatives.  
}

The significance of this work is that the new architecture of networks enables solving problems with small parameters in a simple fashion. 
For example, we do not need to 
determine how many equations in APNN, or levels in multilevel neural networks \cite{multilvNN23}, which are problem-dependent.  We also avoid the use of Fourier features, a lot of which may be required in two dimensions and even higher dimensions. 

We compare our method with multilevel neural networks \cite{multilvNN23}, underscoring our method's capability to capture large derivatives without tuning scaling factors.
{The design of the two-scale neural networks enables the neural networks to implicitly pick up the scale expansions to accommodate the multiple scales. In contrast, APNN accommodate the multiple scales via explicit asymptotic expansions, resulting in decomposition strategies such as micro-macro and even-odd decomposition. The selection of decomposition strategies and truncation terms of asymptotic expansions is critical to the APNN model's outcomes. Owing to the variability these choices introduce, we refrain from using APNN results as a benchmark to be compared within our study.}

The effectiveness of the proposed two-scale approach in handling large derivatives is demonstrated by comprehensive numerical examples in Section \ref{sec:numerical}. These include complex cases such as problems with dual boundary layers in 1D and 2D, viscous Burgers equations, and Helmholtz equations exhibiting oscillations along $xy$ and radial directions, where the corresponding one-scale method falls short in delivering accurate simulations without special treatments.
 %

 \section{Formulation and algorithms}

Under the framework of PINNs, the neural network solution is obtained by minimizing the loss function of residuals. The loss function of the considered problem \eqref{eq:equation-of-interest} at the continuous level is 
\begin{equation}\label{eq:loss_conti}
\int_{[0,T]\times D}r^2(t,x) + \alpha \int_{\partial D } (u-g)^2 +\alpha_1 \int_D \vert \nabla u\vert^2 + \beta \cdot\text{loss of initial conditions}
\end{equation}
where $\alpha \geq 1, \alpha_1 \geq 0$ and $\beta\geq 0$ and 
\begin{equation*}
 r(t,x)= 
{ p}\partial_t u - \epsilon L_2u +  L_0u - f.
\end{equation*}
The discrete formulation corresponding to \eqref{eq:loss_conti} is
\begin{align}
&\frac{1}{N_c}\sum_{i=1}^{N_c}\left|r\left(t_r^i, \boldsymbol{x}_r^i\right)\right|^2+\frac{\alpha}{N_b}\sum_{i=1}^{N_b}\left|u\left(t_b^i, \boldsymbol{x}_b\right)-g^i\right|^2+ \frac{\beta}{N_0}\sum_{i=1}^{N_0}\left|u\left(t_0, \boldsymbol{x}_0^i\right)-u_0^i\right|^2 
+\frac{\alpha_1}{N_c}\sum_{i=1}^{N_c}\vert \nabla u(t_r^i, \boldsymbol{x}_r^i) \vert^2 \label{eq:loss-general},
\end{align}
where $\{{t}_r^i, \boldsymbol{x}_r^i \}_{i=1}^{N_c}$ specify the collocation points in the interior time $(0,T]$ and spatial domain $D$, $\{ t_b^i\}$ are time collocation points at $\boldsymbol{x}=\boldsymbol{x}_b$, $\boldsymbol{x}_0^i$ are spacial collocation points at $t=t_0$.

{In this paper, we utilize the Adam \cite{kingma2014adam} optimizer to train the neural network parameters.}
{To enhance the accuracy of capturing solutions to problems with small parameters $\epsilon$, we employ a successive training strategy. This approach progressively optimizes the weights and biases of the neural networks, starting with modestly larger $\epsilon_0$ as the initial guess.} We summarize the 
successive training method in Algorithm 
\ref{alg:succesive-training-two-scale}. 
{We remark that this training strategy is similar to continual learning in \cite{howard2023stacked} and sequential training in \cite{HanLee23hierarchical}. However, we fix the size of neural networks, while in \cite{HanLee23hierarchical,howard2023stacked}, the network size increases during different stages of  training.} 

\begin{algorithm}[!htb]

\SetAlgoLined 
\KwData{Training set, adaptive learning rates, $\epsilon$ and other parameters from the PDE} 
\KwResult{Optimized weights and bias of the neural networks}

{Pick $\epsilon_0= O(0.1)$ if $\epsilon$ is very small.}
Initialize the weights and biases randomly.

Step 1.  {Solve the PDE with parameter $\epsilon_0$ instead of $\epsilon$. 
Use the loss function above
 with $\epsilon_0$ in place of $\epsilon$. 
 Train to obtain weights and biases of the neural networks 
 \eqref{eq:nn-with-aux-general}
  with  $\epsilon=\epsilon_0$.

 Step 2. Pick a positive integer $\ell>1$. 
 While $\epsilon_0\geq \ell \epsilon$\\
 \qquad \qquad set $\epsilon_0=\epsilon_0/\ell$ and go to Step 1;\\
\qquad end
 
 If 
 $\epsilon_0<\ell\epsilon$, set 
 $\epsilon_0=\epsilon$ and go to Step 1. In these iterations, initialize the weights and biases obtained from the last iteration.}

 Step 3. (Fine-tuning). Use Adam with a smaller learning rate {or a second-order optimizer}  to further optimize the weights and biases. 

 Step 4. Stop the process if the maximal epoch number is reached.

\caption{Successive training of two-scale neural networks for PDEs with small parameters}
\label{alg:succesive-training-two-scale}
\end{algorithm}

%

\section{Numerical results} \label{sec:numerical}
In this section, we employ the two-scale neural network method for extensive numerical tests on problems where small parameters lead to solutions with large derivatives. The numerical tests cover one-dimensional (1D) ODEs with one (Examples \ref{exm:1p1}{, \ref{exm:burgers1d-steady} and \ref{exm:1d_nonlinear}}) and two (Example \ref{exm:1p4}) boundary layers, 1D viscous Burgers equations with an inner layer (Example \ref{exm:burgers}), two-dimensional (2D) steady-state convection-diffusion problems featuring two boundary layers (Example \ref{exm:zzm}), and 2D Helmholtz problems characterized by oscillations acting as inner layers (Example \ref{exm:1p14}). The numerical results indicate that the two-scale neural network method can provide reasonable accuracy in capturing features arising from large derivatives in solutions. These features include boundary layers, inner layers, and oscillations. 

For all the numerical tests, we let $\gamma=-\frac{1}{2}$ in \eqref{eq:nn-with-aux-general}. 
{ Unless stated otherwise, regardless of whether Algorithm \ref{alg:succesive-training-two-scale}
 is used},  we adopt the { default setup up of the} piecewise constant learning rates scheduler, the uniform distribution for collocation points, and the ${ \tanh}$ activation function, as listed in Table \ref{tab:lr_pt}. 
{All examples are implemented with JAX and the Adam optimizer unless stated otherwise.}
In some tests, we employ the successive training strategy listed in Algorithm \ref{alg:succesive-training-two-scale}.  {The parameters $\epsilon_0$, $\ell$ are manually set and are problem-dependent. We let the starting $\epsilon_0=0.1$ and $\ell=10$ unless stated otherwise. In Examples \ref{exm:1p1} and \ref{exm:1p4}, when $\epsilon=10^{-5}$, we use $\ell=10$ for the successive  training until $\epsilon_0\geq 10^{-4}$ and 
then $\ell=2$ to reach $\epsilon=10^{-5}$
.}

{ To further demonstrate the robustness of numerical results, we report the statistical metric $\overline{e}\pm \sigma$ for Examples \ref{exm:1p1} and \ref{exm:1p4}. Here, $\overline{e}$ is the average absolute error and $\sigma$ is the standard deviation.
}

\begin{table}[!ht]
    \begin{tabular}{c|c|c} \hline 
        Learning Rates ($\eta$) & Collocation Points & Activation Function \\ \hline 
        $\leq 10000$ steps: $\eta=10^{-3}$ (1D problems), $\eta=10^{-2}$ (2D) &  &  \\      
        $10000$ to $30000$ steps: $\eta=5\times 10^{-3}$ & uniform & $\tanh$ \\
        $30000$ to $50000$ steps: $\eta=10^{-3}$ & distribution &  \\
        $50000$ to $70000$ steps: $\eta=5\times 10^{-4}$ & &  \\
        $\geq 70000$ steps: $\eta=10^{-4}$ & &  \\
    \end{tabular}
    \caption{{ Default setup (unless stated otherwise) of the }learning rates scheduler, collocation points distribution, and activation function}
    \label{tab:lr_pt}
\end{table}


\begin{exm}[1D ODE with one boundary layer]\label{exm:1p1}
	\begin{equation*}
	-\epsilon u'' + 2  u'  =3, \quad  { 0<x<1,}   \quad 
	u(0)=0, \quad u (1)=0. 
	\end{equation*}	
\end{exm}  

\begin{figure} [!ht]
\centering
\subfigure[exact and NN solutions]
{
\includegraphics[width=0.22\textwidth]{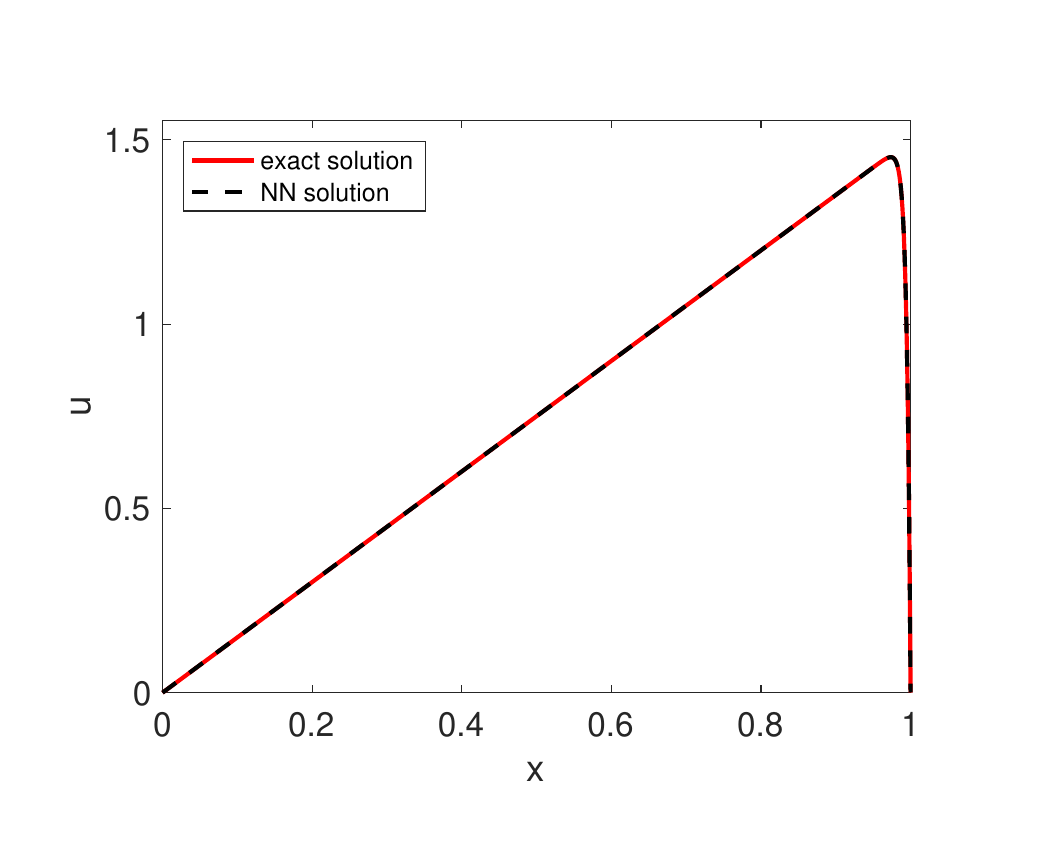}}
\subfigure[absolute error]
 { 
 \includegraphics[width=0.22\textwidth]{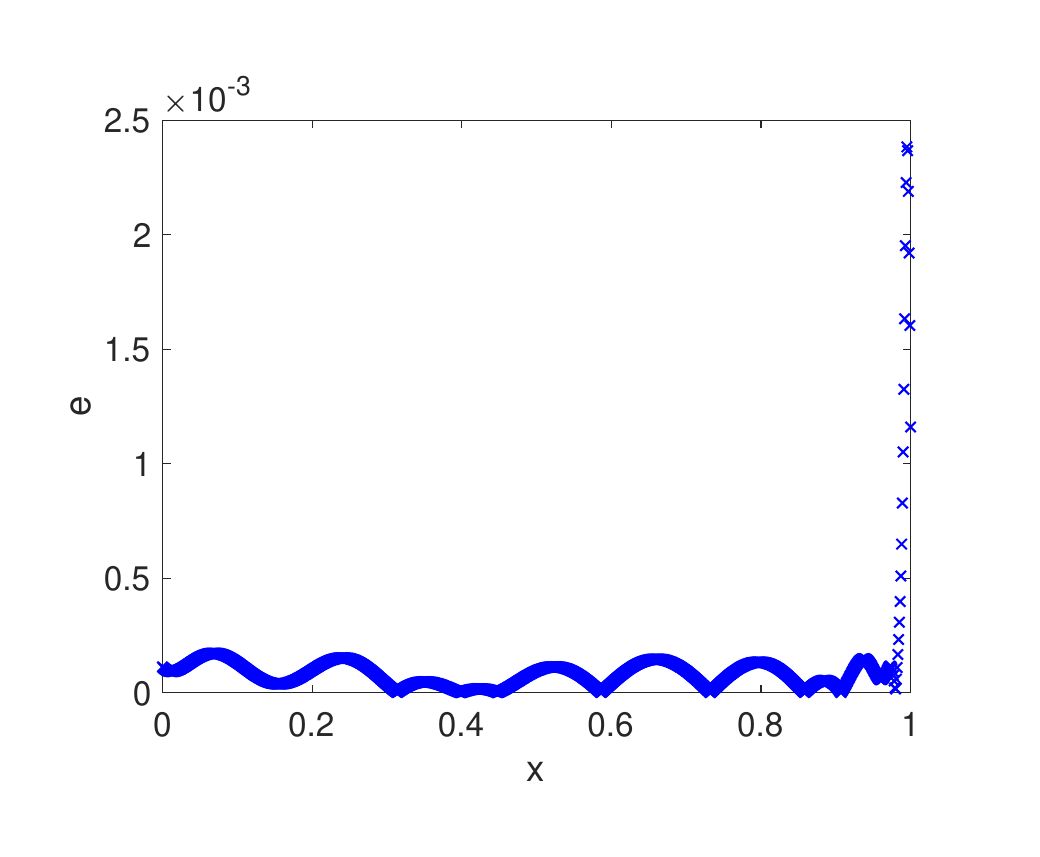}}
 \subfigure[relative error around the boundary layer]
{
\includegraphics[width=0.22\textwidth]{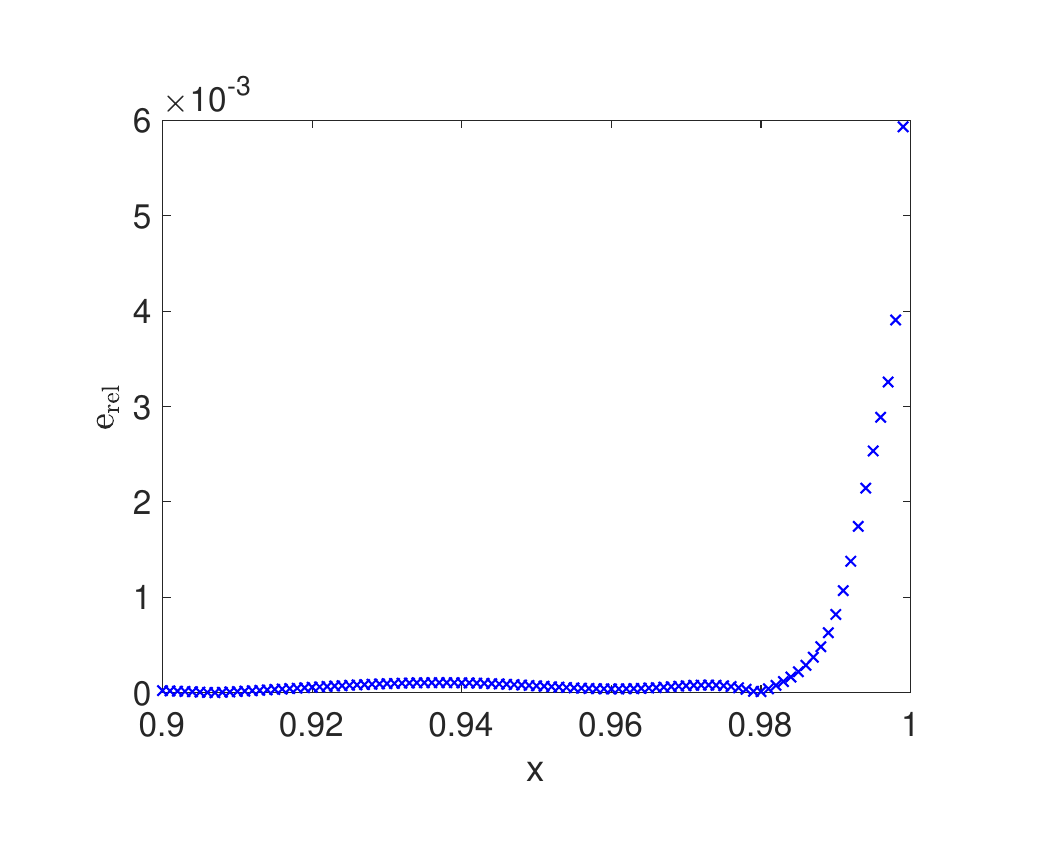}}
\subfigure[relative error away from the boundary layer]
 { 
\includegraphics[width=0.22\textwidth]{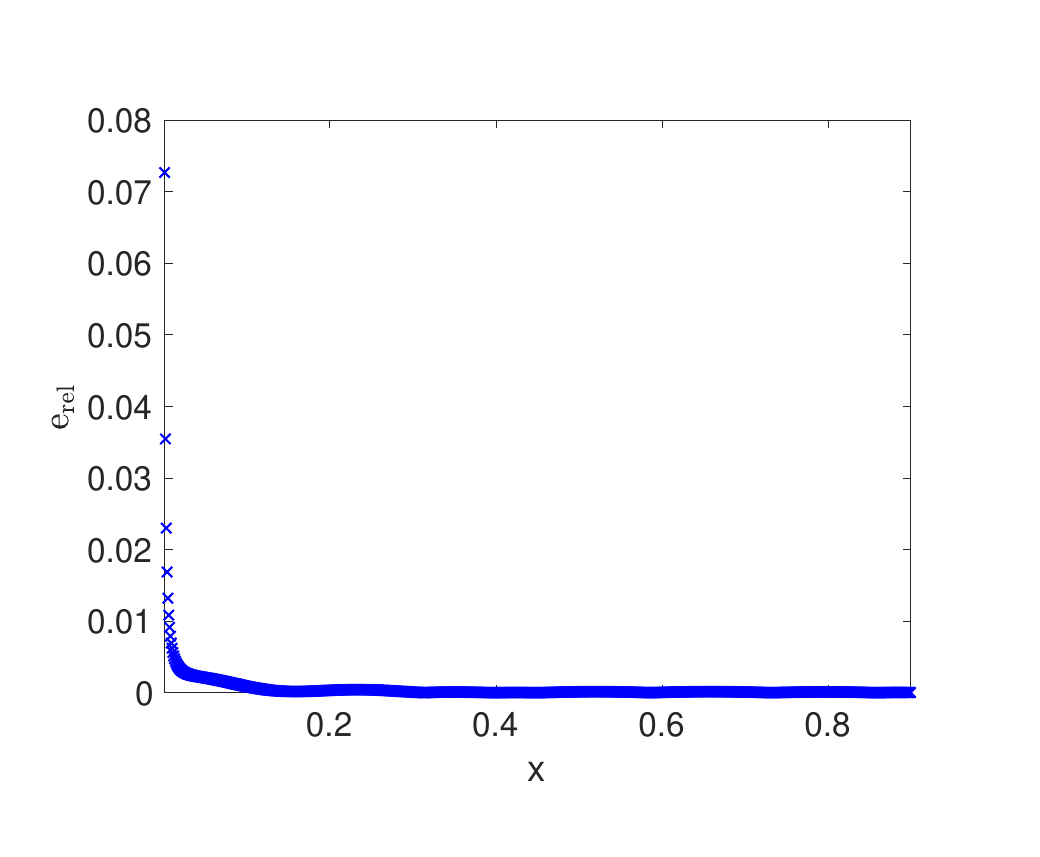}}
 \caption{Numerical results for Example \ref{exm:1p1} when $\epsilon=10^{-2}$ using { $N(x, (x-0.5)/\sqrt{\epsilon}, 1/\sqrt{\epsilon})$}, { with parameters specified in Table \ref{tab:exm1p1_para}  }.}
 \label{fig:results_exm1p1_1e-2} 
 \end{figure}

{
The solution has a boundary layer at $x=1$.}
The exact solution to this problem is  $$ u(x)= \frac{3}{2}\left(x-\frac{\exp(-2(1-x)/\epsilon)
		-\exp(-2/\epsilon)}{1-\exp(-2/\epsilon)}\right).$$ 
  We use the two-scale neural networks $N(x, (x-0.5)/\sqrt{\epsilon},1/\sqrt{\epsilon})$ to solve the ODE problem. The neural network (NN) size is $(3,20,20,20,20,1)$.
  {We employ the successive training strategy in Algorithm \ref{alg:succesive-training-two-scale} starting with $\epsilon_0=10^{-1}$, and successively solve the problem for $\epsilon$ as small as $10^{-5}$, following the sequence $\epsilon_0=10^{-1},10^{-2},10^{-3},10^{-4},
  5\times10^{-5},2.5\times 10^{-5},1.25\times 10^{-5}, 10^{-5}$. The parameters in the discrete loss function in \eqref{eq:loss-general} and hyper-parameters of the training are summarized in Table \ref{tab:exm1p1_para} }.
  We collect the numerical results in Figures \ref{fig:results_exm1p1_1e-2} { and \ref{fig:exm1p1_1p4_loss}(a)} for the case where $\epsilon=10^{-2}$. Observing from Figure \ref{fig:results_exm1p1_1e-2}(b) to (d), the NN  solution matches the exact solution very well, with the relative errors around the boundary layer $x=1$ only at the magnitude of $10^{-3}$. 

\begin{figure} [!ht]
\centering
\subfigure[exact and NN solutions]
{
\includegraphics[width=0.23\textwidth]{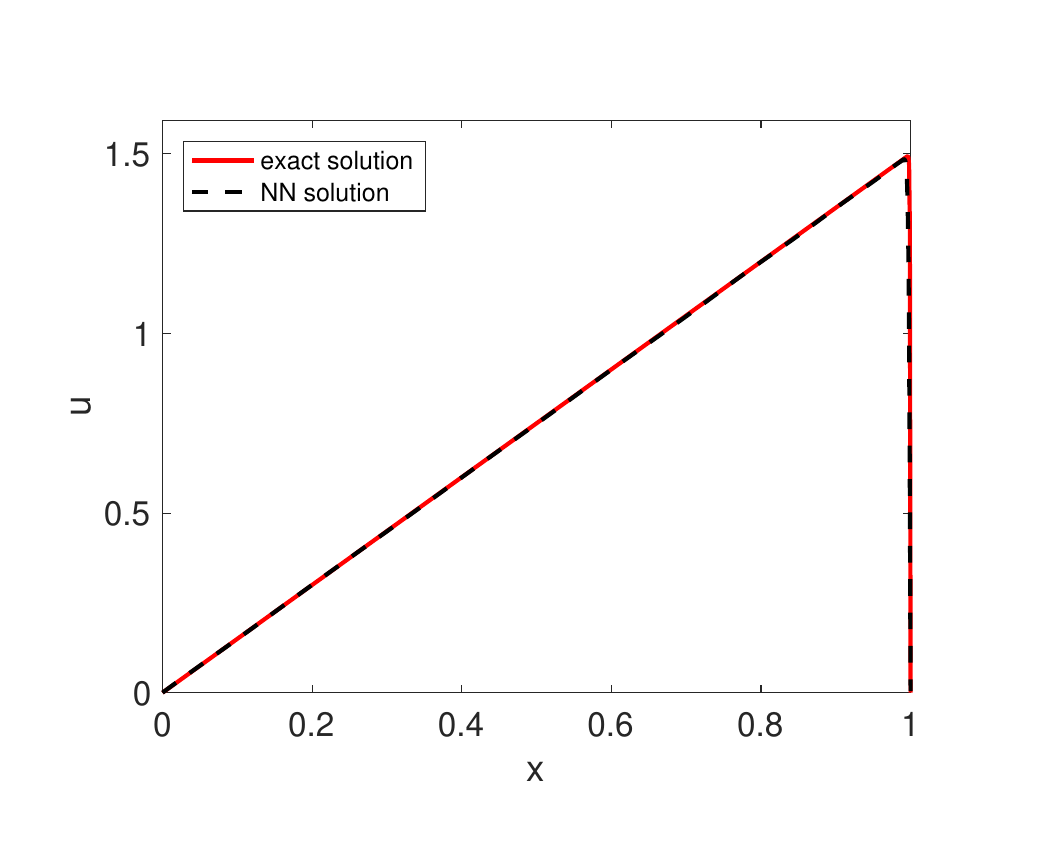}}
\subfigure[absolute error]
 { 
\includegraphics[width=0.23\textwidth]{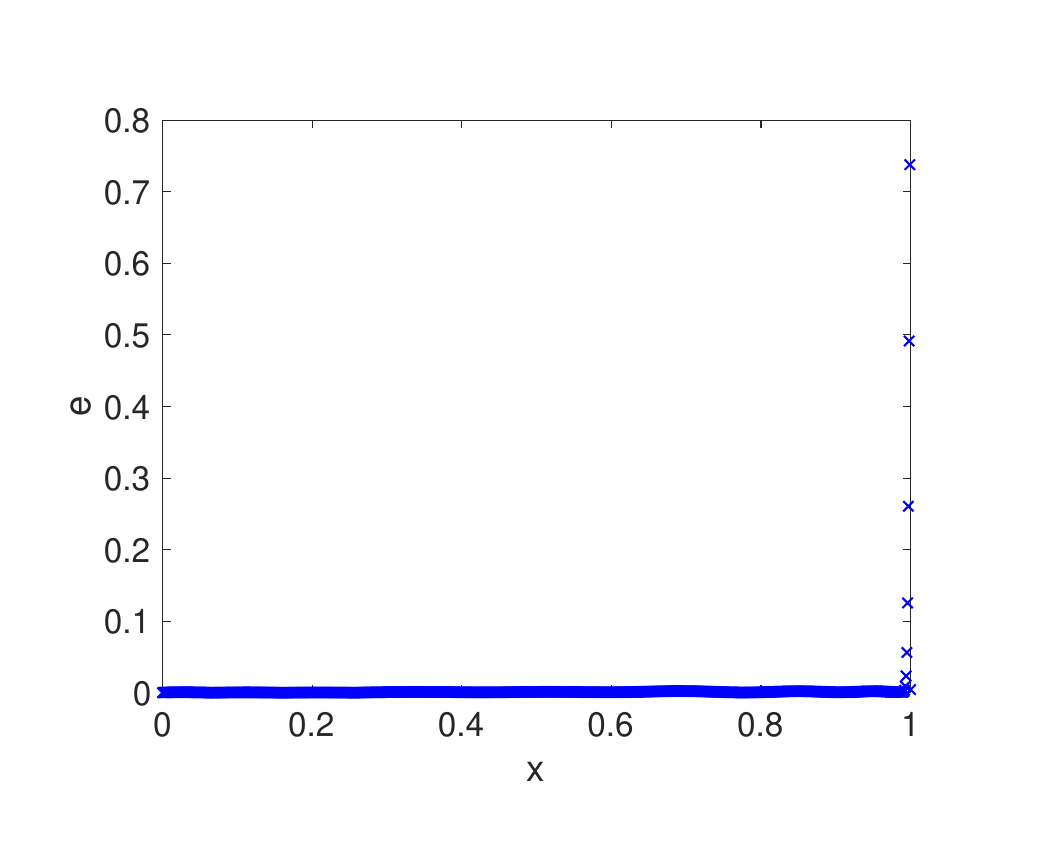}}
 \subfigure[relative error around the boundary layer]
{
\includegraphics[width=0.23\textwidth]{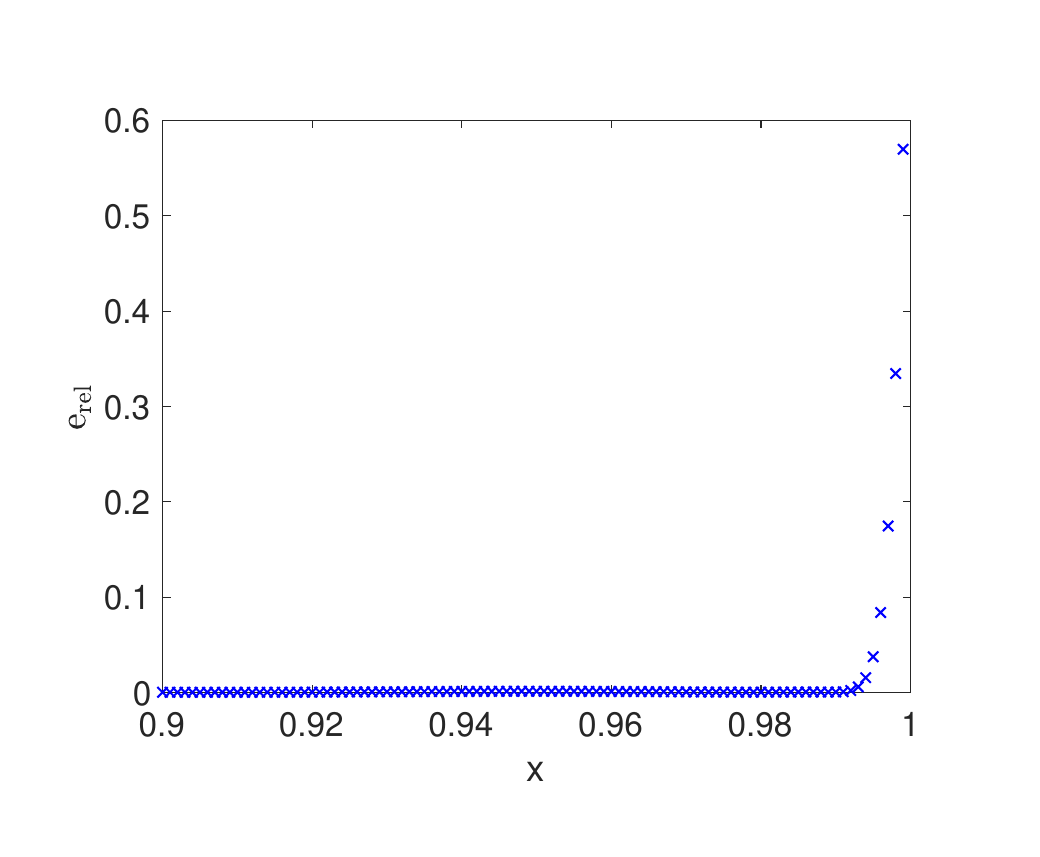}}
\subfigure[relative error away from the boundary layer]
 { 
\includegraphics[width=0.23\textwidth]{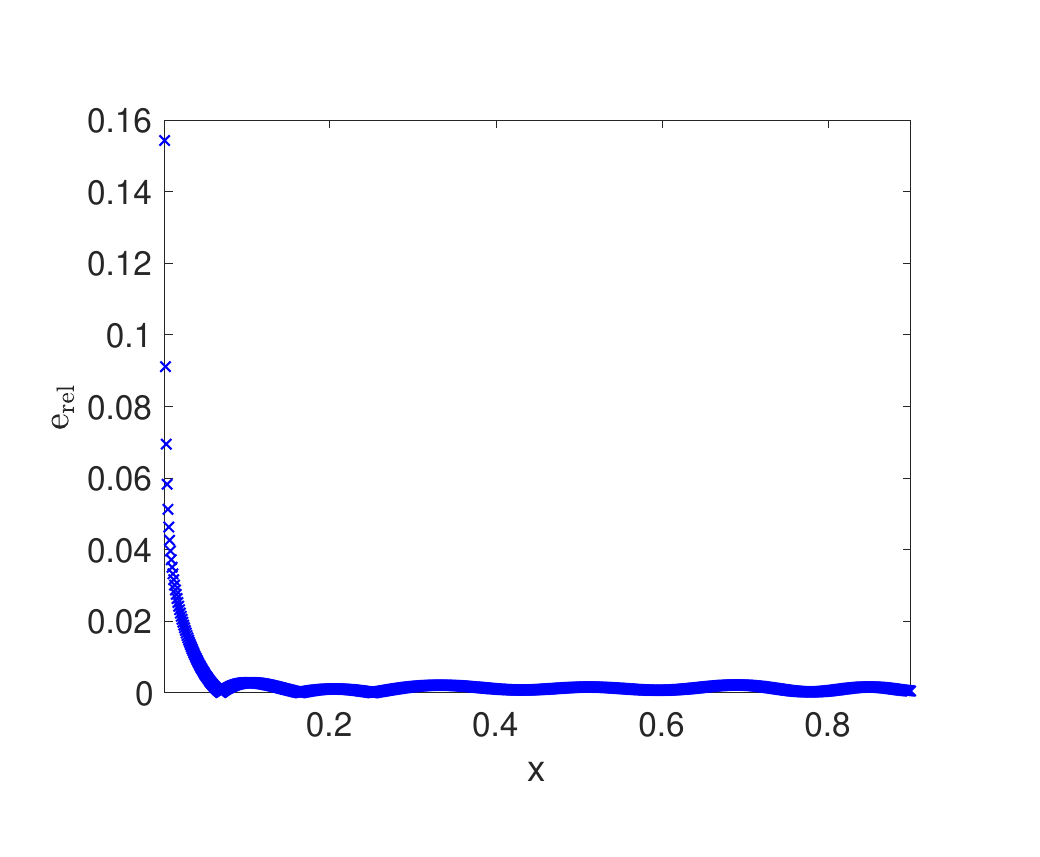}}
 
 \caption{Numerical results for Example \ref{exm:1p1}  when $\epsilon=10^{-3}$, using { $N(x, (x-0.5)/\sqrt{\epsilon}, 1/\sqrt{\epsilon})$}, { with parameters specified in Table \ref{tab:exm1p1_para}  }.}
 \label{fig:results_exm1p1_1e-3} 
 \end{figure}

\begin{figure} [!ht]
\centering
\subfigure[exact and NN solutions for $\epsilon=10^{-4}$]
{
\includegraphics[width=0.23\textwidth]{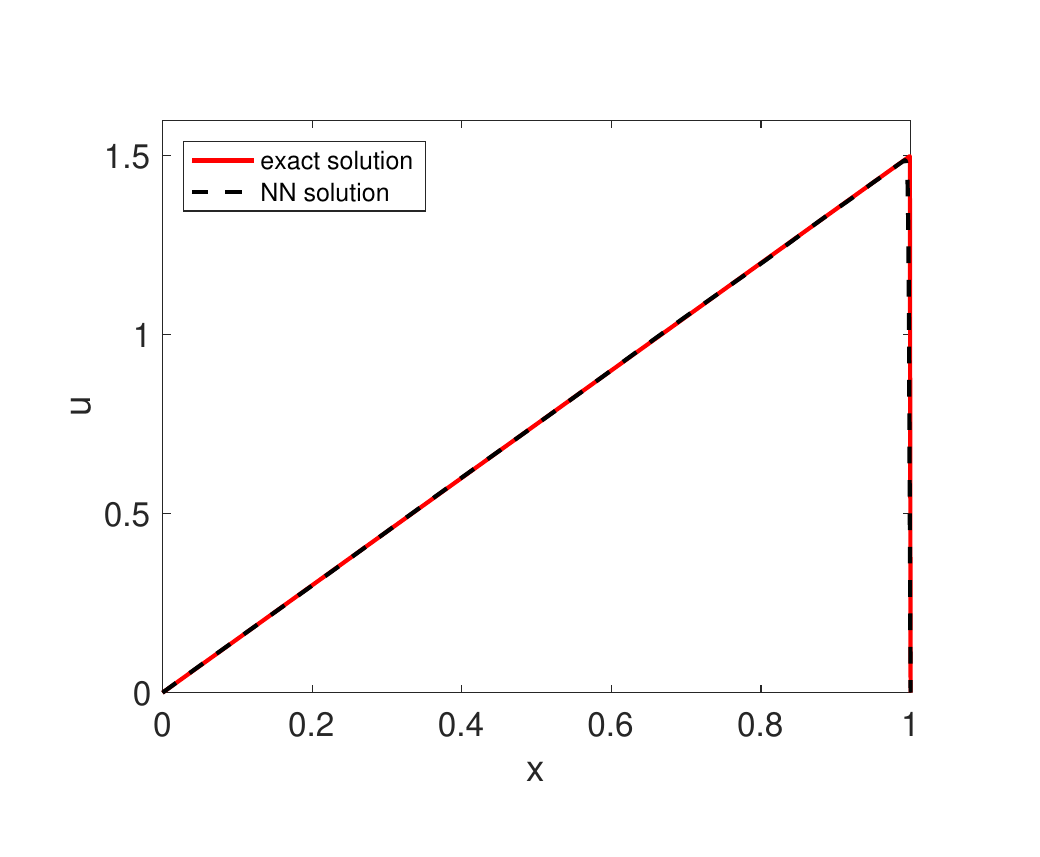}}
\subfigure[absolute error]
 { 
\includegraphics[width=0.23\textwidth]{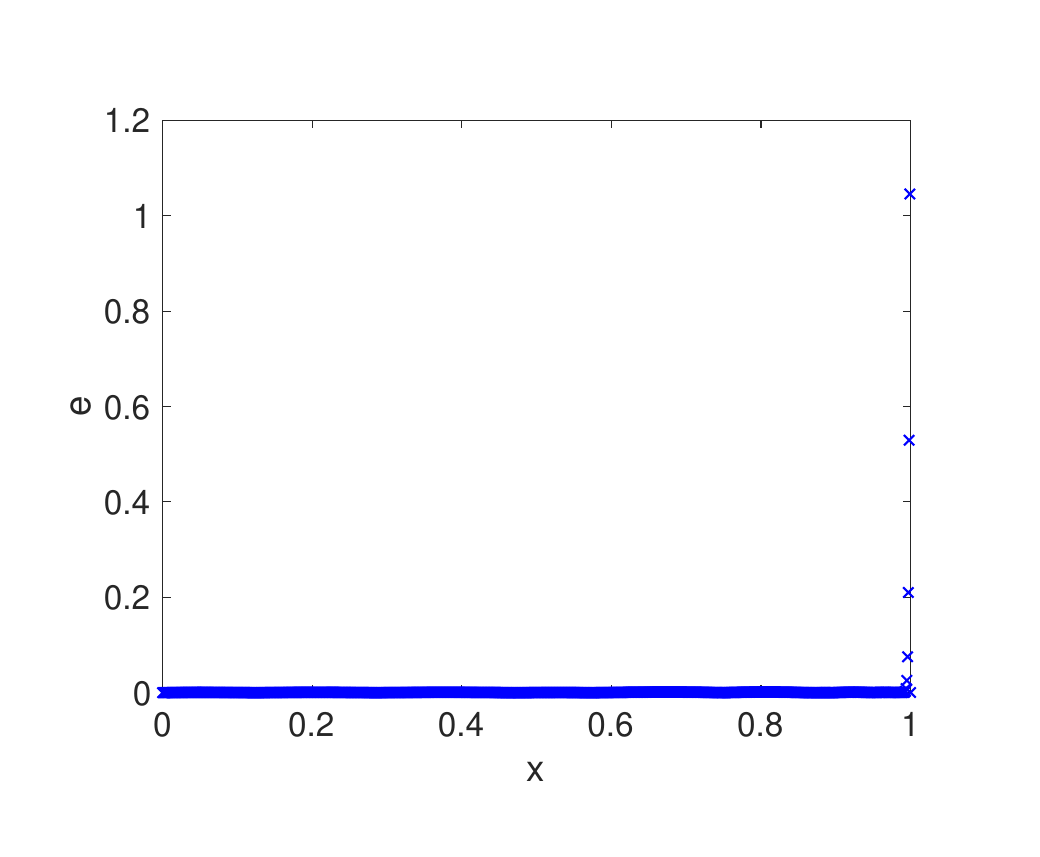}}
 \subfigure[relative error around the boundary layer]
{
\includegraphics[width=0.23\textwidth]{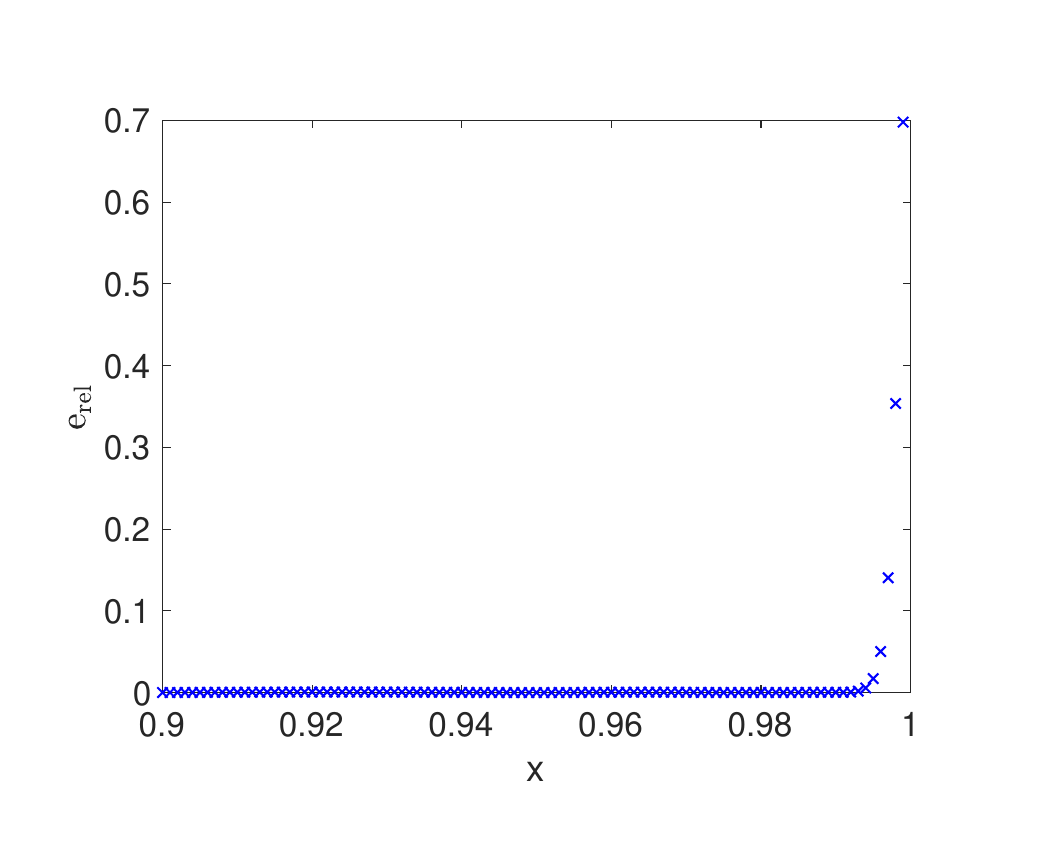}}
\subfigure[relative error away from the boundary layer]
 { 
\includegraphics[width=0.23\textwidth]{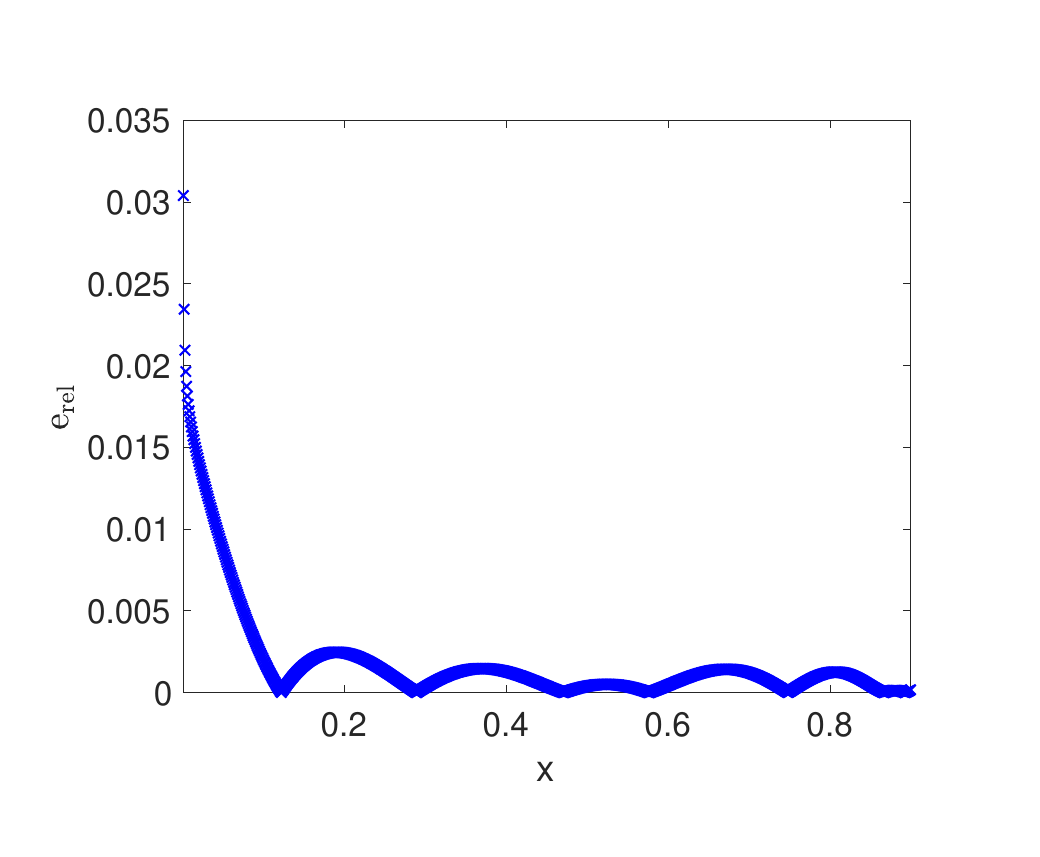}}

\subfigure[exact and NN solutions for $\epsilon=10^{-5}$]
{
\includegraphics[width=0.23\textwidth]{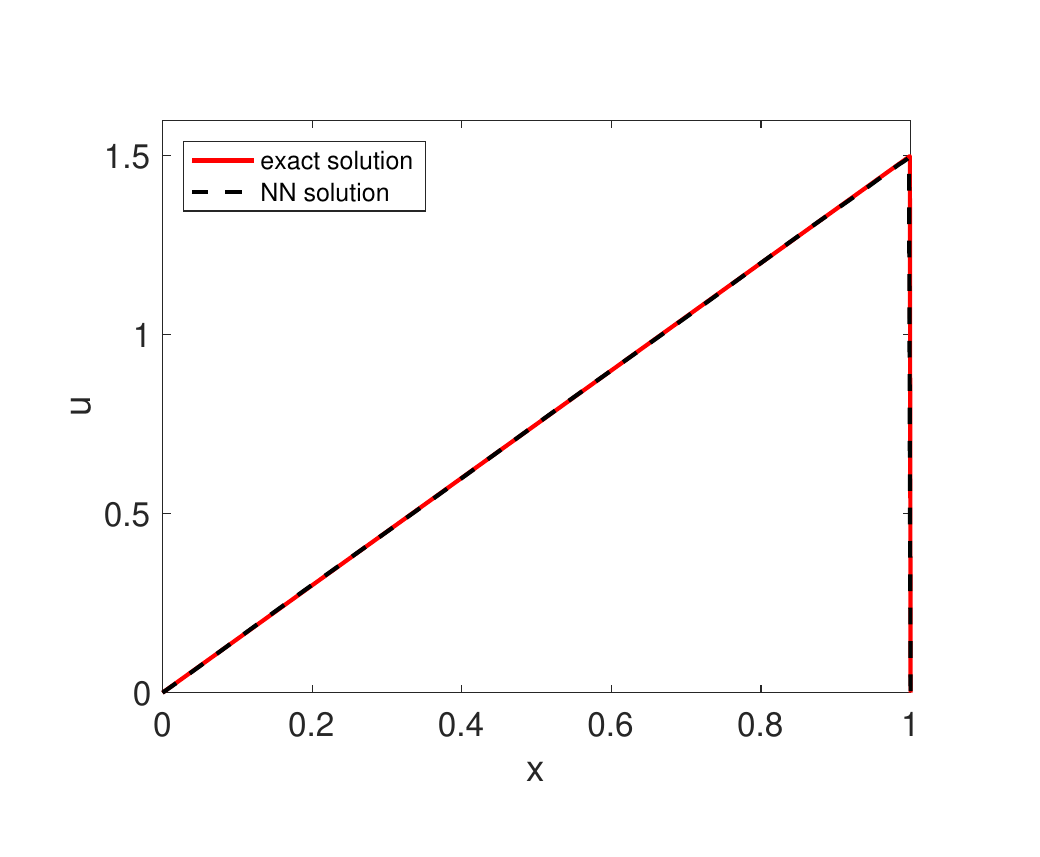}}
\subfigure[absolute error]
 { 
\includegraphics[width=0.23\textwidth]{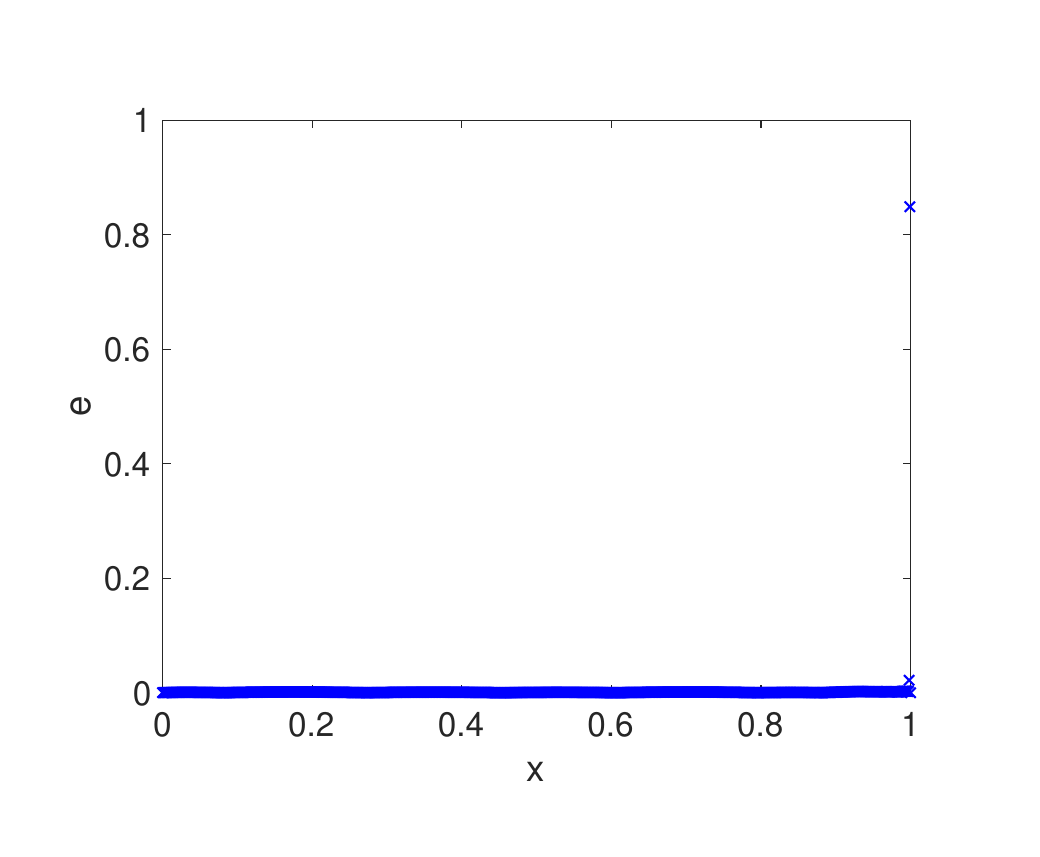}}
 \subfigure[relative error around the boundary layer]
{
\includegraphics[width=0.23\textwidth]{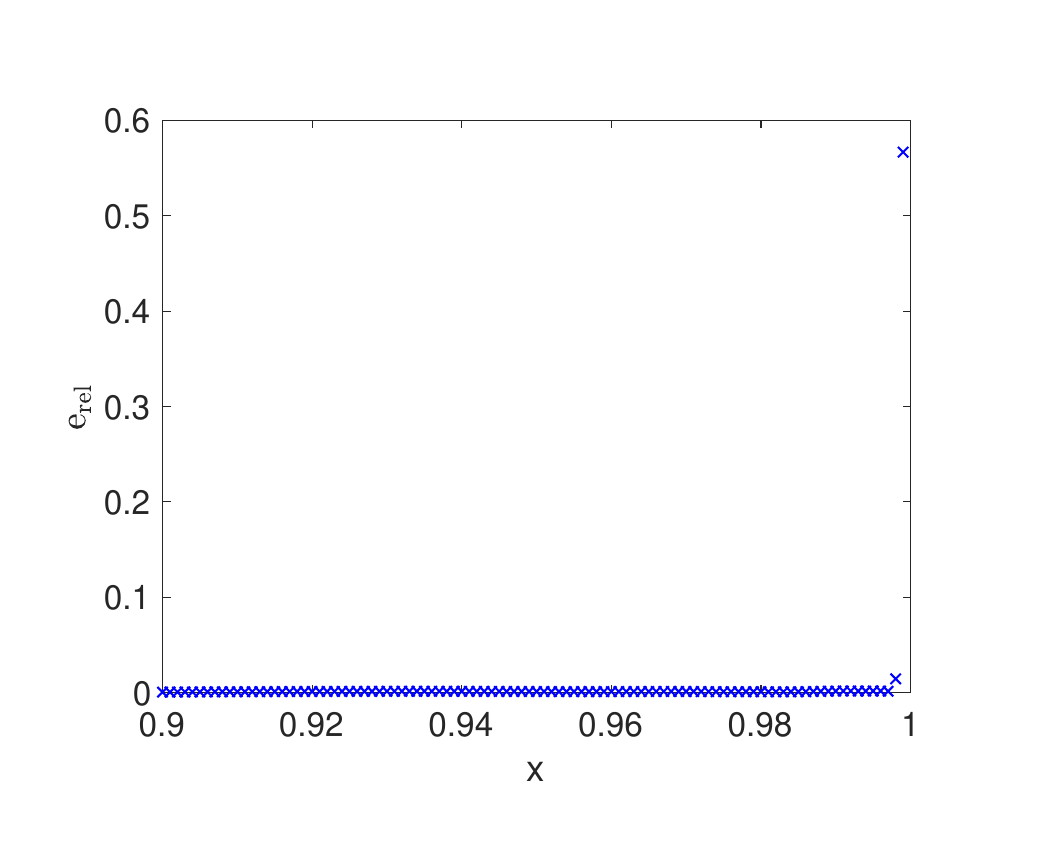}}
\subfigure[relative error away from the boundary layer]
 { 
\includegraphics[width=0.23\textwidth]{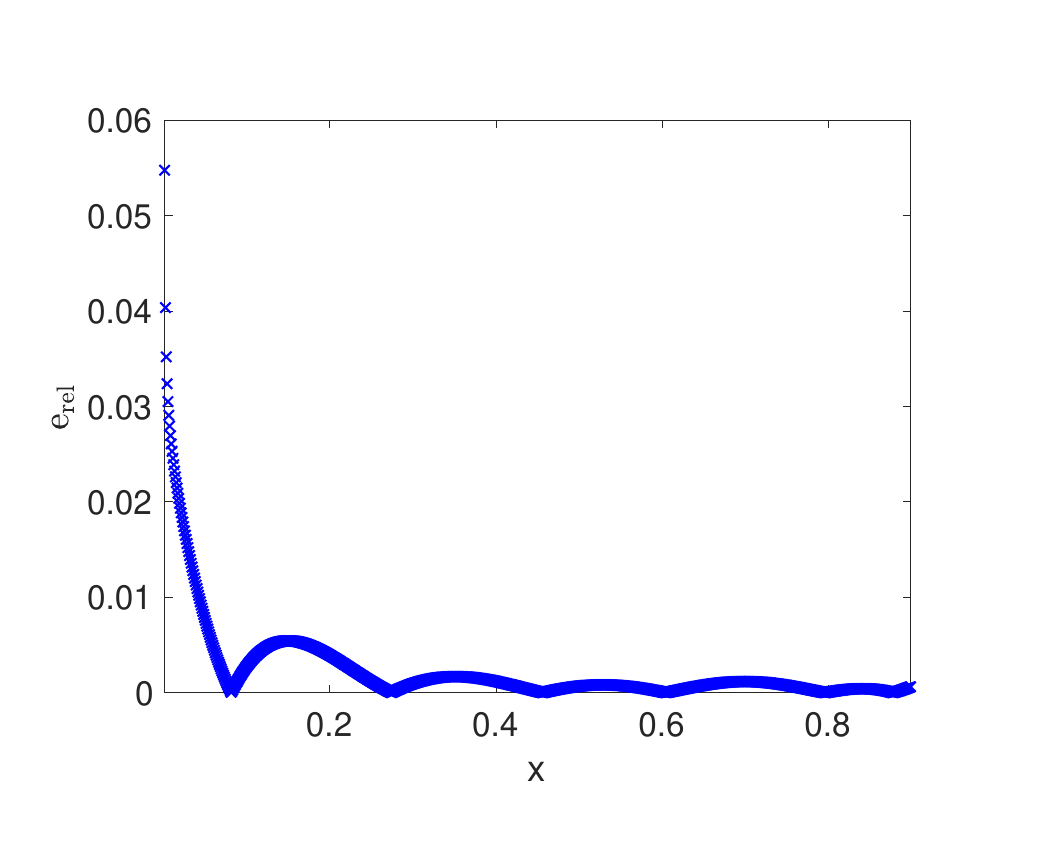}}
 
 \caption{ Numerical results for Example \ref{exm:1p1}  when $\epsilon=10^{-4}$ ((a) to (d)) and $10^{-5} $ ((e) to (h)), with parameters specified in Table \ref{tab:exm1p1_para}.}
 \label{fig:results_exm1p1_1e-4_1e-5} %
 \end{figure}

For the case where $\epsilon=10^{-3}$, we collect the results in Figures \ref{fig:results_exm1p1_1e-3} { and \ref{fig:exm1p1_1p4_loss}(b)}. As depicted in Figure \ref{fig:results_exm1p1_1e-3}(a), the NN solution captures the behavior of the exact solution. While there is a relatively larger error observed very close to $x=1$ as shown in Figure \ref{fig:results_exm1p1_1e-3}(b) and (c), it is essential to note that this challenge is inherent to the nature of the problem itself. In particular, when $\epsilon=10^{-3}$, the boundary layer essentially takes the form of a vertical line. This increases sensitivity to the shift between the actual and predicted boundary layer around $x=1$, as shown in Figure
\ref{fig:results_exm1p1_1e-3}(c). 
{
We also present the results for $\epsilon=10^{-4}$ and $10^{-5}$ in Figure \ref{fig:results_exm1p1_1e-4_1e-5}, where we observe comparable accuracy as in the case $\epsilon=10^{-3}$. }

\begin{table}[ht]
\centering

\setlength{\tabcolsep}{11pt}
\begin{tabular}{c c c c c c}
\hline
$\epsilon$ & $10^{-1}$ (starting $\epsilon_0$)    & $10^{-2}$ & $10^{-3}$ & $10^{-4}$ & $\boldsymbol{\epsilon}_{s}$ \\
\hline
$\alpha$ & 1 & 100 & 1000 & 1000 & 1000 \\
$\alpha_1$ & 0 & 0 & $10^{-6}$ & $10^{-6}$ & $10^{-5}$ \\
$N_c$ & 200 & 200 & 450 & 450 & 300 \\
LR & PC & $10^{-4}$ & $10^{-4}$ & $10^{-4}$ & $10^{-4}$ \\
epochs & 20000 & 30000 & 50000 & 50000 & 30000 each \\
\hline
\end{tabular}
\caption{ Parameters in the loss function \eqref{eq:loss-general} and hyper-parameters of the successive training for Example \ref{exm:1p1}. Here, $\boldsymbol{\epsilon}_s=10^{-5}\times [5, 2.5, 1.25, 1]$, LR is the abbreviation for learning rate, PC is the piecewise constant scheduler in Table \ref{tab:lr_pt}.}
\label{tab:exm1p1_para}
\end{table}



\begin{figure} [!ht]
\centering
\subfigure[exact and NN solutions]
{
\includegraphics[width=0.23\textwidth]{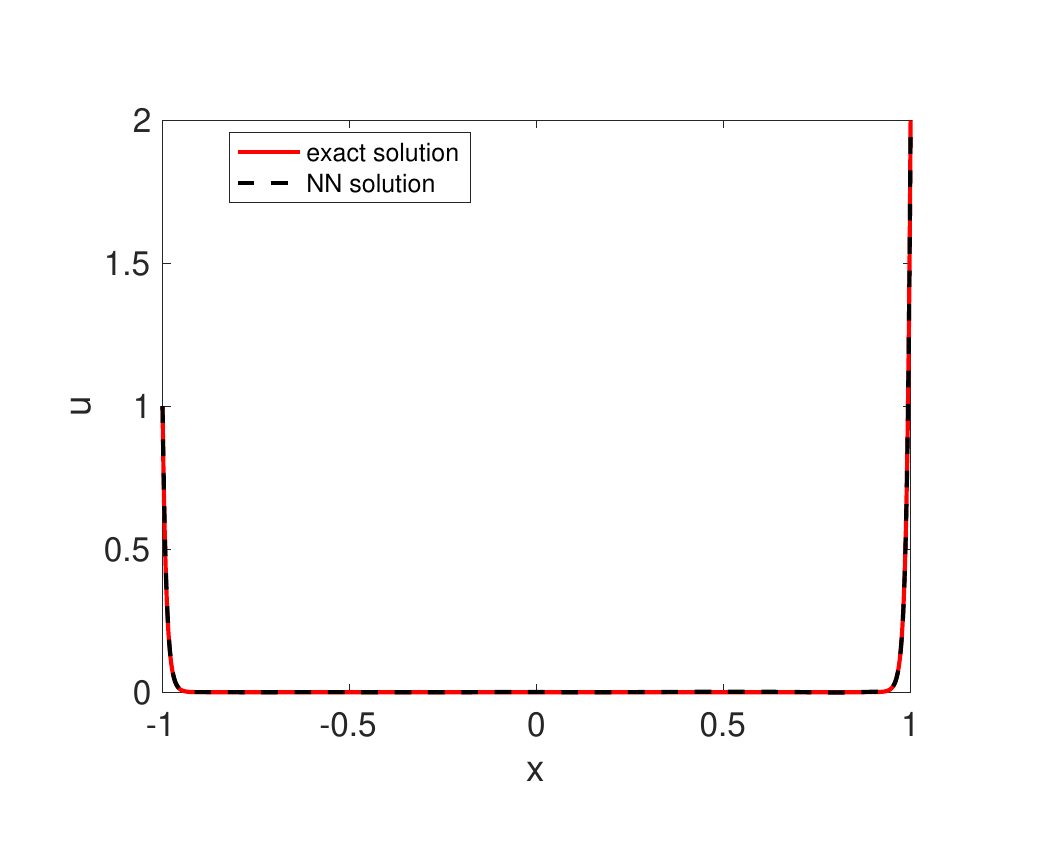}}
\subfigure[absolute error]
 { 
\includegraphics[width=0.23\textwidth]{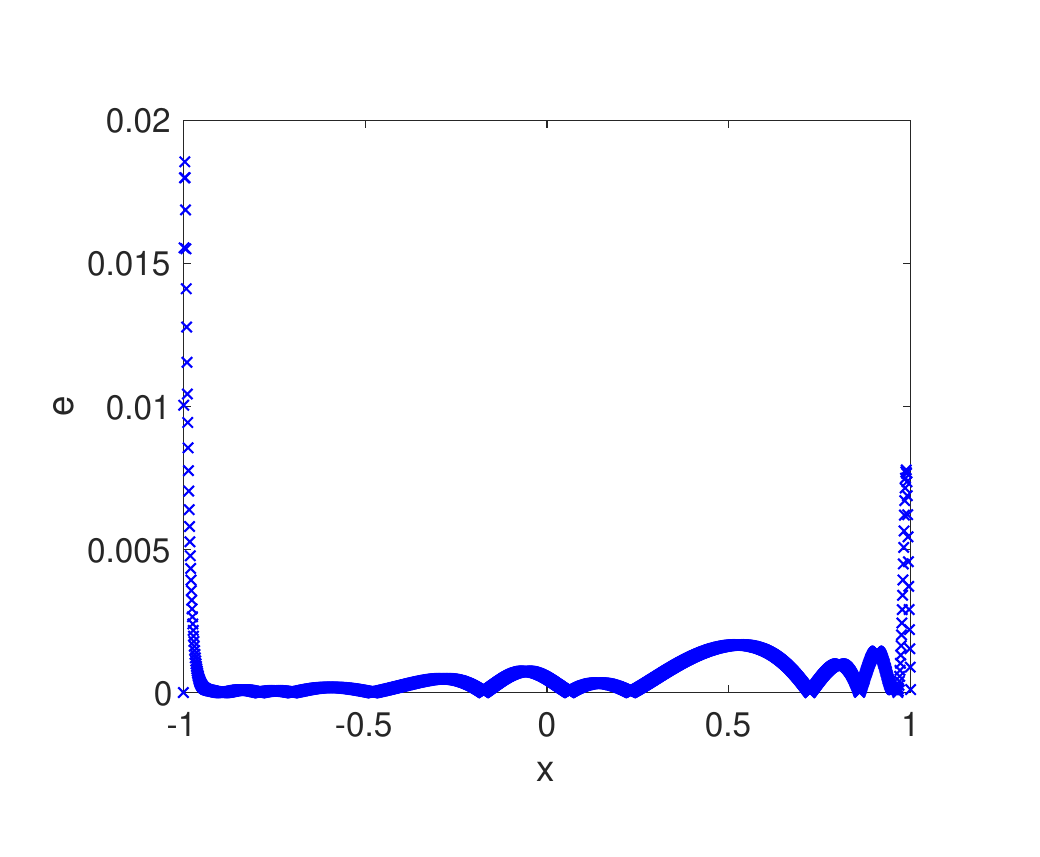}}
 \subfigure[relative error around the left boundary layer]
{
\includegraphics[width=0.23\textwidth]{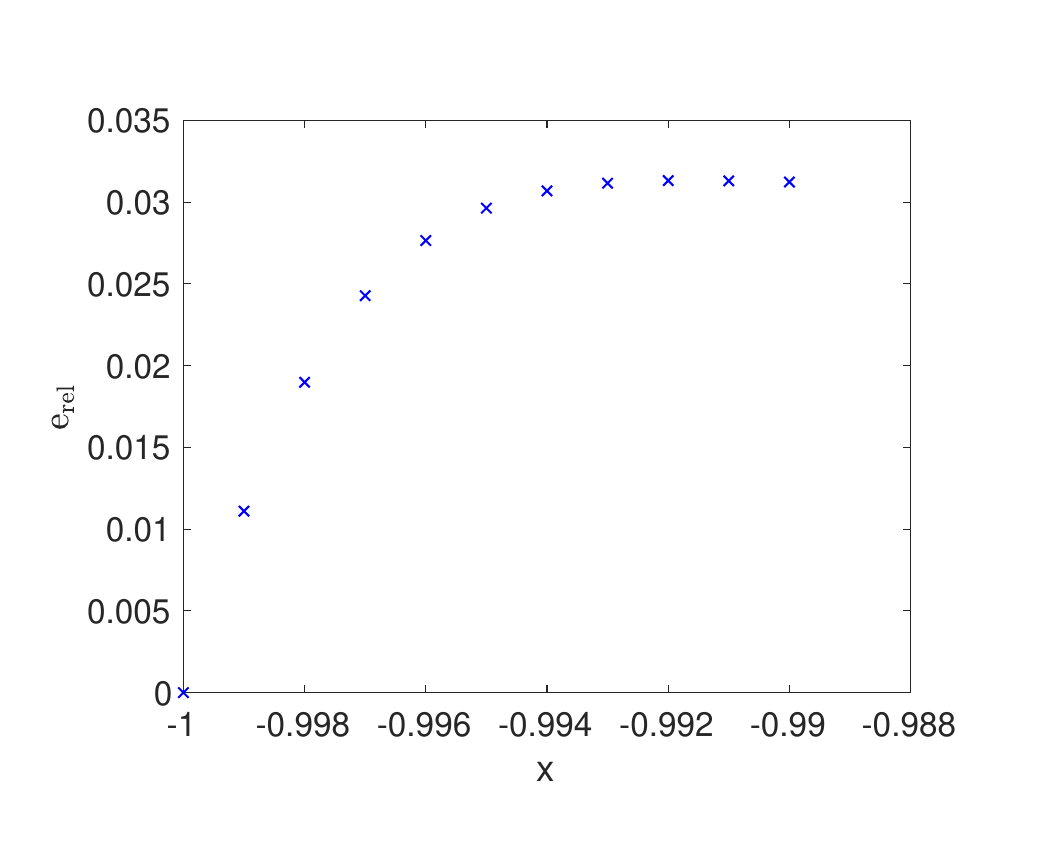}}
\subfigure[relative error around the right boundary layer]
 {  \includegraphics[width=0.23\textwidth]{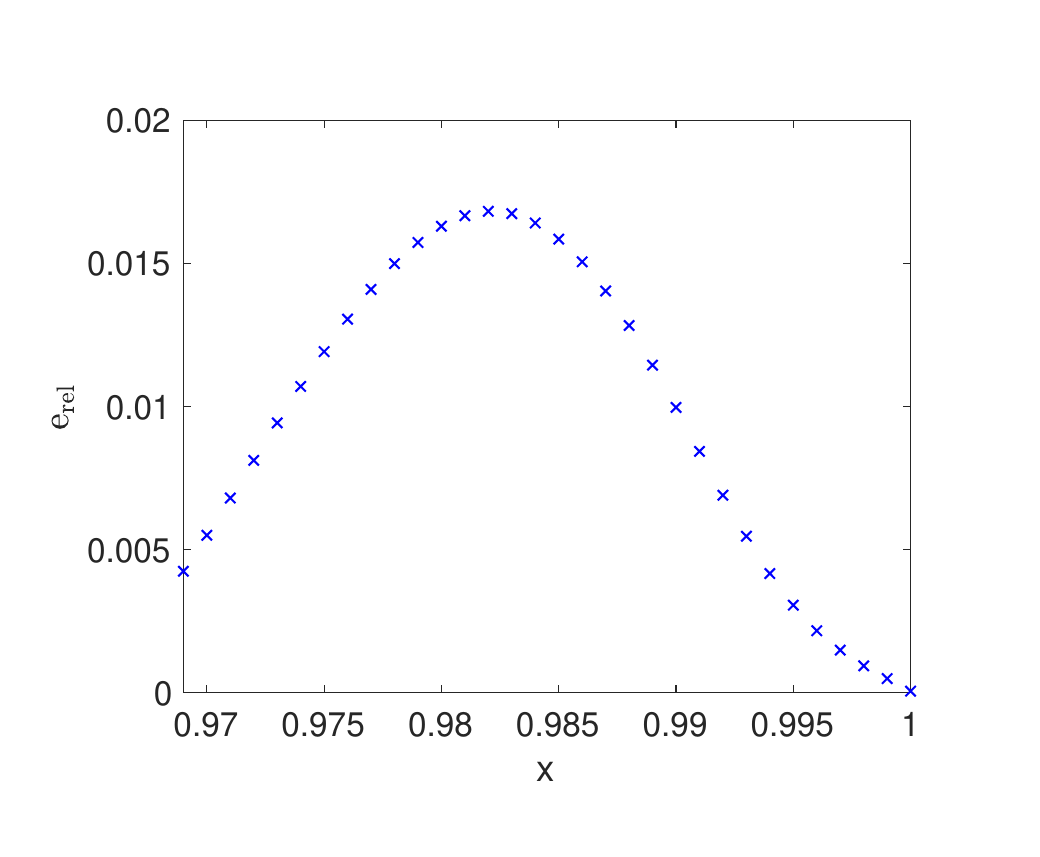}}
 \caption{ Numerical results for Example \ref{exm:1p4}  when $\epsilon=10^{-2}$ using $N(x, x/\sqrt{\epsilon}, 1/\sqrt{\epsilon})$, with parameters specified in Table \ref{tab:exm1p4_para}.}
 

\label{fig:results_exm1p4_1e-2}
 \end{figure}

\begin{figure}[!hb]
	\centering
	     \subfigure[$\epsilon=10^{-2}$ for Example \ref{exm:1p1}]{
        \includegraphics[width=0.26\textwidth]{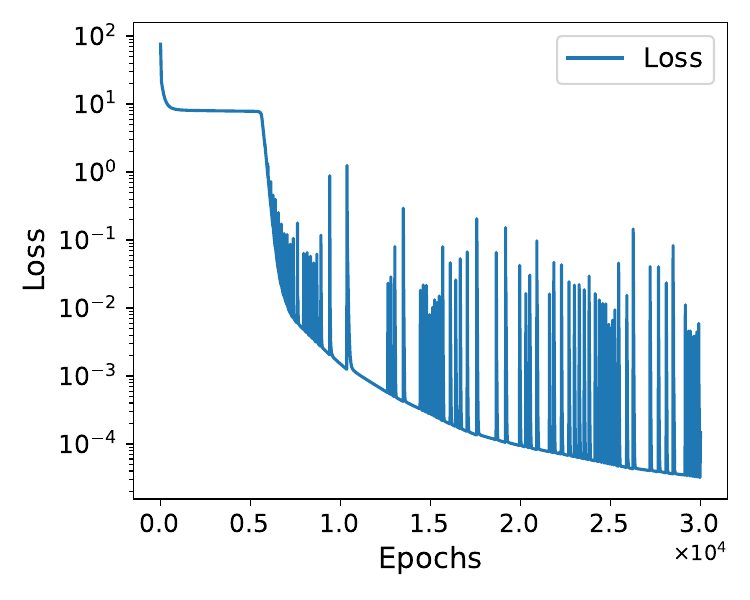}
        \includegraphics[width=0.26\textwidth]{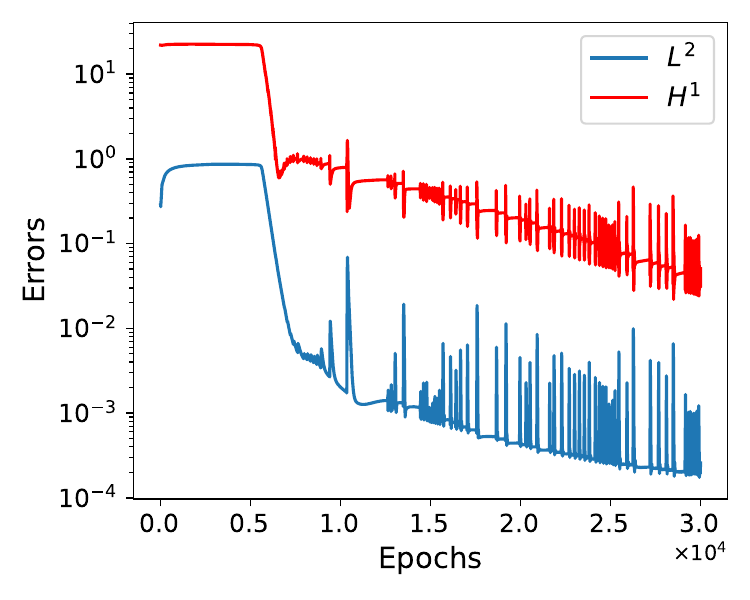}
        \includegraphics[width=0.29\textwidth]{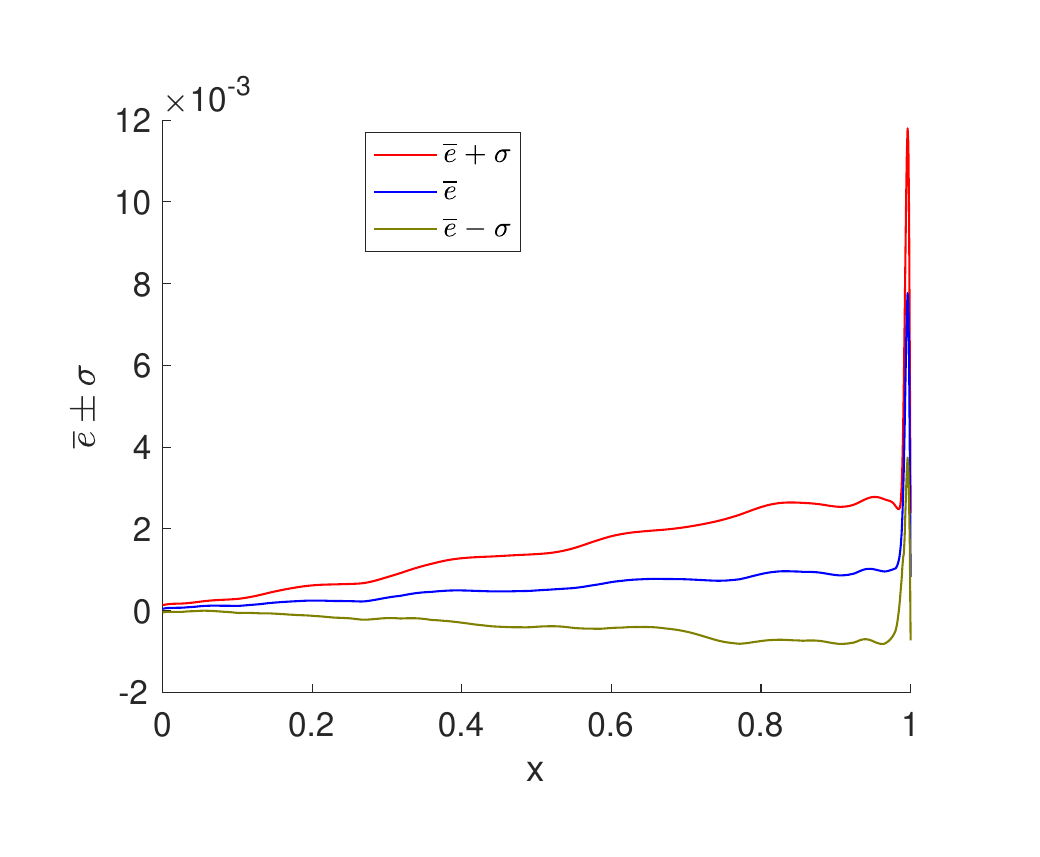}
        }
        
        \subfigure[$\epsilon=10^{-3}$ for Example \ref{exm:1p1}]{
           \includegraphics[width=0.25\textwidth]{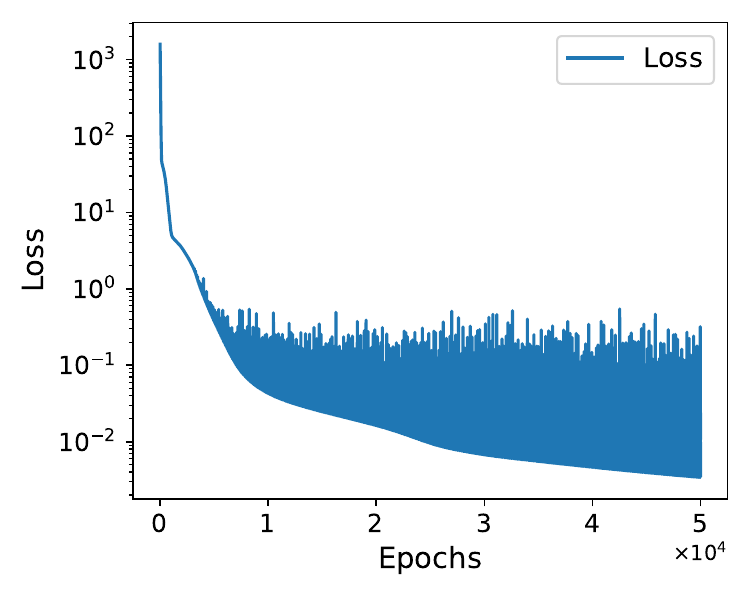}
        \includegraphics[width=0.25\textwidth]{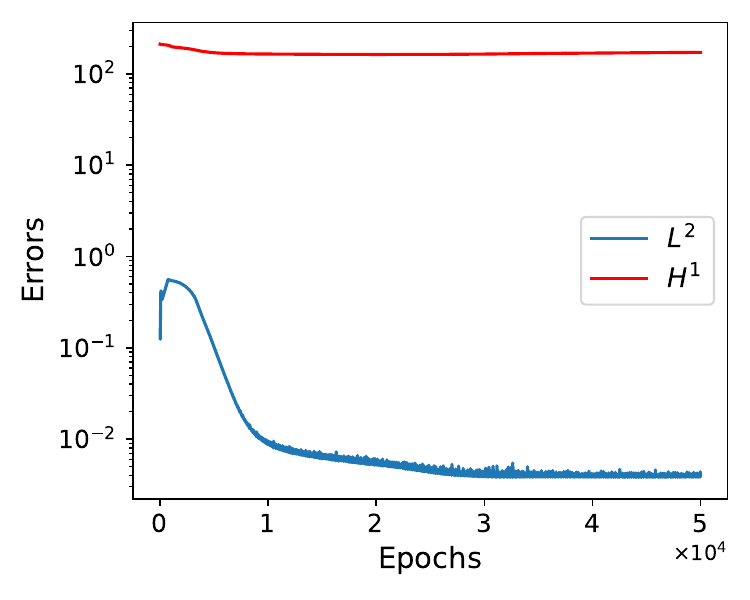}
         \includegraphics[width=0.29\textwidth]{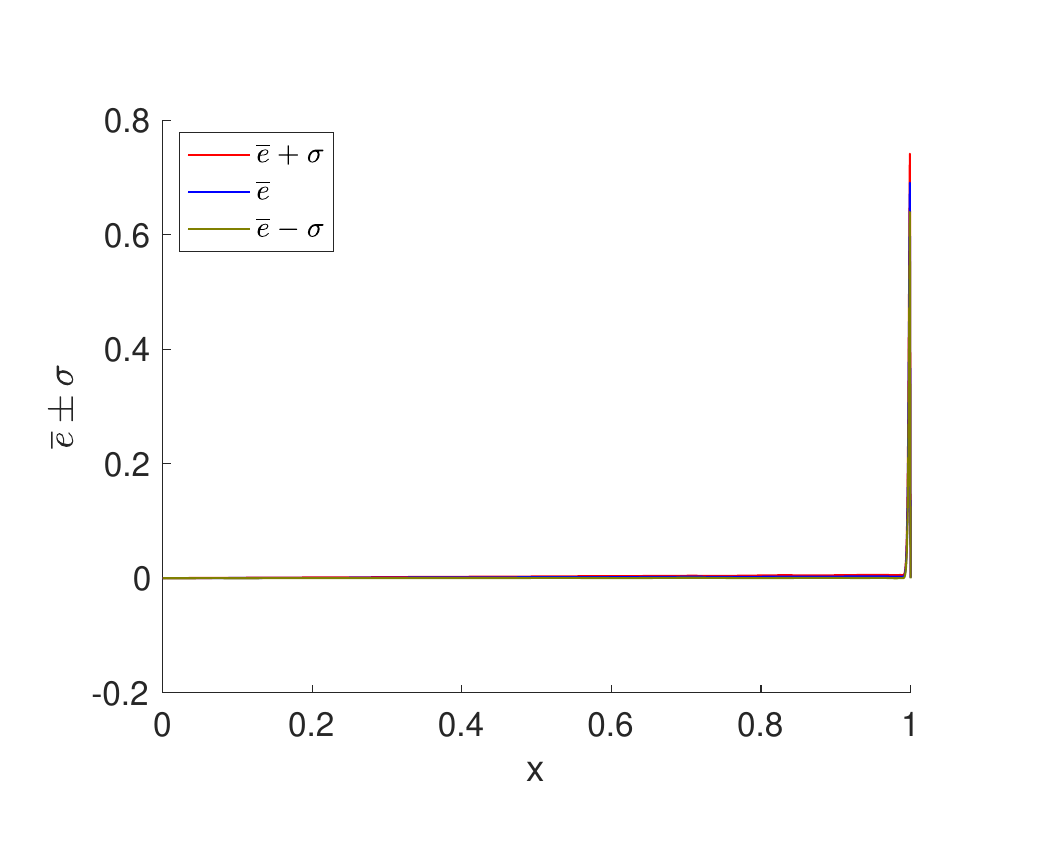}

           }
           
           \subfigure[$\epsilon=10^{-2}$ for Example \ref{exm:1p4}]{
          \includegraphics[width=0.26\textwidth]{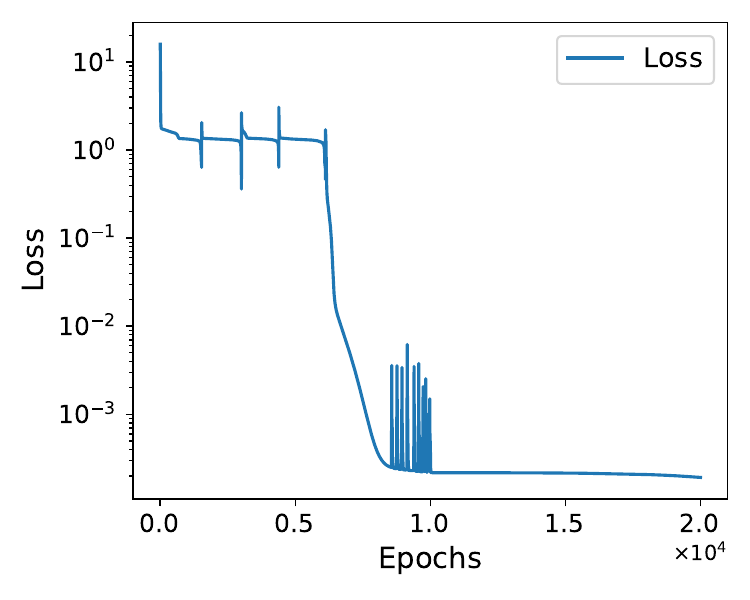}

          \includegraphics[width=0.26\textwidth]{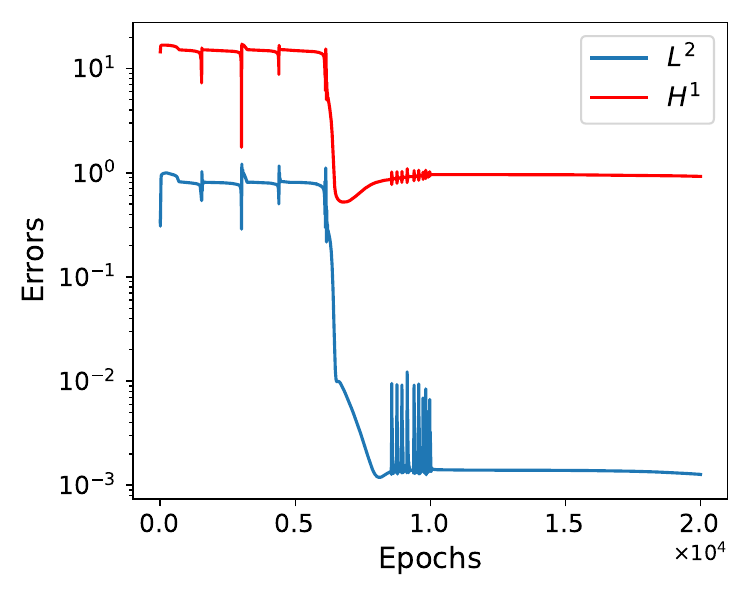}

          \includegraphics[width=0.29\textwidth]{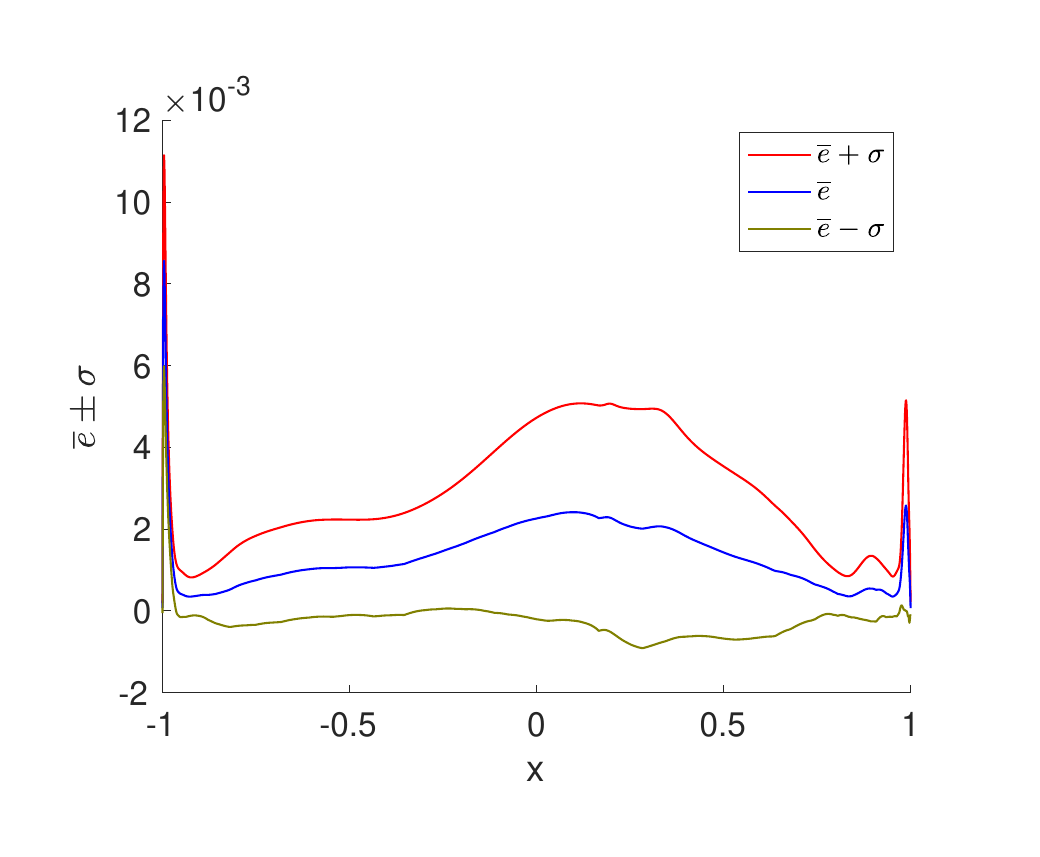}

          }
          \subfigure[$\epsilon=10^{-3}$ for Example \ref{exm:1p4}]{
          \includegraphics[width=0.26\textwidth]{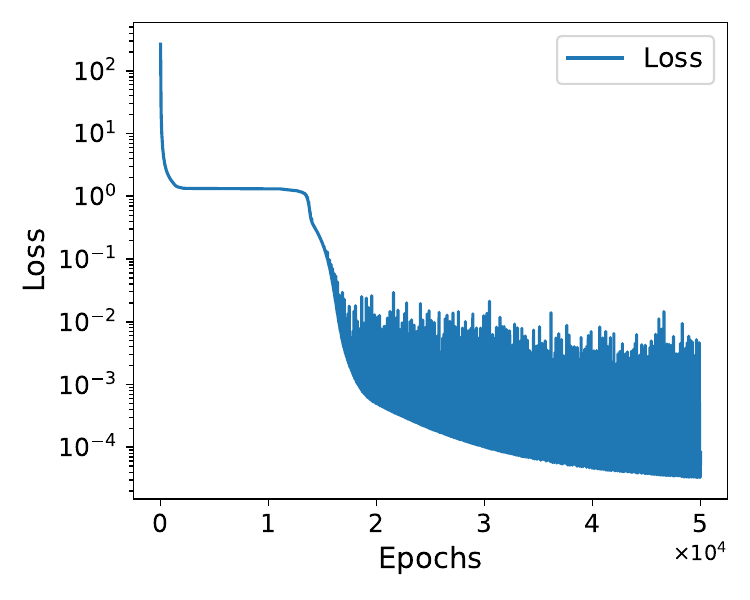}

          \includegraphics[width=0.26\textwidth]{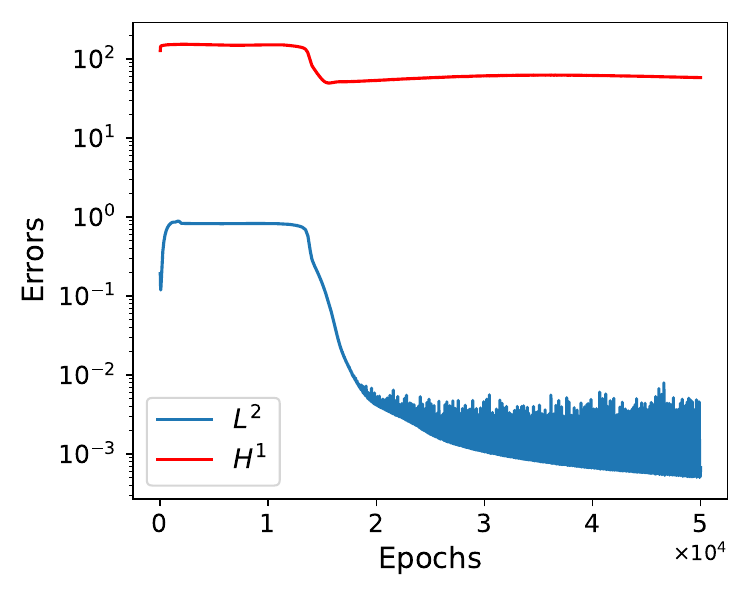}

          \includegraphics[width=0.29\textwidth]{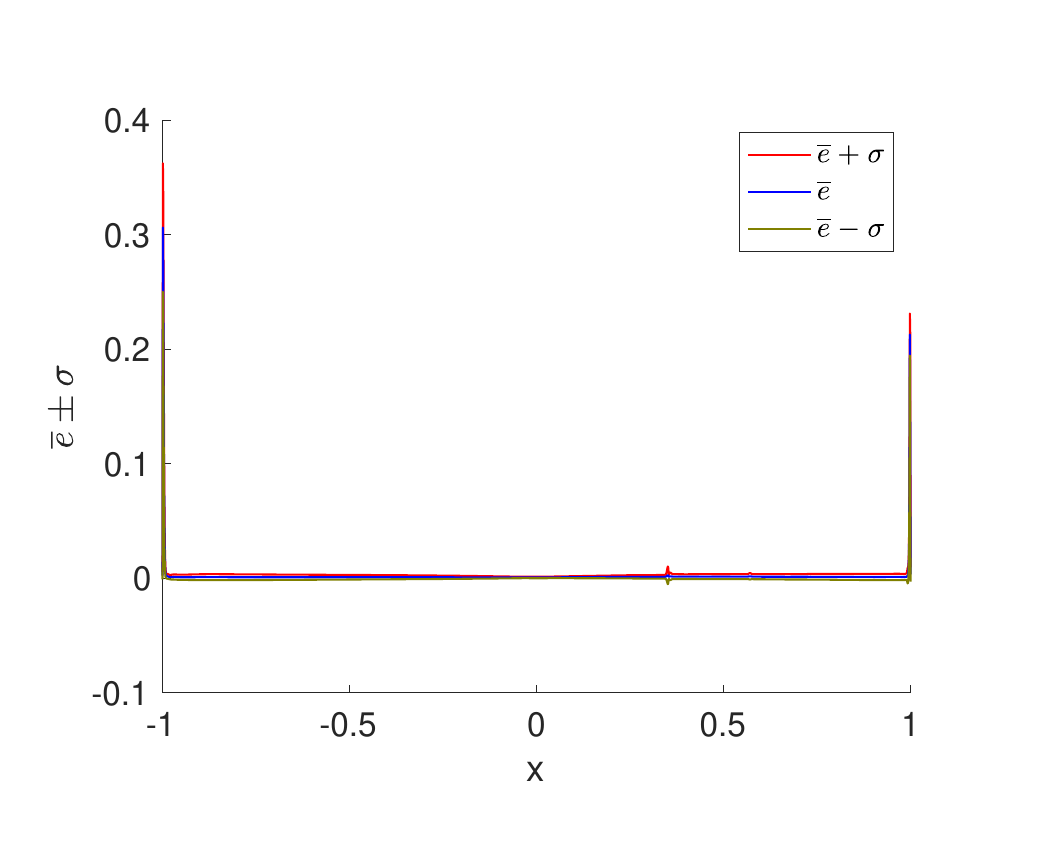}

          }
           
		\caption{ Loss and errors history, as well as the { statistical metrics $\overline{e}$ and $\overline{e}\pm \sigma$} for Examples \ref{exm:1p1} and \ref{exm:1p4} when $\epsilon=10^{-2}$ and $\epsilon=10^{-3}$. 
}
		 \label{fig:exm1p1_1p4_loss}
\end{figure}

\begin{exm}[1D ODE with two boundary layers]\label{exm:1p4}
	\begin{equation}
		\epsilon u'' - x u'-u =0, \quad -1<x<1, \quad u(-1)=1, \quad u(1)=2.
	\end{equation}
\end{exm}
{The solution has two boundary layers, at $x=\pm 1$. 
The exact solution to the problem is} 
\[\exp(\frac{x^2-1}{2\epsilon}) \frac{\text{Erf} (\frac{x}{\sqrt{2\epsilon}})+ 3 \text{Erf} (\frac{1}{\sqrt{2\epsilon}})}{2\text{Erf} (\frac{1}{\sqrt{2\epsilon}})},\quad 
\text{Erf} (z)=\frac{2}{\sqrt{\pi}}\int_0^z \exp(-t^2)\,dt.\]
We solve the considered problem using the two-scale NN $N(x, x/\sqrt{\epsilon}, 1/\sqrt{\epsilon})$ {. 
We employ the successive training strategy in Algorithm \ref{alg:succesive-training-two-scale} starting with $\epsilon_0=10^{-1}$, and successively solve the problem for $\epsilon$ as small as $10^{-5}$, following the sequence $\epsilon_0=10^{-1},10^{-2},10^{-3},10^{-4},
  5\times10^{-5},2.5\times 10^{-5},1.25\times 10^{-5}, 10^{-5}$.
The parameters in the discrete loss function in \eqref{eq:loss-general} and hyper-parameters of the training are summarized in Table \ref{tab:exm1p4_para}}.
The NN size is $(3,20,20,20,20,1)$. We collect the numerical results in Figures \ref{fig:results_exm1p4_1e-2} { and \ref{fig:exm1p1_1p4_loss}(c)} for the case when $\epsilon=10^{-2}$. The figures illustrate the excellent agreement between the NN solution and the exact solution, even at the boundary layers. 
As shown in Figure \ref{fig:results_exm1p4_1e-2} (c) and (d), the relative error is at the level of $10^{-2}$ in the vicinity of the { left and right boundary layers}. These results underscore the capability of the two-scale NN method in capturing multiple boundary layers in the exact solution. 

\begin{figure} [!ht]
\centering
\subfigure[exact and NN solutions]
{
\includegraphics[width=0.22\textwidth]{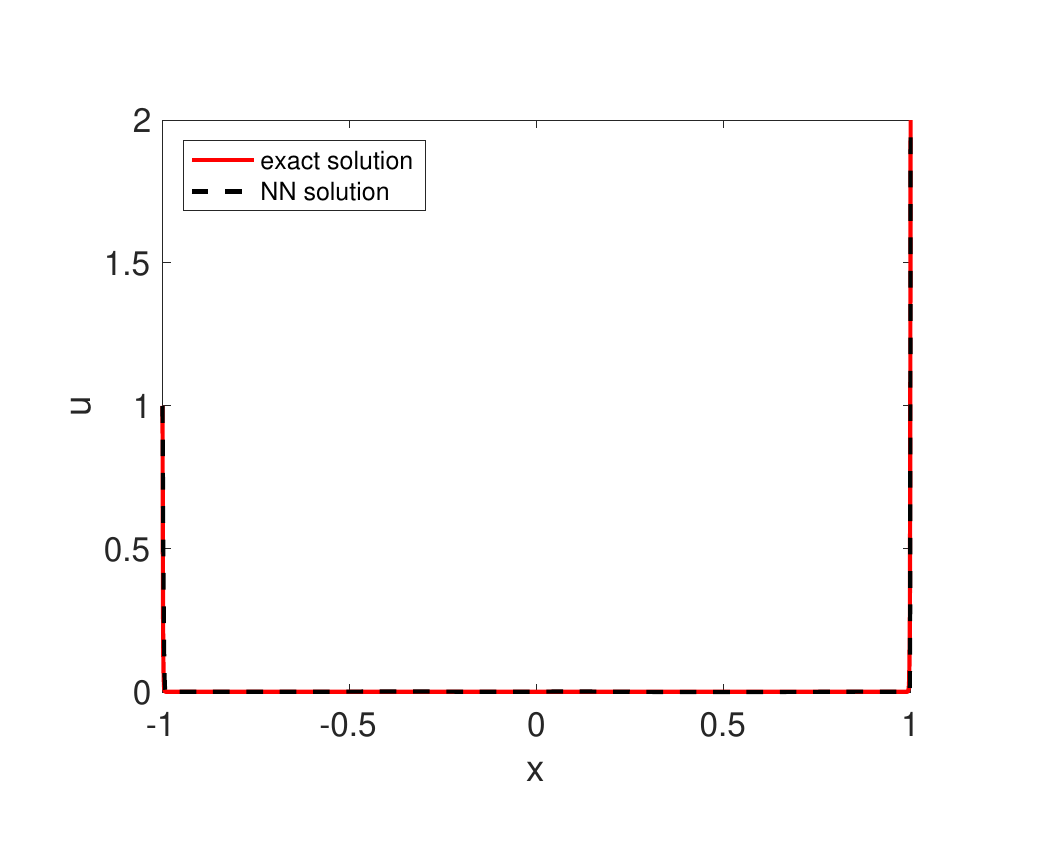}}
\subfigure[absolute error]
 { 
\includegraphics[width=0.22\textwidth]{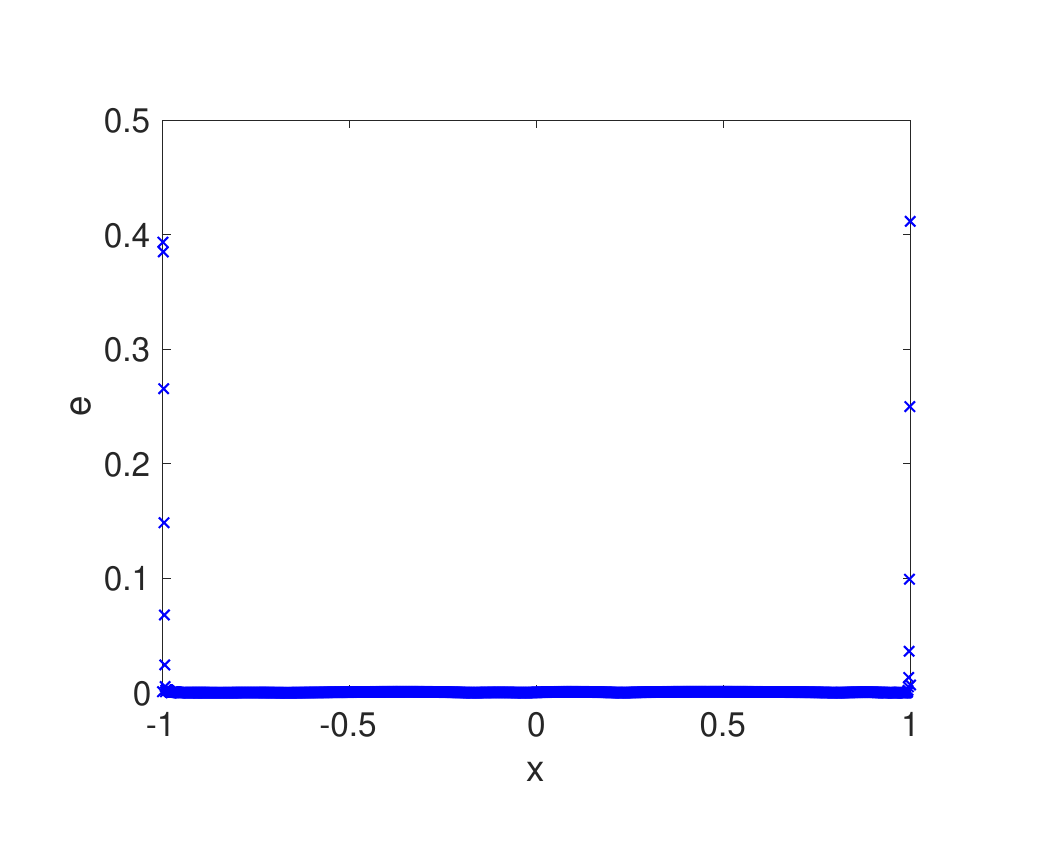}}
 \subfigure[exact and NN solutions around the right boundary]
{
\includegraphics[width=0.22\textwidth]{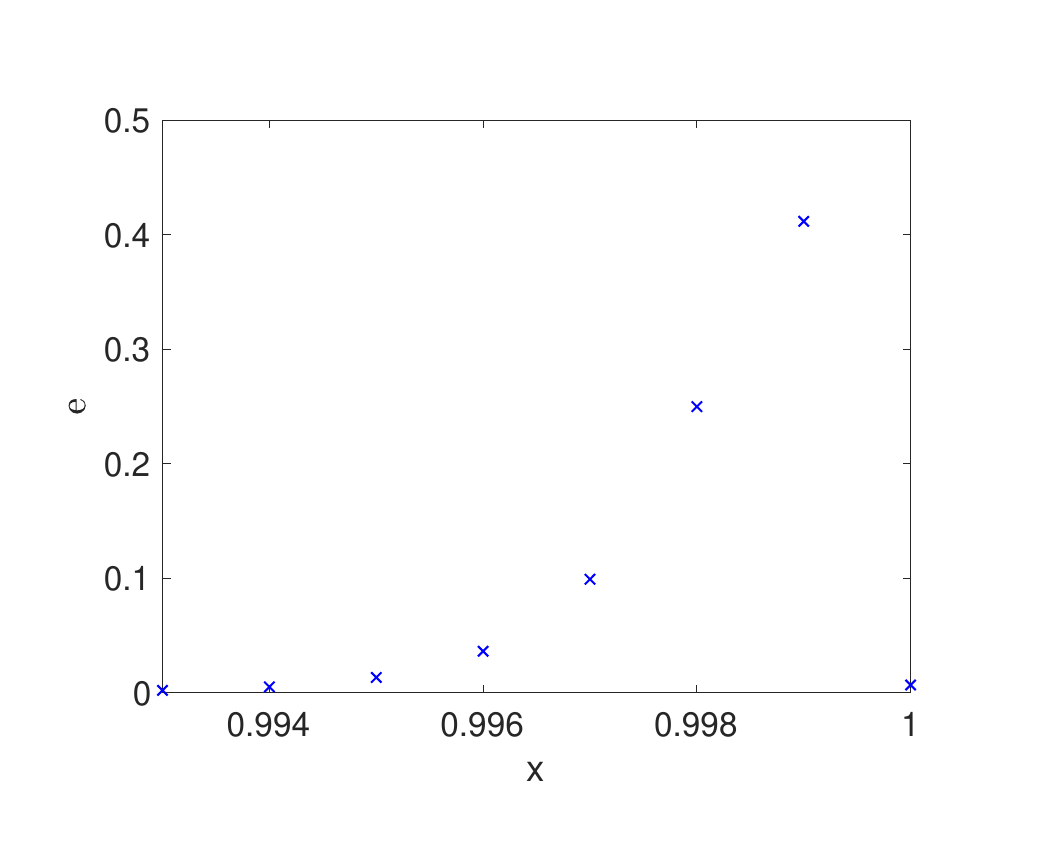}}
\subfigure[absolute error around the right boundary]
 { 
 \includegraphics[width=0.22\textwidth]{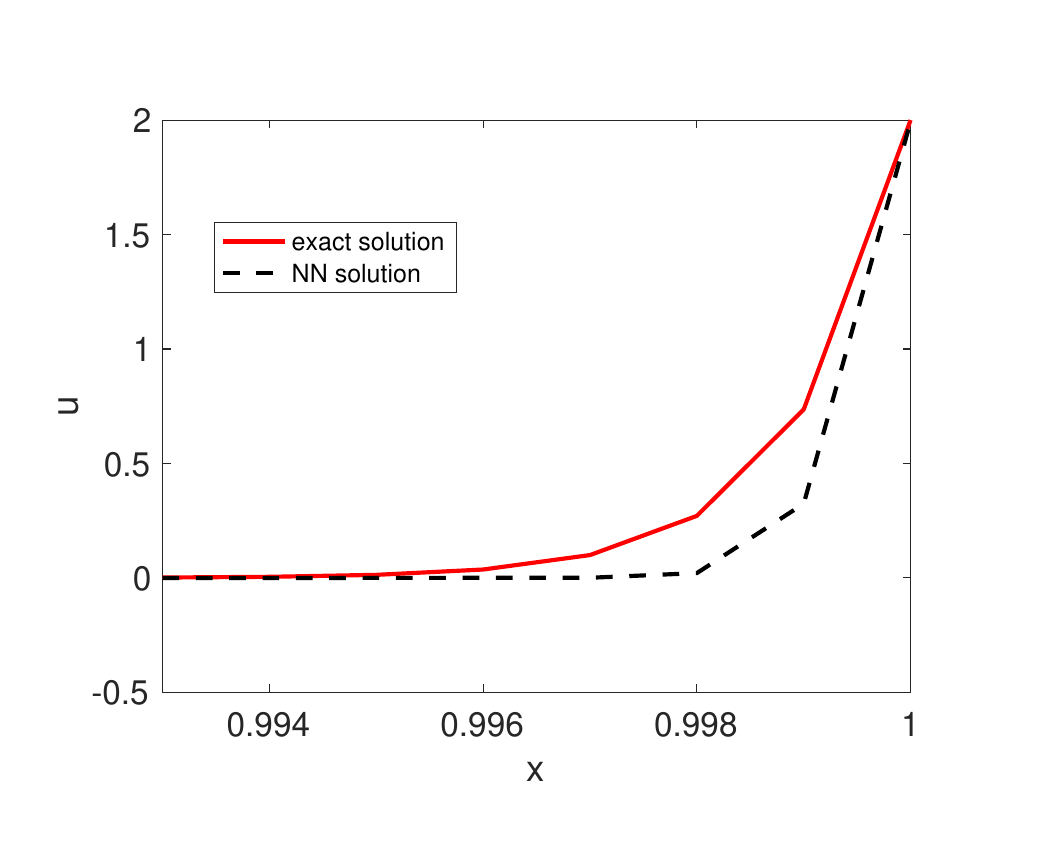}}
 
 \caption{  Numerical results for Example \ref{exm:1p4} when $\epsilon=10^{-3}$ using $N(x, x/\sqrt{\epsilon}, 1/\sqrt{\epsilon})$, 
 with parameters specified in Table \ref{tab:exm1p4_para}.}
 \label{fig:results_exm1p4_1e-3} 
 \end{figure}

\begin{figure} [!ht]
\centering
\subfigure[exact and NN solutions]
{
\includegraphics[width=0.22\textwidth]{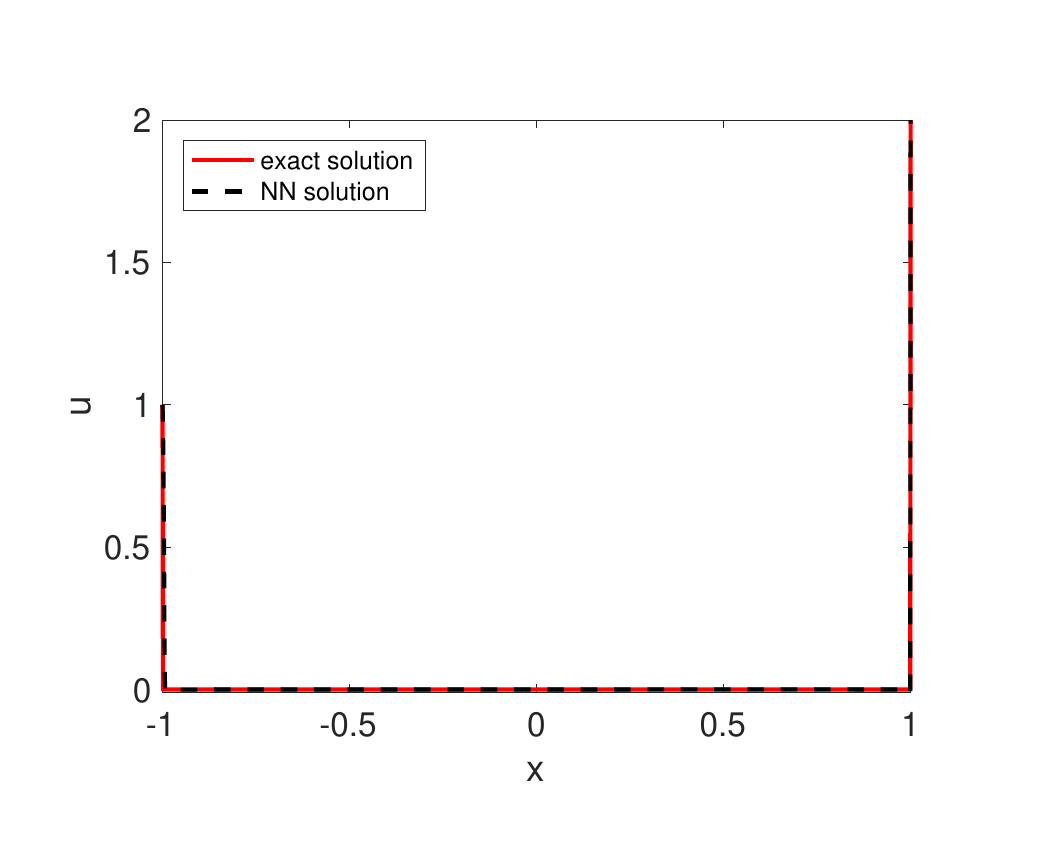}}
\subfigure[absolute error]
 { 
\includegraphics[width=0.22\textwidth]{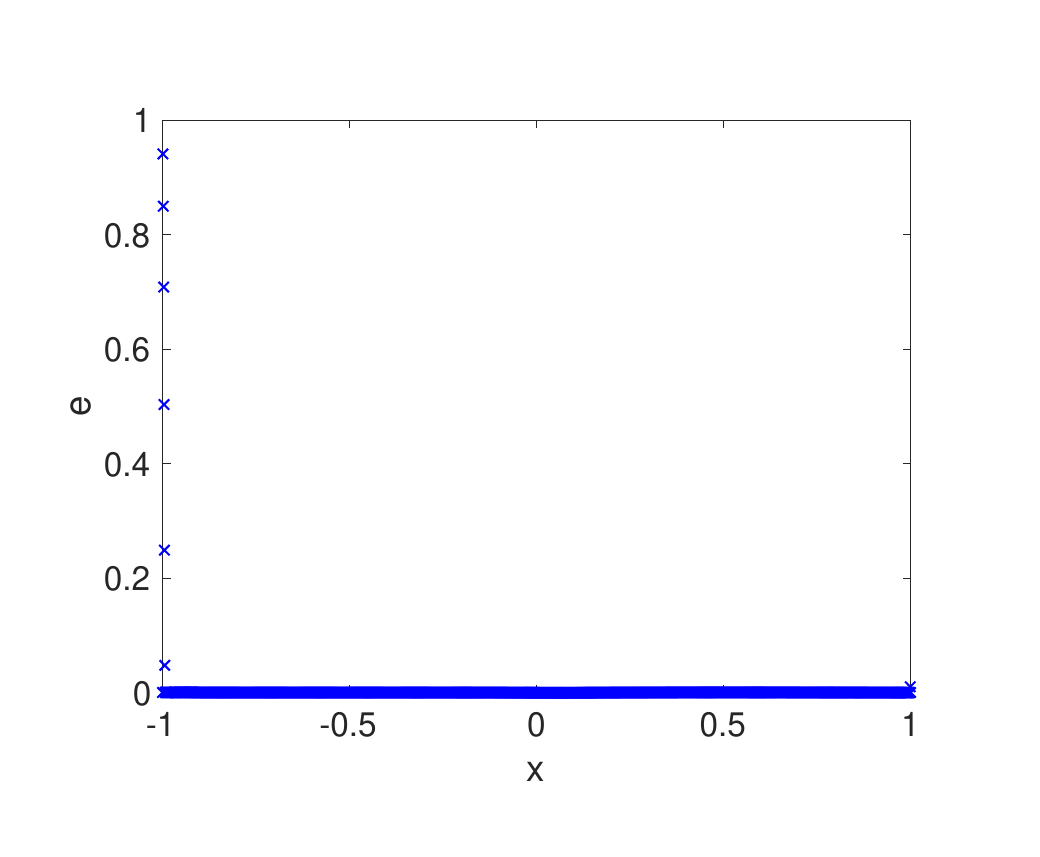}}
 \subfigure[exact and NN solutions around the left boundary]
{
\includegraphics[width=0.22\textwidth]{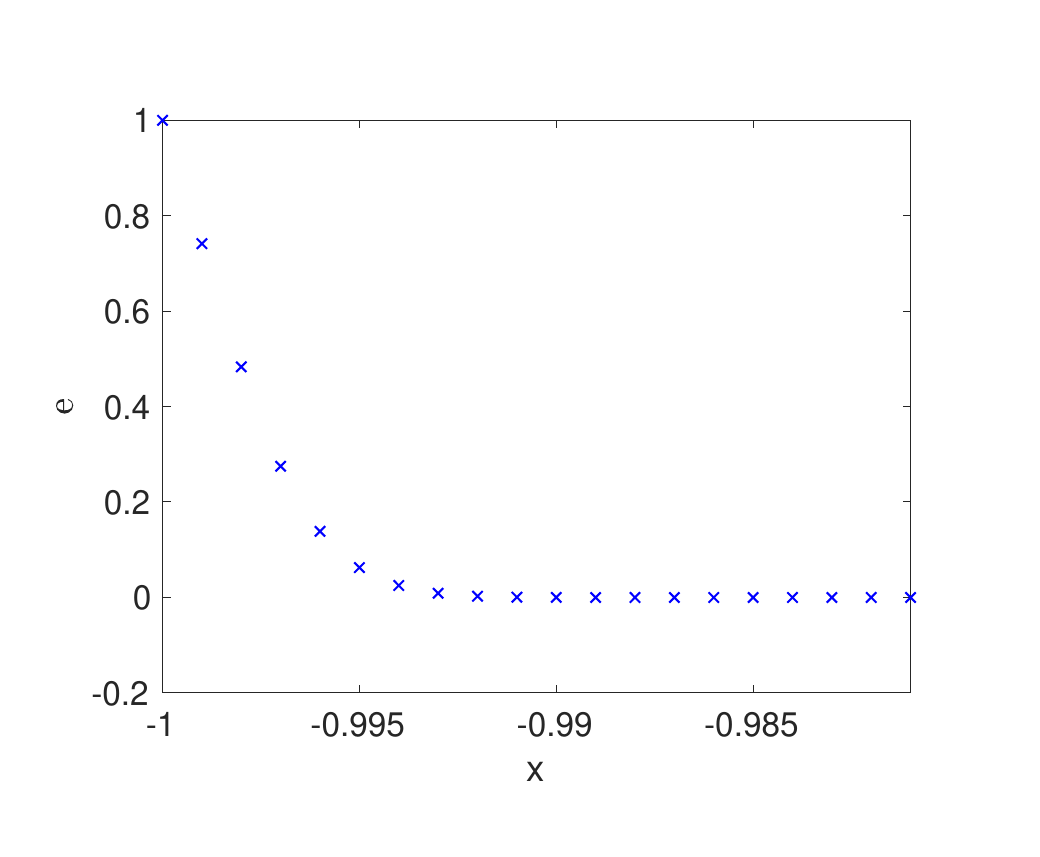}}
\subfigure[absolute error around the left boundary]
 { 
 \includegraphics[width=0.22\textwidth]{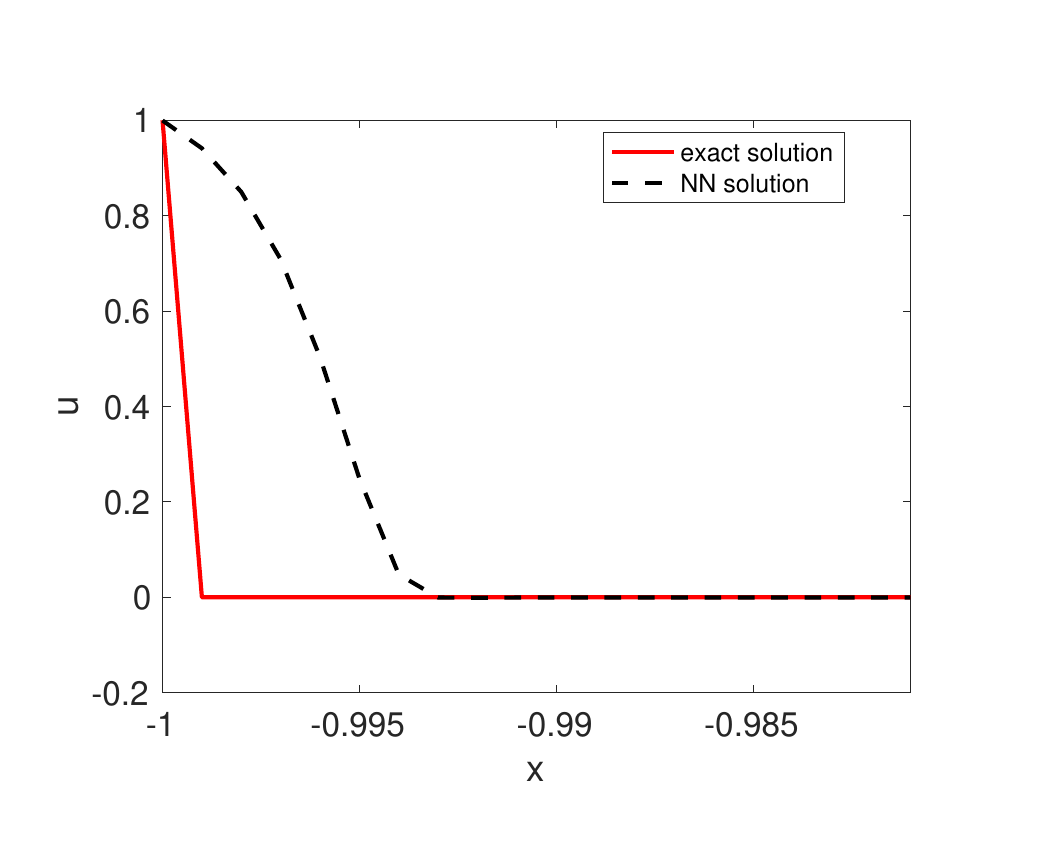}}
 
 \caption{ Numerical results for Example \ref{exm:1p4} when $\epsilon=10^{-5}$ using $N(x, x/\sqrt{\epsilon}, 1/\sqrt{\epsilon})$, with parameters specified in Table \ref{tab:exm1p4_para} .} 
 
 %
 
 \label{fig:results_exm1p4_1e-5} 
 \end{figure}

For the scenario when \( \epsilon = 10^{-3} \), the boundary layers practically become vertical lines, as shown in Figure \ref{fig:results_exm1p4_1e-3}(a). 
We collect the results in Figures \ref{fig:results_exm1p4_1e-3} { and \ref{fig:exm1p1_1p4_loss}(d)}. 
While the NN solution does not achieve the same level of precision as when $\epsilon=10^{-2}$, it is crucial to note that this discrepancy is primarily due to the inherent difficulty of the problem itself: the two boundary layers essentially become vertical lines. Despite the challenge, the two-scale NN solution still captures the key features of the exact solution, as shown in Figure \ref{fig:results_exm1p4_1e-3}(a).  In addition, we provide a visualization of the shift between the actual and predicted solutions in Figure \ref{fig:results_exm1p4_1e-3}(c) and (d), illustrating that the shift occurs within a very narrow band ({ width$\approx 4\times 10^{-3}$}) around the right boundary. { We also present the results for $\epsilon=10^{-5}$ in Figure \ref{fig:results_exm1p4_1e-5}, where the two-scale NN solution provides comparable accuracy as in the case $\epsilon=10^{-3}$. 
We closely examine the left boundary, where the largest errors occur, and observe in Figure \ref{fig:results_exm1p4_1e-5} (d) that the shift between the exact and NN solutions is confined to a very narrow band (width $\approx 10^{-2}$). 

}

\begin{table}[ht]
\centering

\setlength{\tabcolsep}{11pt}
\begin{tabular}{c c c c c c}
\hline
$\epsilon$ & $10^{-1}$ (starting $\epsilon_0$)  & $10^{-2}$ & $10^{-3}$ & $10^{-4}$ & $\boldsymbol{\epsilon}_{s}$ \\
\hline
$\alpha$ & 1 & 1 & 1 & 1 & 1 \\
$\alpha_1$ & 0 & 0 & 0 & $10^{-4}$ & $10^{-4}$ \\
$N_c$ & 500 & 500 & 430 & 430 & 430 \\
LR & PC & PC & $10^{-4}$ & $10^{-4}$ & $10^{-4}$ \\
epochs & 20000 & 20000 & 50000 & 50000 & 20000 each \\
\hline
\end{tabular}
\caption{ Parameters in the loss function \eqref{eq:loss-general} and hyper-parameters of the successive training for Example \ref{exm:1p4}. Here, $\boldsymbol{\epsilon}_s=10^{-5}\times [5, 2.5, 1.25, 1]$, LR is the abbreviation for learning rate, PC is the piecewise constant scheduler in Table \ref{tab:lr_pt}.}
\label{tab:exm1p4_para}
\end{table}






\begin{exm} [1D viscous Burgers' equation]\label{exm:burgers}
\begin{equation}
\partial_t u + u u_x - \epsilon u_{xx} =0
\end{equation}
with initial conditions: $u(0,x)=-\sin(\pi (x-x_0))$ and boundary conditions  $u(t,-1)=u(t,1)=g(t)$, where $g(t)$ is determined from the exact solution
\begin{equation}\label{eq:exact_burgers}
u(x,t)=-2\epsilon \frac{\partial_x v}{v},\quad \text{ where }
v= \int_{\Real}  \exp (-\frac{\cos \left(\pi(x-x_0-2\sqrt{\epsilon t} s)\right)}{2\pi\epsilon}) e^{-s^2} ds.
\end{equation}
\end{exm}

{}

\begin{figure} [!ht]
\centering
\subfigure[exact and NN solutions when $x_0=0$]
{
\includegraphics[width=0.32\textwidth]{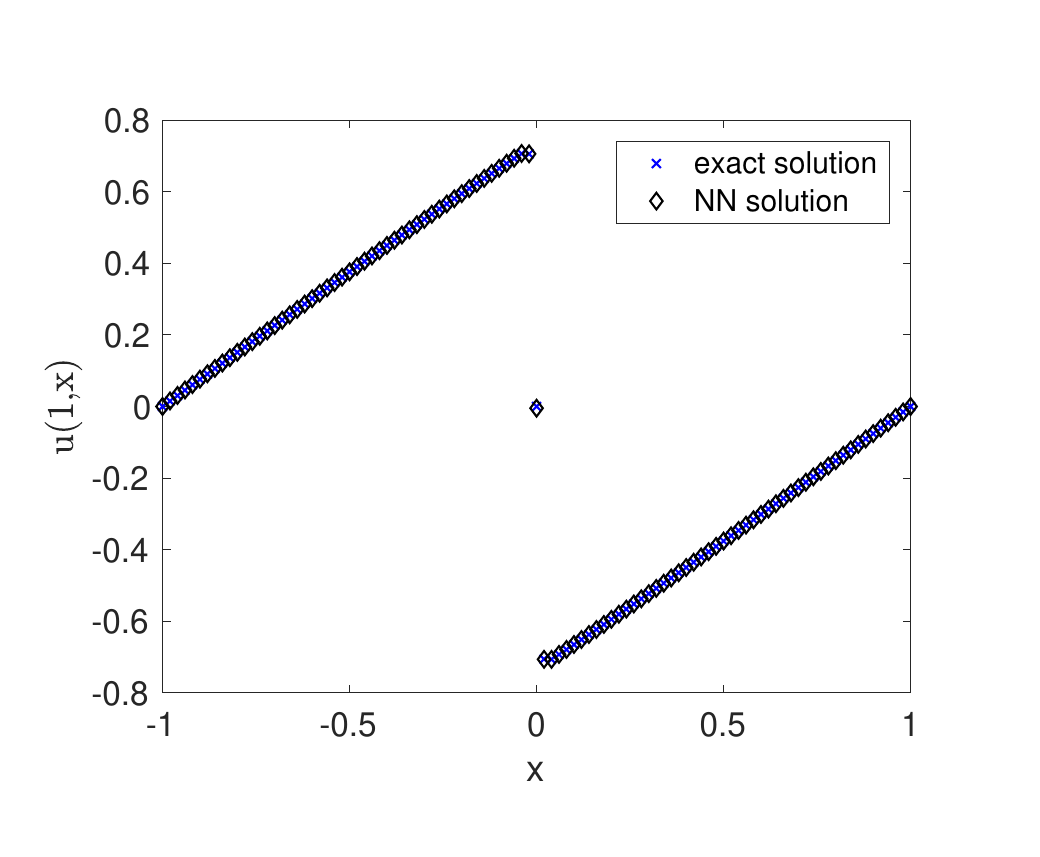}}
\subfigure[exact and NN solutions when \newline
$x_0=0.5$]
{ 
\includegraphics[width=0.32\textwidth]{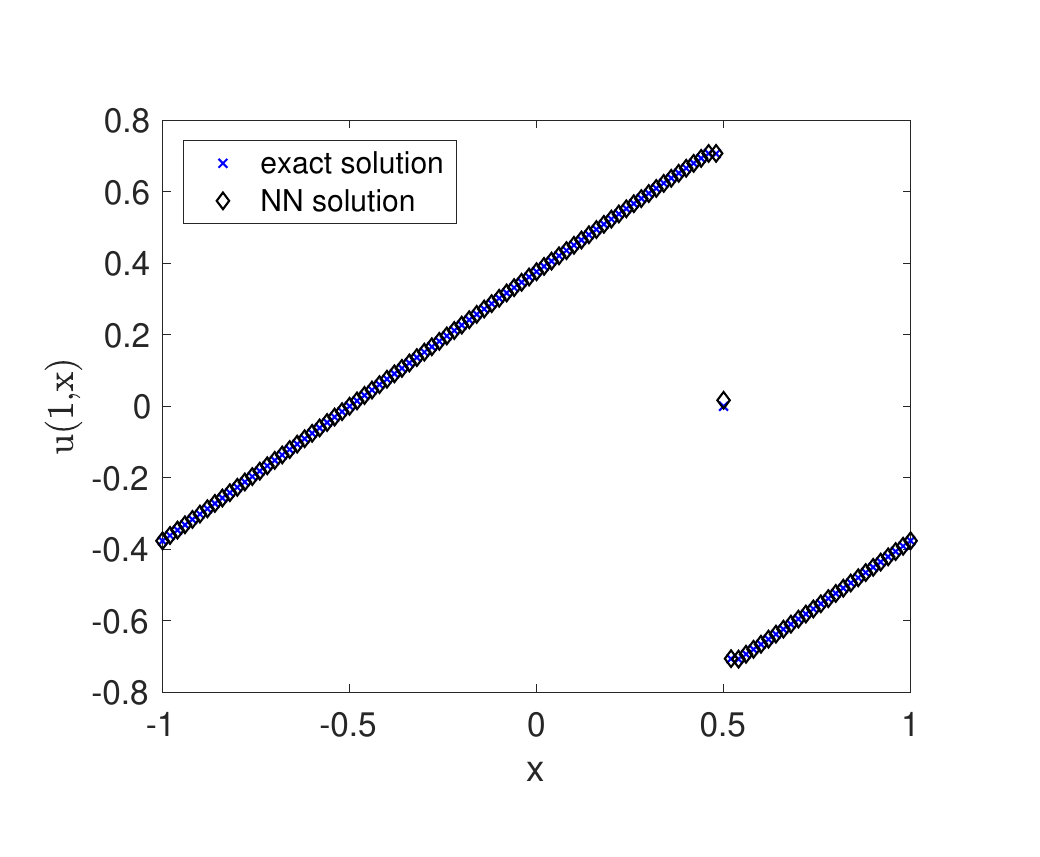}}
\subfigure[{ loss history when $x_0=0$}]
 { 
\includegraphics[width=0.29\textwidth]{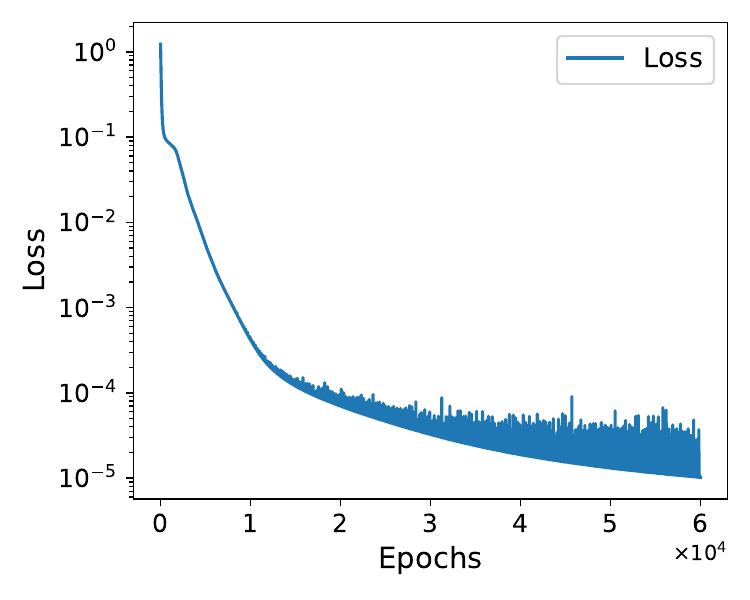}}

%
 \subfigure[absolute error when $x_0=0$]
{
\includegraphics[width=0.32\textwidth]{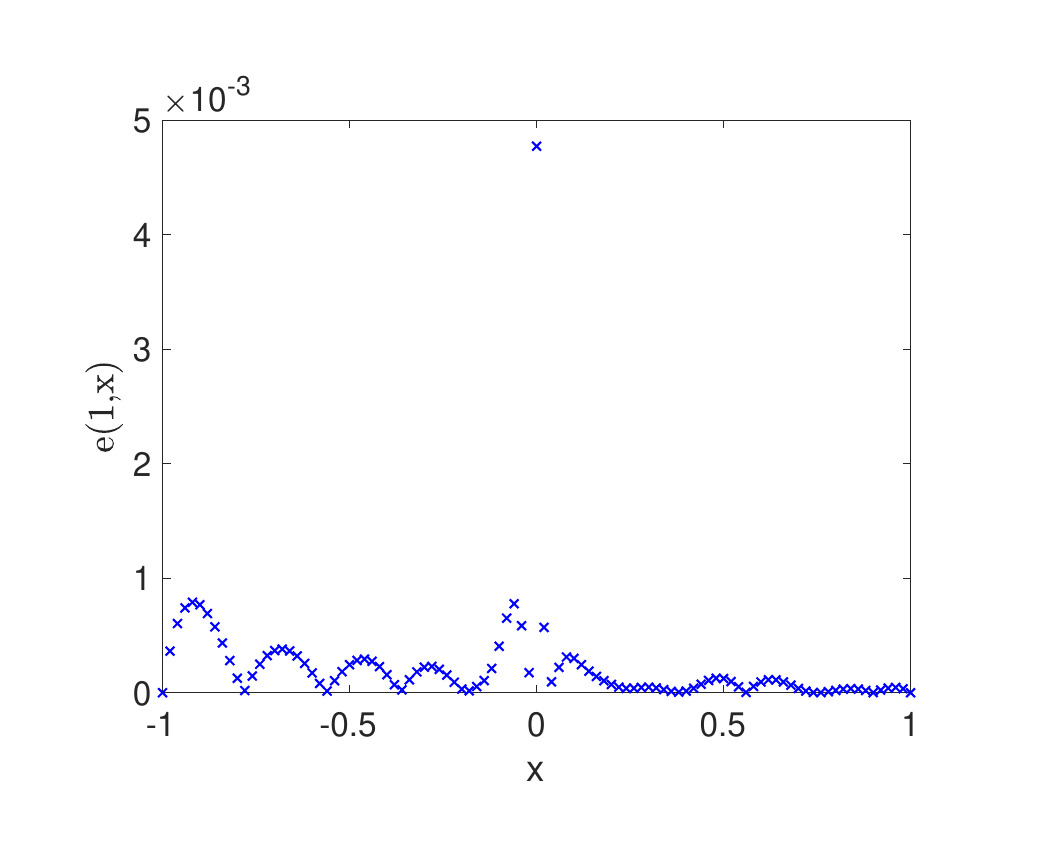}}
\subfigure[relative error when $x_0=0$]
 {  \includegraphics[width=0.32\textwidth]{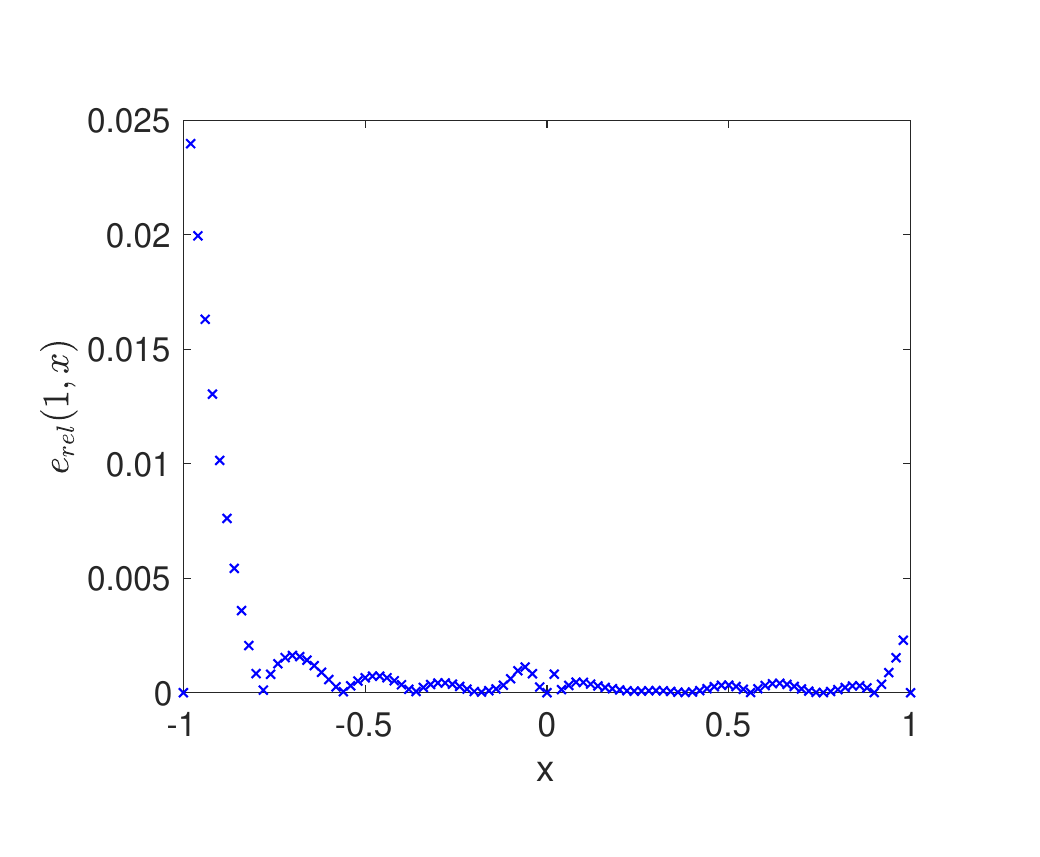}}
  \subfigure[absolute error when $x_0=0.5$]
{
\includegraphics[width=0.32\textwidth]{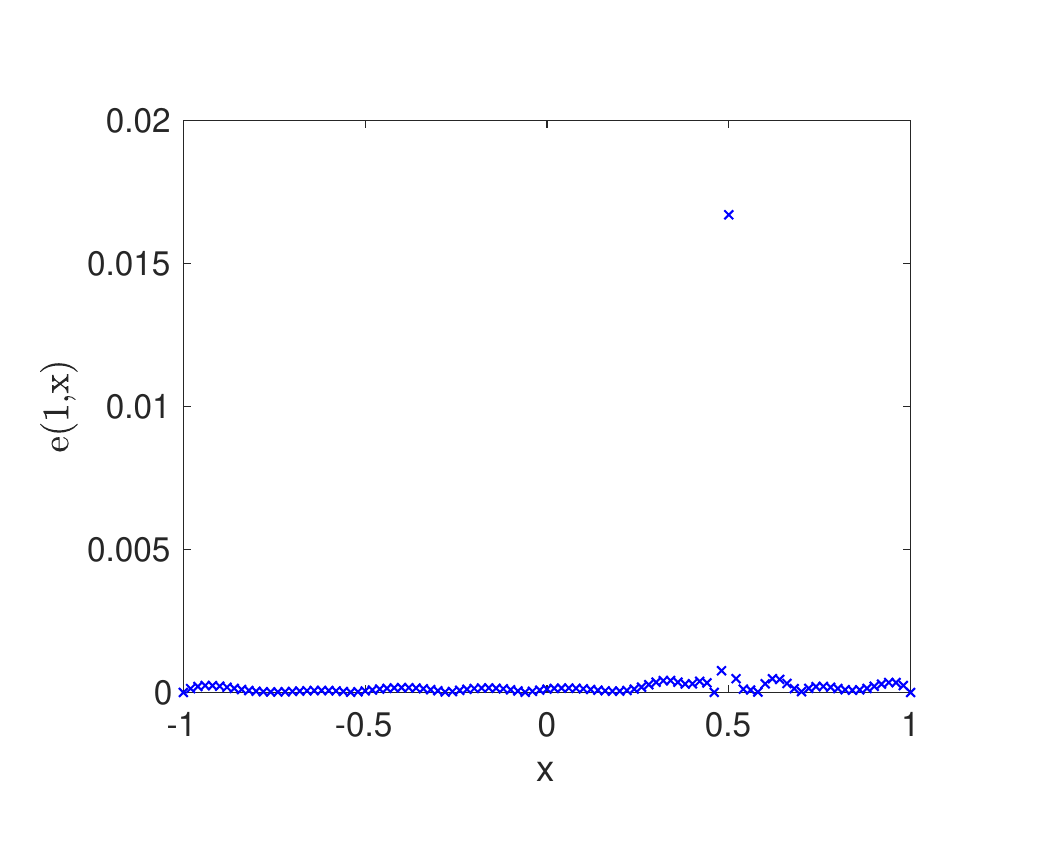}}
\subfigure[relative error when $x_0=0.5$]
 { 
\includegraphics[width=0.32\textwidth]{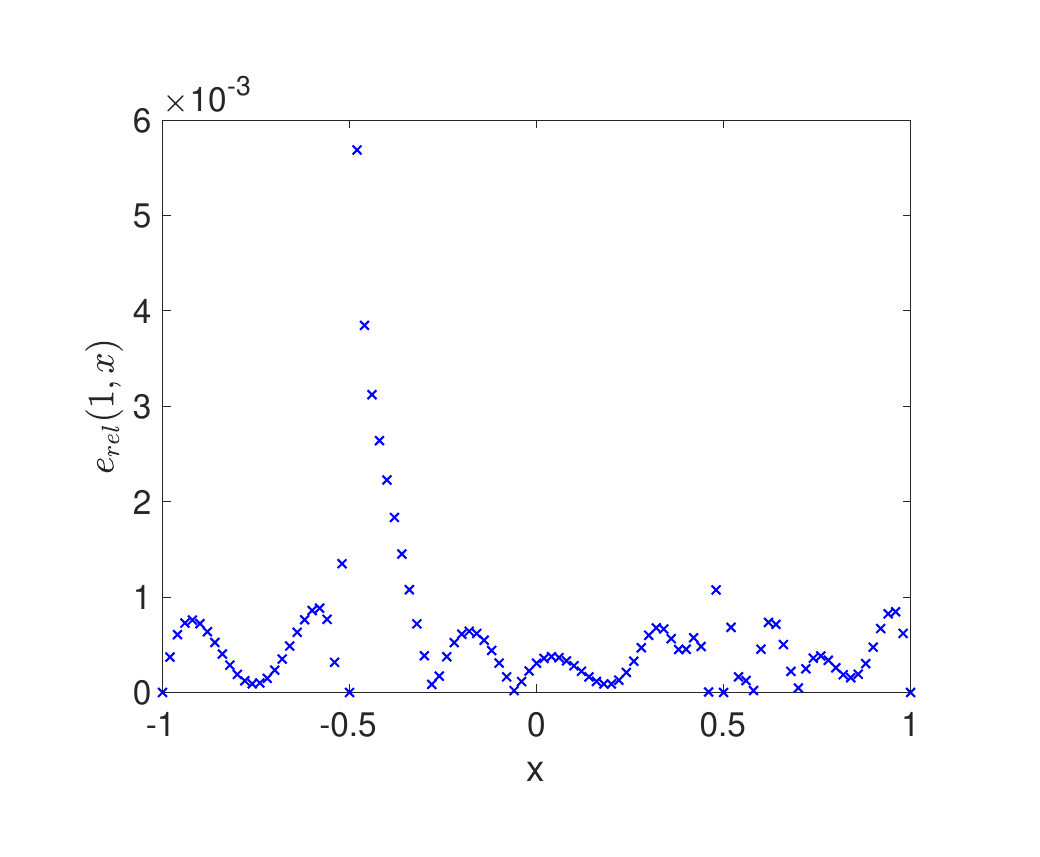}}
\subfigure[absolute error around $x=0.5$ when $x_0=0.5$]
 { 
\includegraphics[width=0.32\textwidth]{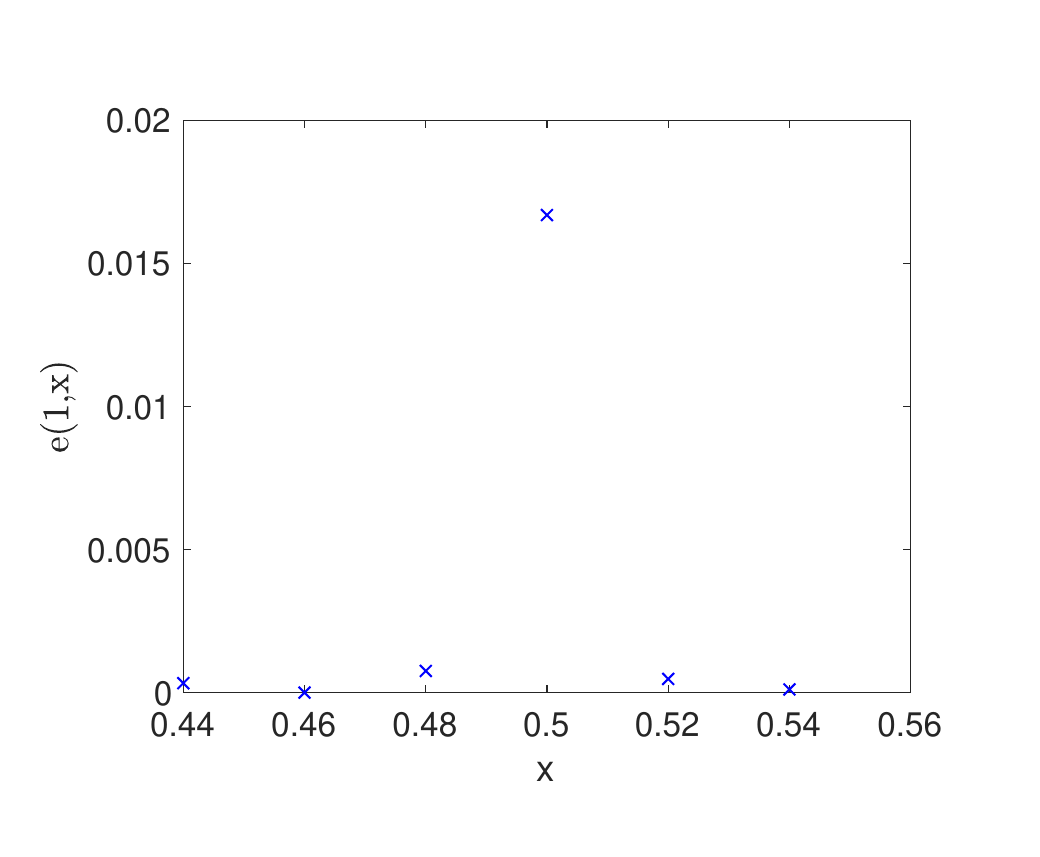}}
\subfigure[relative error around $x=0.5$ when $x_0=0.5$]
 { 
\includegraphics[width=0.32\textwidth]{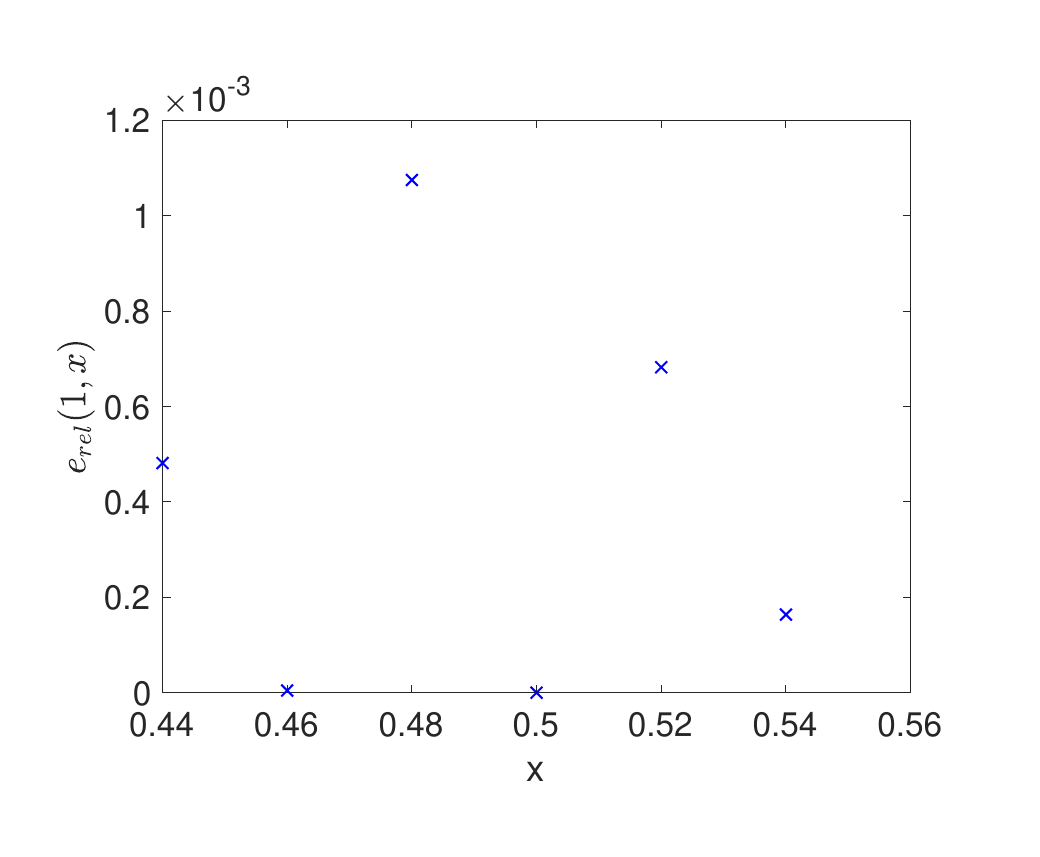}}

%
 \caption{Numerical results for Example \ref{exm:burgers} when {${\epsilon=10^{-2}/\pi}$}  using $(x^2-1)N(t, x,x/ \sqrt{\epsilon}, 1/\sqrt{\epsilon})+g(t)$ at $t=1$. $\alpha=1, \alpha_1=0$, $N_c=22500$ on the plane $(t,x)$, epochs=60000, with constant learning rate $10^{-4}$ (pre-trained with ${ \epsilon_0}=10^{-1}/\pi$, epochs=30000). The relative error is set to zero at $x=-1,0,1$ (when $x_0=0)$ and at $x=\pm 0.5$ (when $x_0=0.5$), where $\vert u(1,x)\vert<10^{-14}$. } 
 \label{fig:burgers_1e-2_pi} 
 \end{figure}

We employ the two-scale NN that enforces the Dirichlet boundary condition, i.e., $(x^2-1) N(t,x,x/\sqrt{\epsilon}, 1/\sqrt{\epsilon})+g(t)$ to solve the Burger's equation under a small viscosity parameter $\epsilon$. The NN size is $(4,20,20,20,20,1)$. 

First, we discuss the case when $\epsilon=10^{-2}/\pi$. For $x_0=0$, i.e., the inner layer (nearly a shock profile) is located at the center of the domain. According to \eqref{eq:exact_burgers}, $x_0=0$ leads to $g(t)=0$. We employ the successive training strategy in Algorithm \ref{alg:succesive-training-two-scale} by initializing the parameters of this NN using the pre-trained parameters obtained from training with ${\epsilon_0}=10^{-1}/\pi$. We collect the numerical results in Figure \ref{fig:burgers_1e-2_pi}. From Figure \ref{fig:burgers_1e-2_pi}(e), we can see that the maximal relative error at the final time node $t=1$ is 2.5\%. The absolute error around the inner layer, i.e., around $x=0$, is only at the level of $10^{-3}$, which accounts for a mere 0.4\% of the amplitude of the inner layer. 

In the scenario where the inner layer is not at the center of the domain, i.e., $x_0=0.5$, we collect the results
in Figures \ref{fig:burgers_1e-2_pi}(b), (f), (g), (h) and (i), indicating an excellent agreement between the NN and exact solutions. We can see that the maximal relative error at the final time node $t=1$ is within 0.6\%, and the absolute error around  $x=0.5$ is at the level of $10^{-2}$, representing only 1\% of the amplitude of the inner layer. The results indicate that the two-scale NN method can accurately capture the exact solution of the Burgers' equation when $\epsilon=10^{-2}/\pi$, regardless of whether the inner layer is near the center of the spatial domain or not. 

\begin{figure}[!ht]
	\centering
        \includegraphics[width=0.36\textwidth]{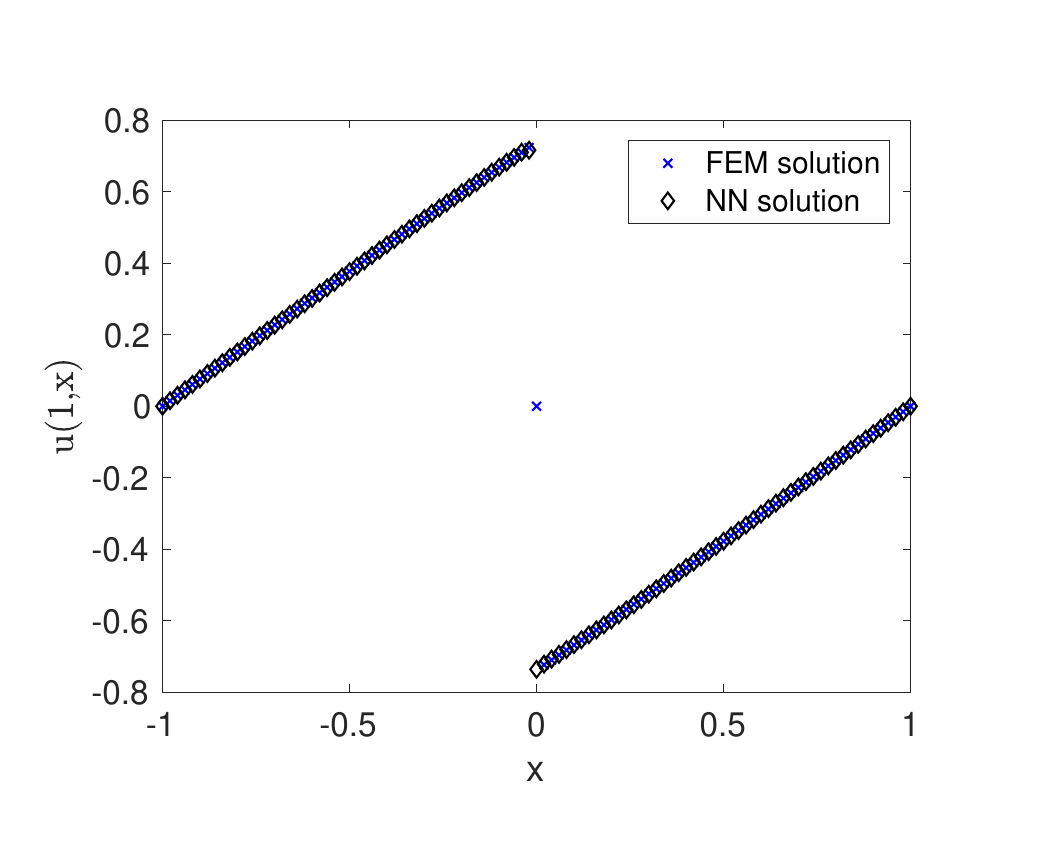}
        \includegraphics[width=0.36\textwidth]{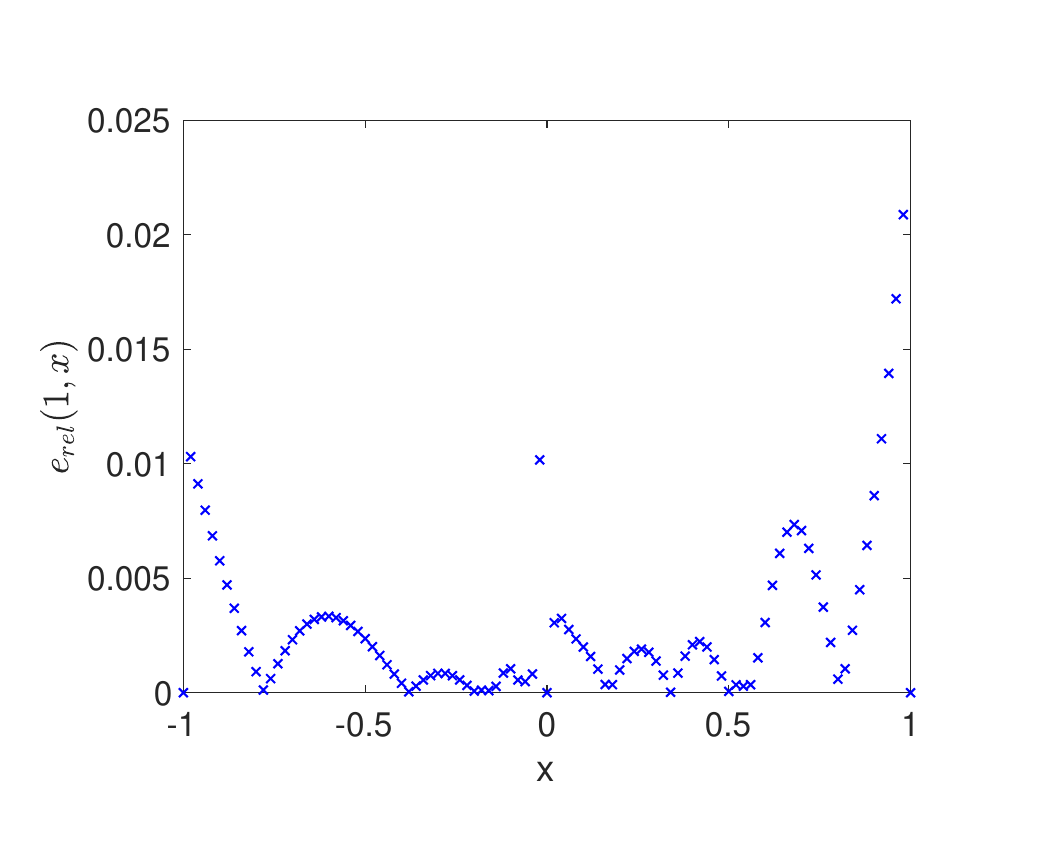}

	\caption{\textbf{left}: high resolution FEM solution with mesh size 1/500, polynomial degree 2 and NN solution obtained with $N(t, x,x/ \sqrt{\epsilon}, 1/\sqrt{\epsilon})$  for Example \ref{exm:burgers} when $\epsilon=10^{-3}/\pi, x_0=0$ at $t=1$. $\alpha=1$, $N_c=22500$ on the plane $(t,x)$, epochs=60000 (two-phase successive training: pre-trained with ${\epsilon_0}=10^{-2}/\pi$, epochs=60000; the pre-trained parameters are obtained from training with ${\epsilon_0}=10^{-1}/\pi$, epochs=30000); \textbf{right}: relative error.  }
\label{fig:burgers_1e-3_pi_soln_pretrain_wNNRBA}
\end{figure}

When $\epsilon=10^{-3}/\pi$, it is difficult to achieve accurate numerical quadrature for the exact solution in \eqref{eq:exact_burgers}  because the integrand can vary from tiny scales to infinity, driven by the very small $\epsilon$. Therefore, we employ the high-resolution finite element solution as the reference solution. The NN solution is obtained through a two-phase successive training process with Algorithm \ref{alg:succesive-training-two-scale}. To be specific, we initialize the NN parameters with the pre-trained parameters obtained from training with ${\epsilon_0}=10^{-2}/\pi$. These pre-trained parameters themselves were initially obtained from the training process with ${ \epsilon_0}=10^{-1}/\pi$. As observed in Figure \ref{fig:burgers_1e-3_pi_soln_pretrain_wNNRBA}, there is a shift between the exact solution and the predicted solution by the two-scale NN method within a very narrow band around the inner layer at $x=0$. The observed deviation primarily arises from the inherent difficulty of the problem itself. Indeed, given $\epsilon=10^{-3}/\pi$, the inner layer practically becomes a vertical line. Nevertheless, the two-scale NN solution still captures the key features of the exact solution in this context. 

\begin{exm}[2D steady-state convection-diffusion problem]\label{exm:zzm}
\begin{eqnarray*}
-\epsilon \Delta u+(u_x+u_y)&=& f ~\text { in } \Omega=(0,1)^2,\\
u&=&0 ~\text{on}~\partial \Omega.
\end{eqnarray*}
The exact solution is given by:
$\displaystyle u = xy\left(1 - \exp\left(-\frac{1-x}{\epsilon}\right) \right) \left(1 - \exp\left(-\frac{1-y}{\epsilon}\right)\right),$
and the right-hand side function is then:
\[ f=(x+y) \left(1-\exp\left(-\frac{1-x}{\epsilon}\right)\exp\left(-\frac{1-y}{\epsilon}\right) \right) + (x-y)\left(\exp\left(-\frac{1-y}{\epsilon}\right)-\exp\left(-\frac{1-x}{\epsilon}\right)\right).
\]
\end{exm}



\begin{figure} [!ht]
\centering
\subfigure[right-hand-side function]
{
\includegraphics[width=0.31\textwidth]{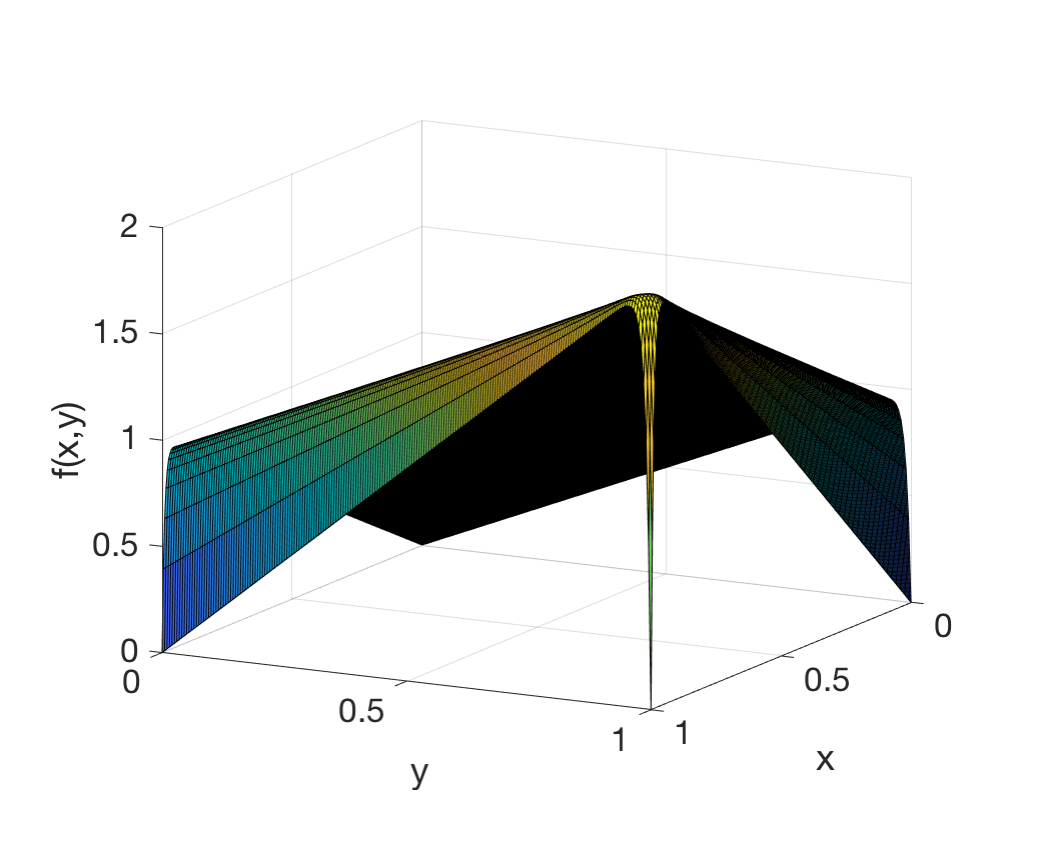}
}
\subfigure[exact solution]
{
\includegraphics[width=0.32\textwidth]{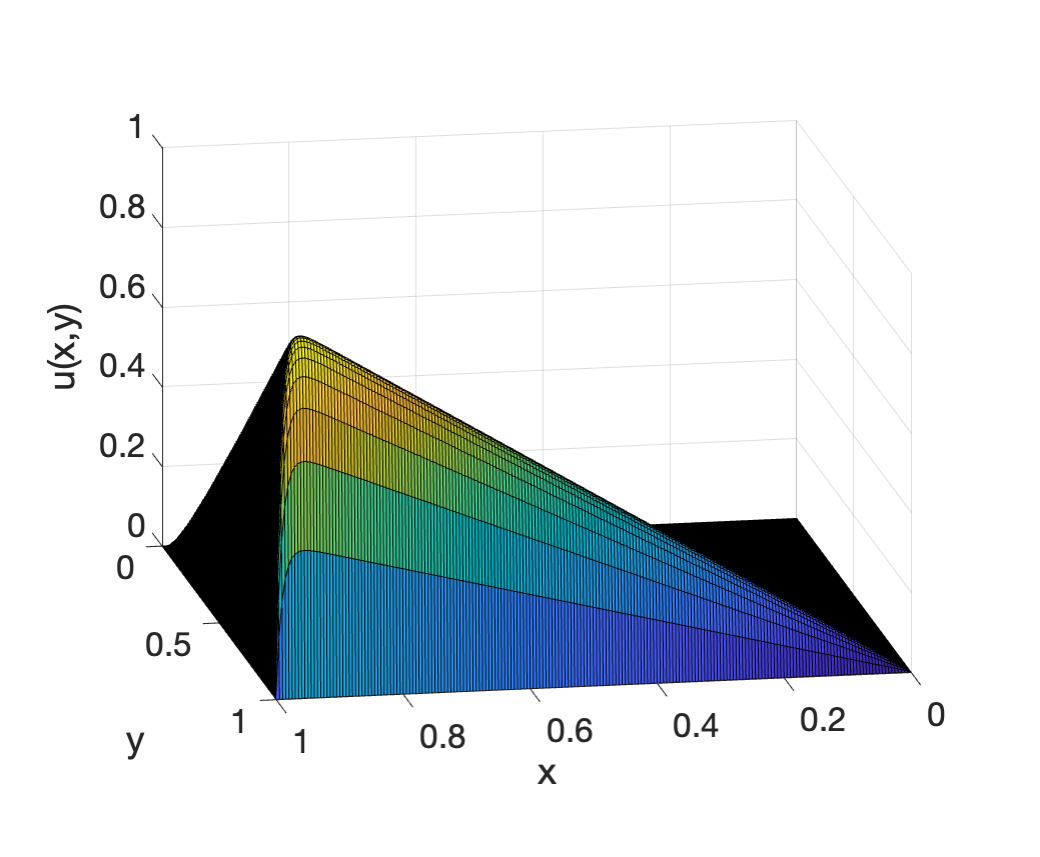}

}
\subfigure[NN solution]
 { 
 \includegraphics[width=0.32\textwidth]{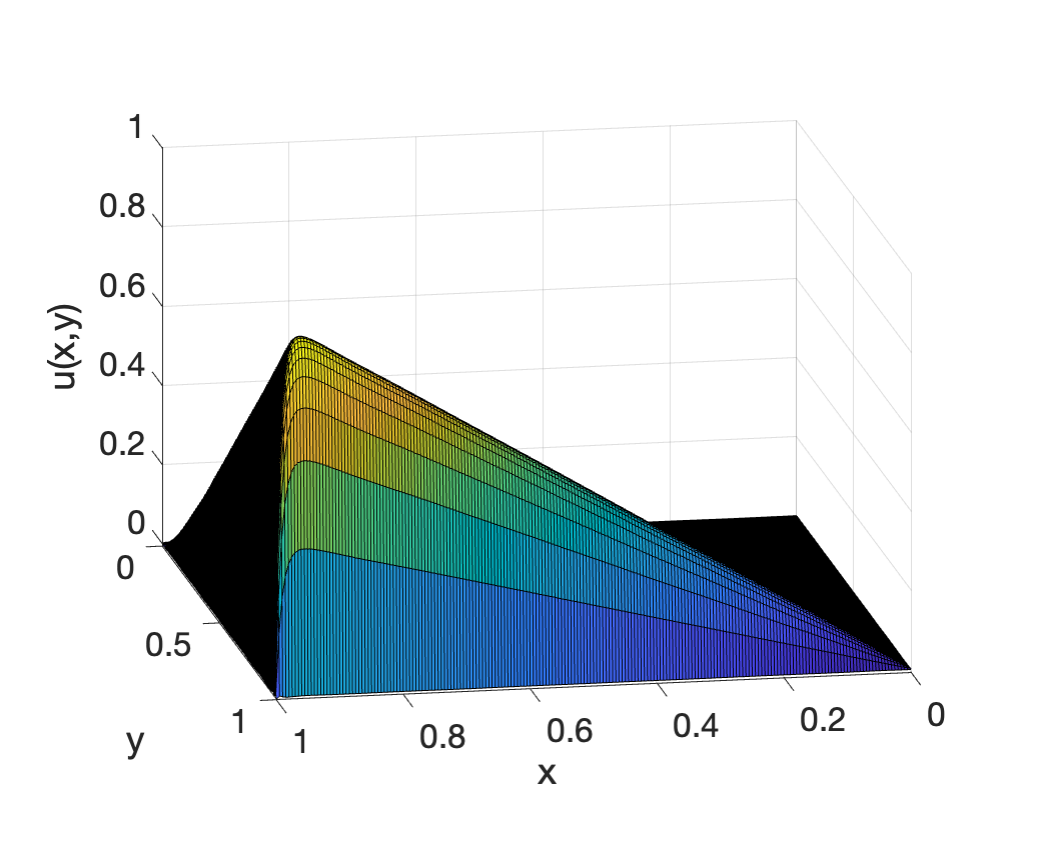}}
 \subfigure[absolute error ]
{
\includegraphics[width=0.31\textwidth]{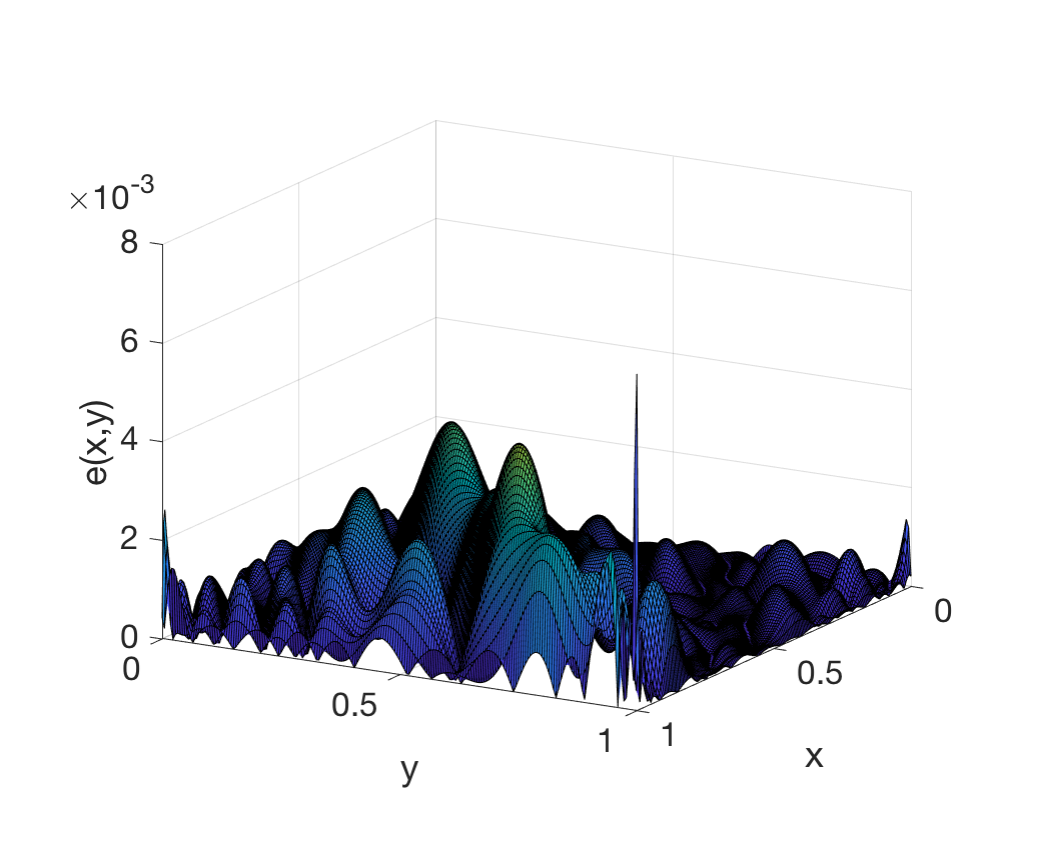}}
\subfigure[relative error]
 { 
 \includegraphics[width=0.31\textwidth]{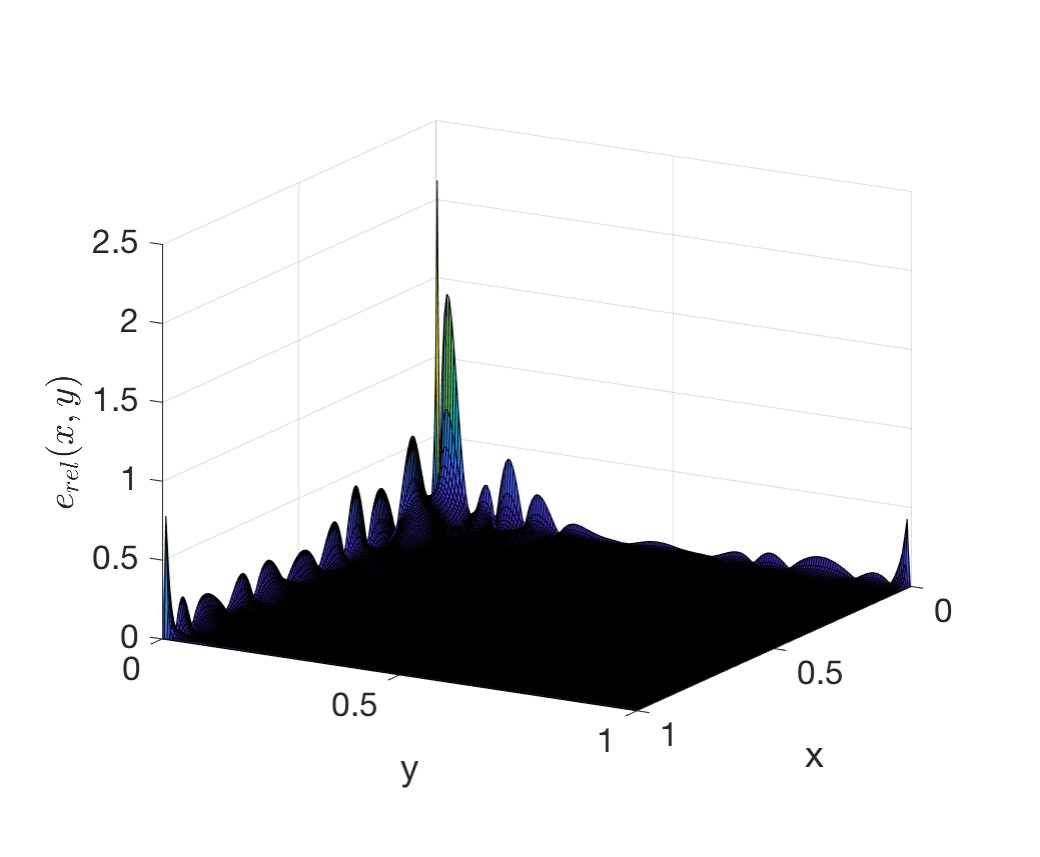}}
\subfigure[{ loss history}]
 {  \includegraphics[width=0.28\textwidth]{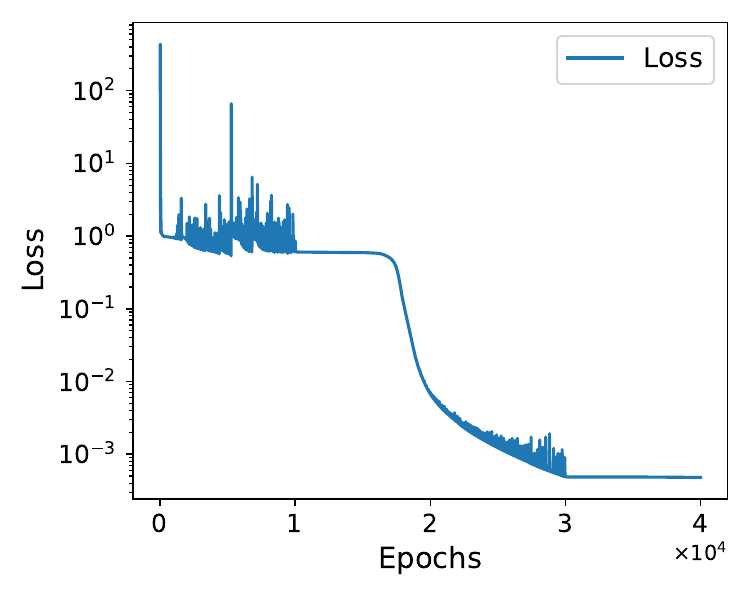}}


   \subfigure[{ errors history}]
   {
\includegraphics[width=0.27\textwidth]{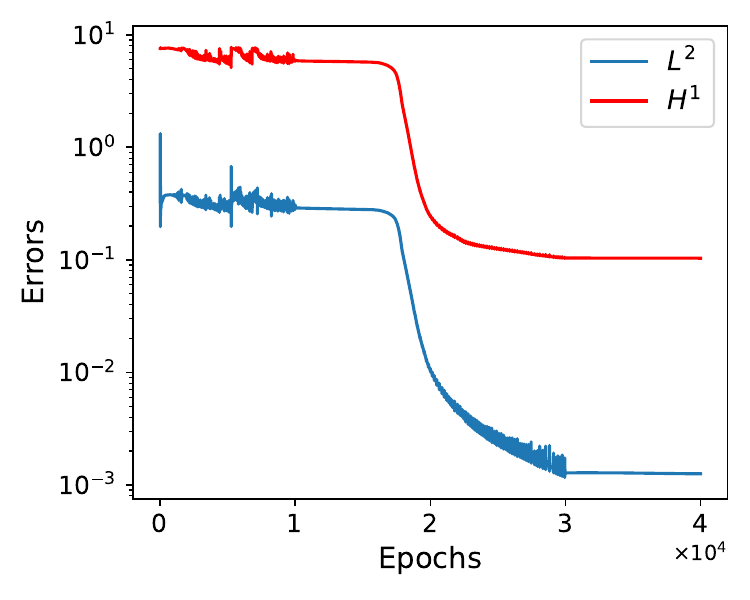}}
 \subfigure[absolute error around $x=1$, $y=1$. (top view)]
{
\includegraphics[width=0.28\textwidth]{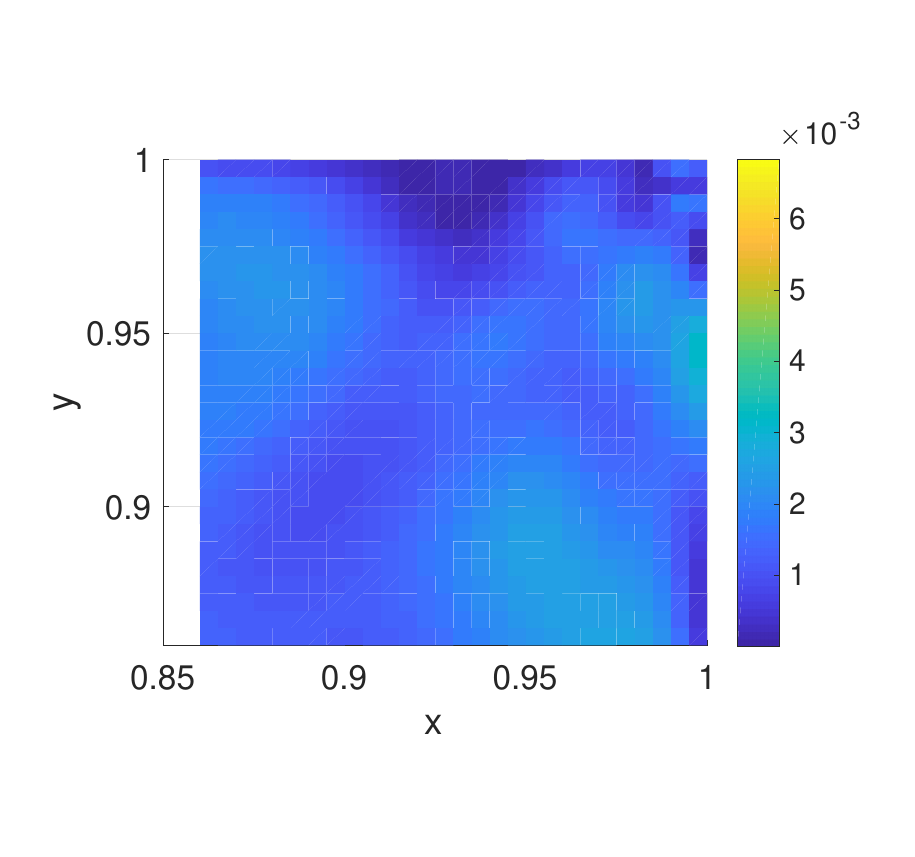}}
 \hspace{ 10pt}
\subfigure[relative error around $x=1$, $y=1$. (top view)]
 {  \includegraphics[width=0.28\textwidth]{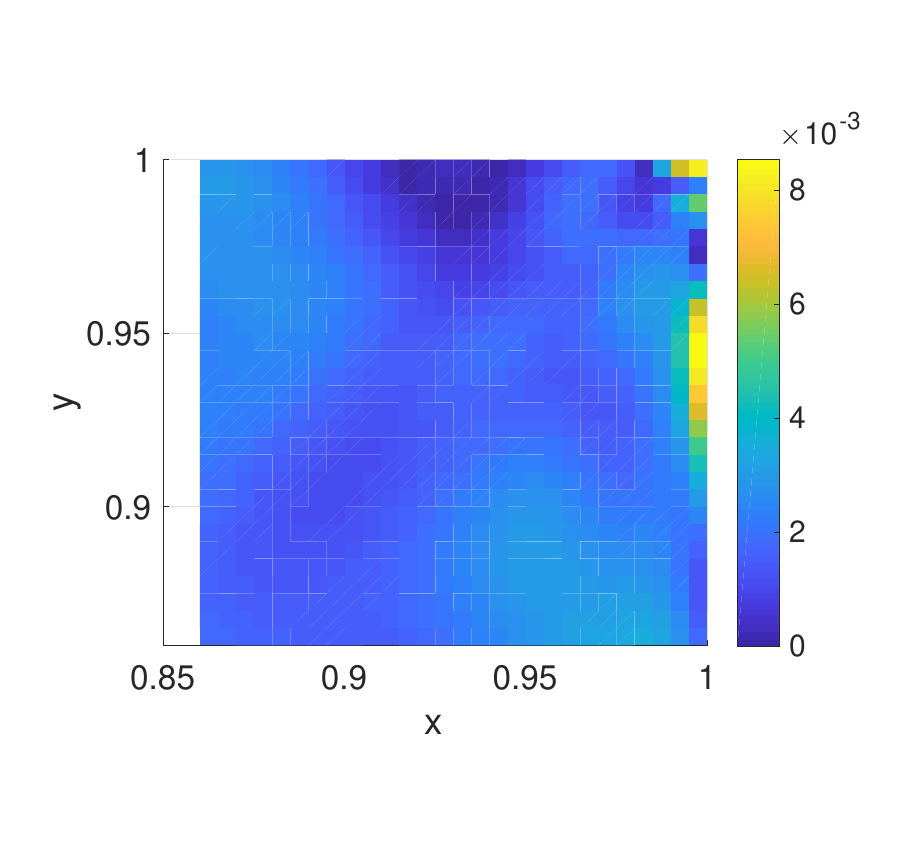}}

 	\caption{ numerical results for Example \ref{exm:zzm} when $\epsilon=10^{-2}$ using  $N(x, (x-0.5)/\sqrt{\epsilon}, 1/\sqrt{\epsilon}, y, (y-0.5)/\sqrt{\epsilon}, 1/\sqrt{\epsilon})$ with $\alpha=100$, $\alpha_1=0, N_c=22500,$ epochs=40000. We enforce $e_{rel}(x,y)=0$ on $\partial \Omega$ where the exact solution $u=0$. The relative $l^2$ error is $e_{rel}^{l^2}=  1.5722\times 10^{-5}.$ }
\label{fig:results_zzm}
 \end{figure}

According to \cite{ZhangFEMSK03}, the exact solution has two boundary layers around the outflow boundaries $x=1$ and $y=1$. The two-scale NN we employ here is $N(x, (0.5-x)/\sqrt{\epsilon}, 1/\sqrt{\epsilon}, y, (0.5-y)/\sqrt{\epsilon}, 1/\sqrt{\epsilon}))$, with the size $(6,64, 64, 64,1)$. The numerical results are collected in Figure \ref{fig:results_zzm} for $\epsilon=10^{-2}$. The figures demonstrate that the NN solution accurately captures the exact solution, including the two boundary layers $x=1$ and $y=1$. Specifically, in Figure \ref{fig:results_zzm}(h), we observe that the relative error around these boundary layers remains below 0.8\%. 

It is important to note that the challenge posed by this problem extends beyond the presence of boundary layers. Another critical aspect lies in the rapid variation of the right-hand-side function $f$, occurring at $x=1$ and $y=1$, as illustrated in Figure \ref{fig:results_zzm}(a). Even so, the prediction by the two-scale NN method is accurate. 

\begin{exm}[2D {Helmholtz} problem] \label{exm:1p14}
Consider the boundary value problem
$$
\begin{aligned}
-\Delta u -k^2 u &=f & & \text { in } \Omega=(-1,1)^2 \\
u & =g & & \text { on } \partial \Omega .
\end{aligned}
$$
We consider the following two cases where the exact solutions are respectively given by
\begin{align*}
a) \quad u  &= \sin \left(a_1 \pi x\right) \sin \left(a_2 \pi y\right).\\
b) \quad u= & \frac{\cos (k r)}{k}-\frac{J_0(kr) (J_0(k)\cos k+ J_1(k)\sin k)}{k\left(J_0^2(k)+ J_1^2(k)\right)},
\end{align*}
where $r=\sqrt{x^2+y^2}$, $J_{\nu} (\cdot)$ are Bessel functions of the first kind. These two cases represent the oscillation pattern of the solutions in the $x-y$ direction and the radial direction, respectively.

\end{exm}

For \textbf{Example \ref{exm:1p14} a)}, we take $\epsilon=1/\max \{ a_1, a_2\}$. In this regard, we address the Helmholtz problem numerically in two scenarios. First, we set $(a_1, a_2)=(1,4)$, which corresponds to a modest derivative in the solution and a relatively lower frequency for the time-harmonic wave. Secondly, we consider the case of $(a_1, a_2)=(2,3)$, representing a larger derivative in the solution and a higher frequency for the wave. We employ a two-scale NN framework $(x^2-1)(y^2-1)N(x, x/\sqrt{\epsilon}, 1/\sqrt{\epsilon}, y, y/\sqrt{\epsilon}, 1/\sqrt{\epsilon})$, where we enforce the homogeneous Dirichlet boundary condition by the construction of the NN.
The size of the NN is $(6,10,10,10,10,1)$.  The results are collected in Figure \ref{fig:exm1p14_k4a14}. Observing from Figure \ref{fig:exm1p14_k4a14}(a), (b), (d), and (e), the NN solution matches the exact solution well. 

As a comparison to the two-scale NN solution, when $(a_1,a_2)=(1,4)$, we present the solution obtained by the one-scale NN, i.e.,  $(x^2-1)(y^2-1)N(x, y)$ in Figure \ref{fig:exm1p14_k4a14}(c), with NN size (2,10,10,10,10,1). The figure shows that the one-scale NN solution fails to capture the oscillation pattern in the exact solution and deviates significantly from the expected outcome.

 { We also compare the two-scale NN solutions to those from multi-level neural networks ({ MLNN} )\cite{multilvNN23}, known for their capability to address high-frequency issues. In the multi-level neural network method, a critical technique for addressing high-frequency issues involves scaling the equation's source term with a factor. To be specific, at the initial neural network level, $u_0$ is solved from the equation with a scaling factor $\mu_0$: $R_0(x,u_0)=\mu_0 f-\mathcal{L}u_0=0$. The NN solution to this equation is denoted as $\tilde{u}_0$, resulting in the initial approximation $\tilde{u}=\tilde{u}_0/\mu_0$.
In the subsequent $i$-th level ($i\geq 1$) of the neural network, $u_i$ is solved from $R_i(x,u_i)=\mu_i R_{i-1} (x,\tilde{u}_{i-1})-\mathcal{L}u_i=0$. The NN solution for this level is denoted as $\tilde{u}_i$.
The updated approximation is then obtained by summing the initial approximation with the cumulative corrections up to the $i$-th level, i.e., $\displaystyle \tilde{u}= \sum_{i=0}^L \frac{\tilde{u}_i} {\prod_{k=0}^{i} \mu_k}$.
For more details of the algorithm of the multi-level neural networks, we refer the readers to Section 4 of \cite{multilvNN23}.  

We specifically investigate the case $(a_1,a_2)=(1,4)$ with the multi-level neural networks. The hyper-parameters used for multi-level NN are summarized in Table \ref{table:multilvNN_hyper}, which are suggested in \cite{multilvNN23} for 2D numerical examples. The networks are first trained with Adam,  followed by refinement with L-BFGS.


{ 
It is recommended in \cite{multilvNN23} to set  {$\mu_i, (i \geq 1)$} as the inverse of the amplitude of $u_{E, i}$. Here, $u_{E, i}$ is a rough approximation of the solution $\hat{u}_i$ to $\mathcal{L}\hat{u}_i-R_{i-1}=0$ by the extreme learning machine (ELM) \cite{ELM2011}. 
}
%
{
For $\tilde{u}_i$, it is of the form:
\begin{equation}\label{fourier_2d}
\sin(\boldsymbol{\omega_M} \pi (1-x))\sin(\boldsymbol{\omega_M} \pi (1-y)) N(\gamma(x), \gamma(y)),
\end{equation}
where the Fourier feature:
$\gamma(x) = [\cos(\boldsymbol{\omega_M} x),  \sin(\boldsymbol{\omega_M} x)],$ and ${\boldsymbol{\omega_M}}=(\omega_1, ..., \omega_M), $
$M$ is specified in Table \ref{table:multilvNN_hyper}.

We compare the two-scale NN solution to the MLNN's initial and first levels of approximated solutions: $\tilde{u}_0/\mu_0$ and $\tilde{u}_0/\mu_0+\tilde{u}_1/\mu_0\mu_1$, with the hyper-parameters in Table \ref{table:multilvNN_hyper}.
The two-scale NN results are obtained using the NN size of $(3,10,1)$, trained for 2500 epochs with Adam, followed by 200 epochs using L-BFGS.
The MLNN results are obtained both without and with the Fourier feature in \eqref{fourier_2d}.
We denote the absolute error between the exact solution and the $i$-th level approximated solution as $e_i=\abs{u-\sum_{k=0}^i \left(\tilde{u}_k/\prod_{j=0}^{k} \mu_j \right)}$. 
The loss history is scaled by dividing $(\prod_{k=0}^{i} \mu_k)^2$, where $i$ denotes the level at which the current loss function values are being calculated.

As we can see in Figure \ref{fig:helm_k14_multilvNN_errors} (a) and (d) (e), the predictions by the two-scale NN and MLNN without Fourier feature provide reasonable accuracy in approximating the exact solution. 
The figure shows that the initial level prediction from the MLNN is on par with the two-scale NN when trained with the same number of epochs. The MLNN's first-level correction refines the accuracy to $10^{-3}$. 
%
%
The results of the MLNN with the Fourier feature \eqref{fourier_2d} shown in Figure \ref{fig:helm_k14_multilvNN_errors} (h) indicate that adding the Fourier feature significantly improves the accuracy. 

The results suggest that in addressing high-frequency issues, the two-scale NN does not necessarily surpass the approach of adding Fourier features or MLNN in accuracy; however, it offers a simpler solution by avoiding special treatments such as selecting wave numbers for Fourier features or scaling factors for MLNN. To further compare the two-scale NN method with MLNN, we consider two more 1D problems with one boundary layer, to be presented in Examples \ref{exm:burgers1d-steady} and \ref{exm:1d_nonlinear}.

}
\begin{table}[ht]
\centering

\begin{tabular}{l|c|c}
\hline
Hyper-parameters    & $\tilde{u}_0$ & $\tilde{u}_1$ \\
\hline
\# Hidden layers  & 2             & 2             \\
Widths of each hidden layer  & 10            & 20            \\
\# Adam iterations  & 2500          & 5000          \\
\# L-BFGS iterations & 200           & 400           \\
{\# Wave numbers M} & { 2}           & { 4}    \\
\hline
\end{tabular}
\caption{Hyper-parameters used in Example \ref{exm:1p14} a) with multi-level neural networks \cite{multilvNN23}}
\label{table:multilvNN_hyper}
\end{table}

\begin{figure} [!ht]
\centering
\subfigure[exact solution]
{
\includegraphics[width=0.23\textwidth]{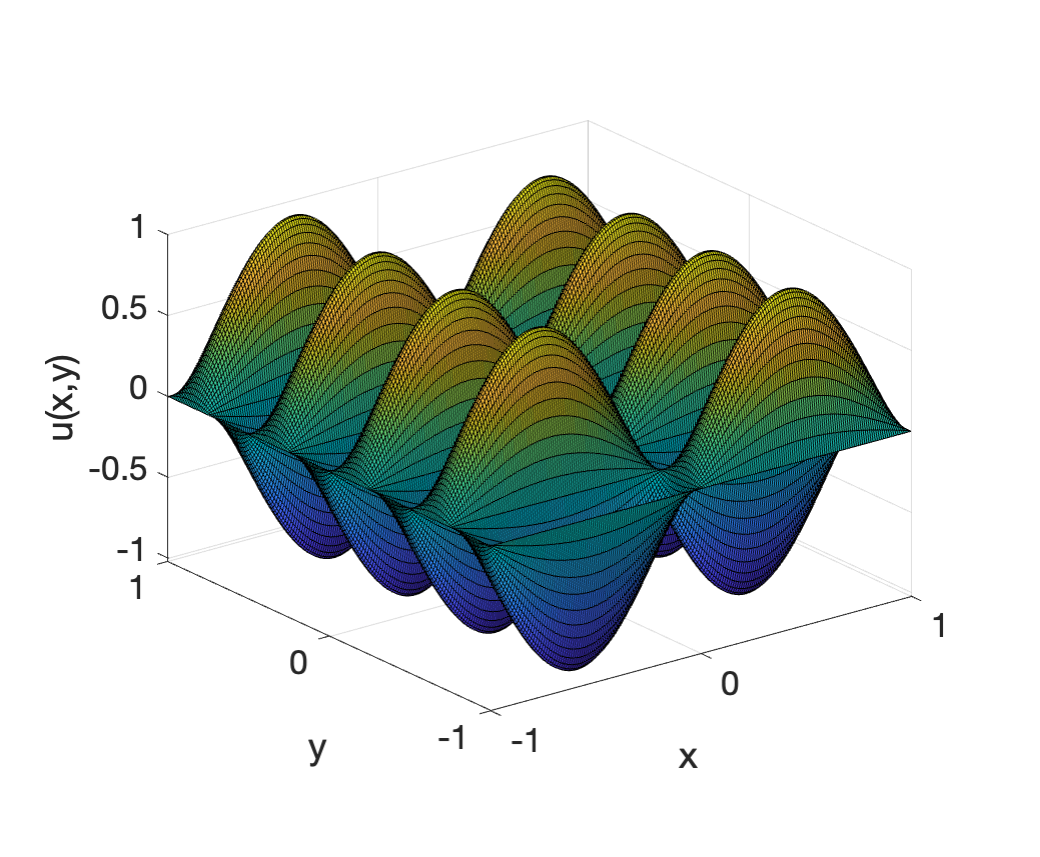}
}
\subfigure[two-scale NN solution]
{
\includegraphics[width=0.23\textwidth]{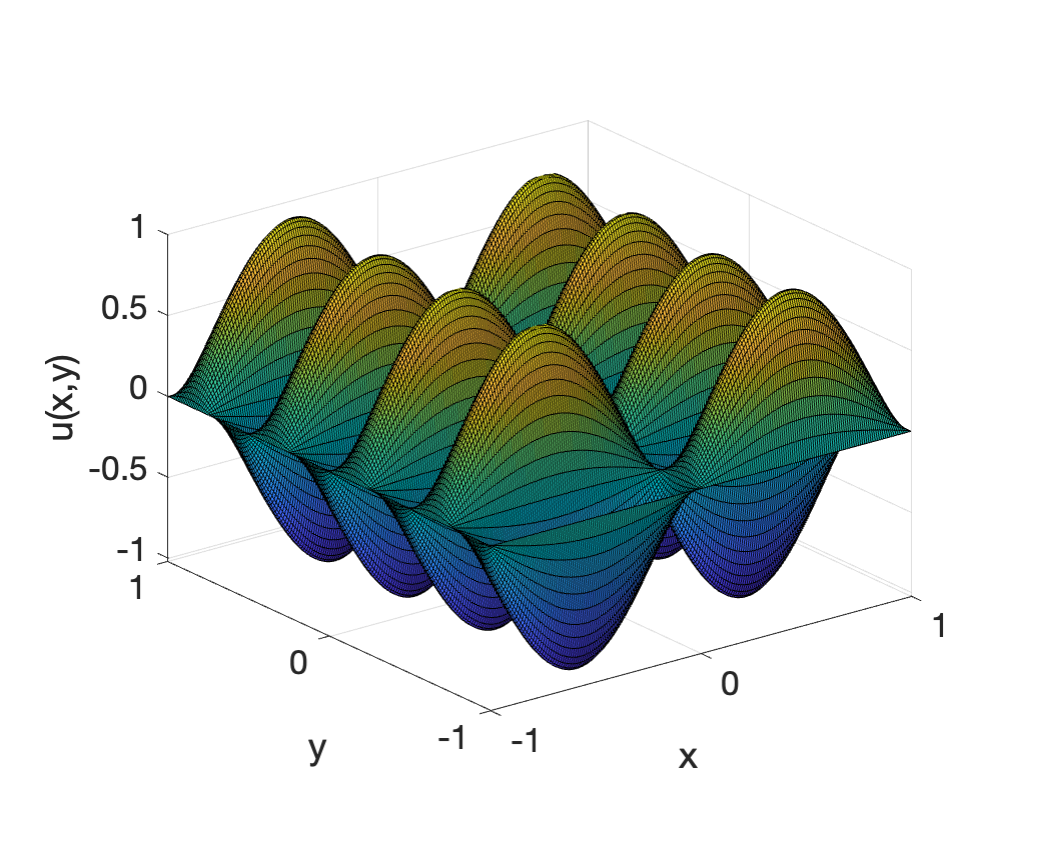}}
\subfigure[one-scale NN solution]
 { 
 \includegraphics[width=0.23\textwidth]{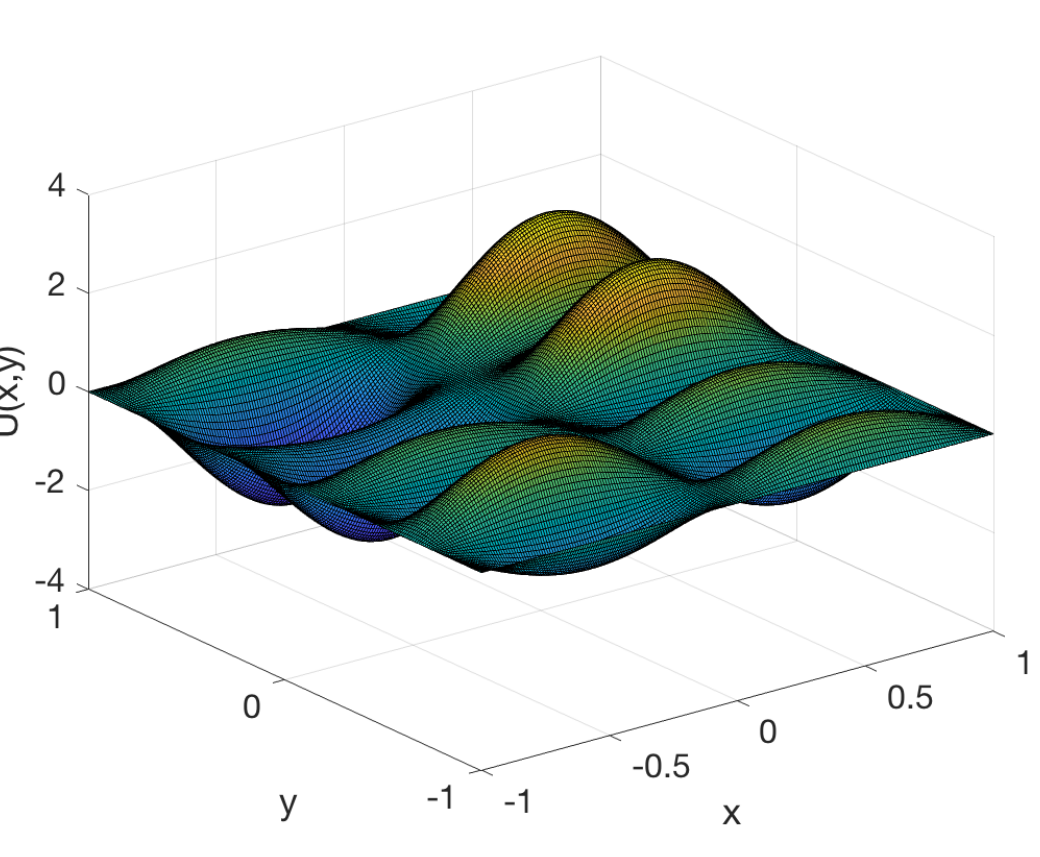}}
\subfigure[absolute error]
{
\includegraphics[width=0.22\textwidth]{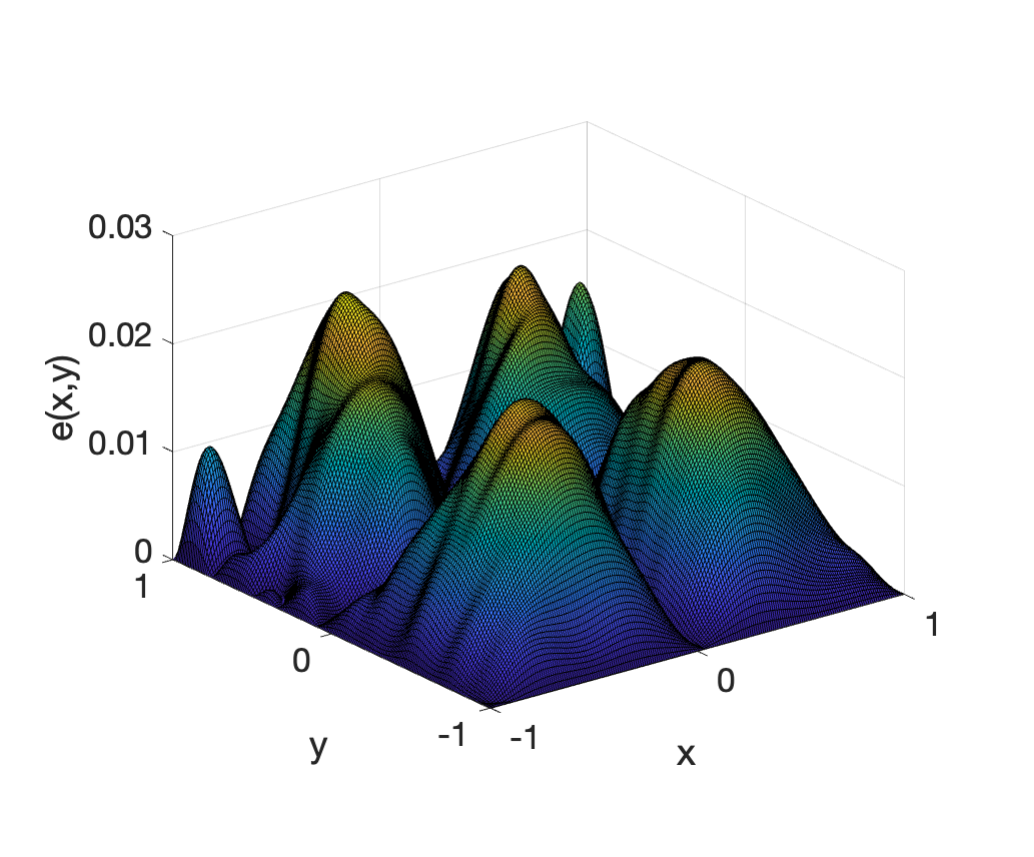}}

\subfigure[relative error]
 { 
 \includegraphics[width=0.3\textwidth]{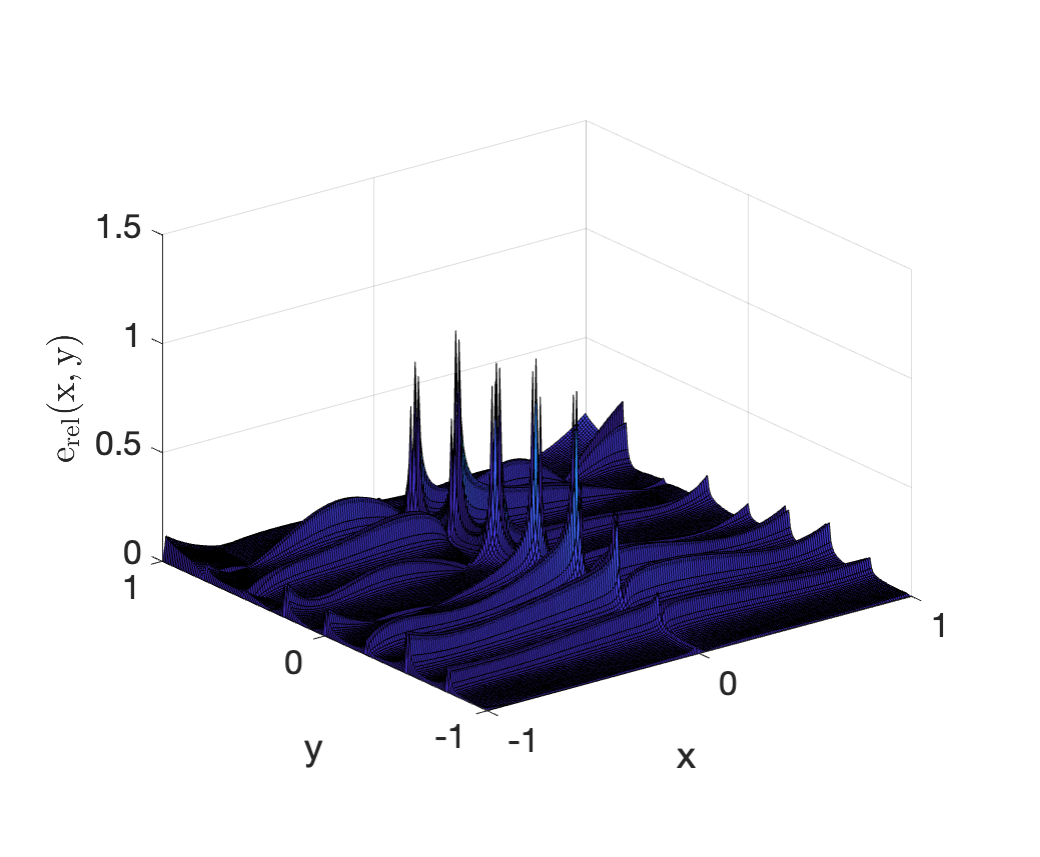}}
\subfigure[{ loss history}]
{ \includegraphics[width=0.28\textwidth]
 {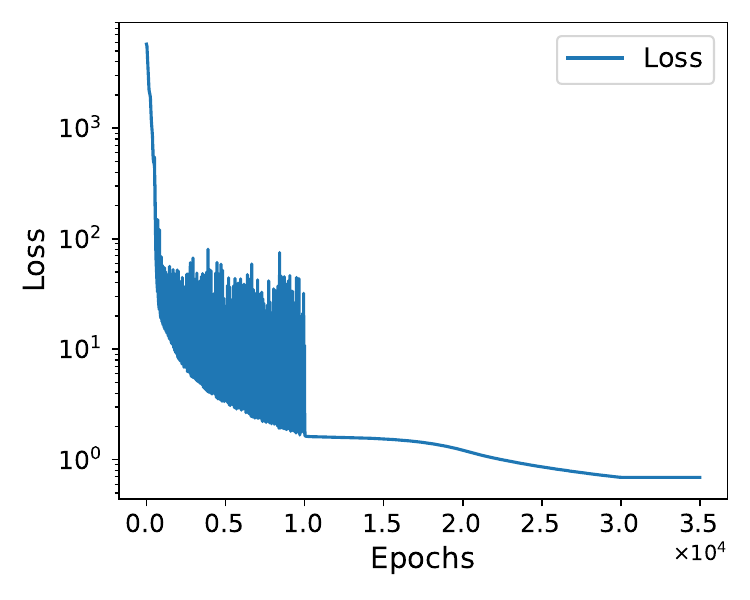}
 }
\subfigure[{ errors history}]
{ \includegraphics[width=0.28\textwidth]
 {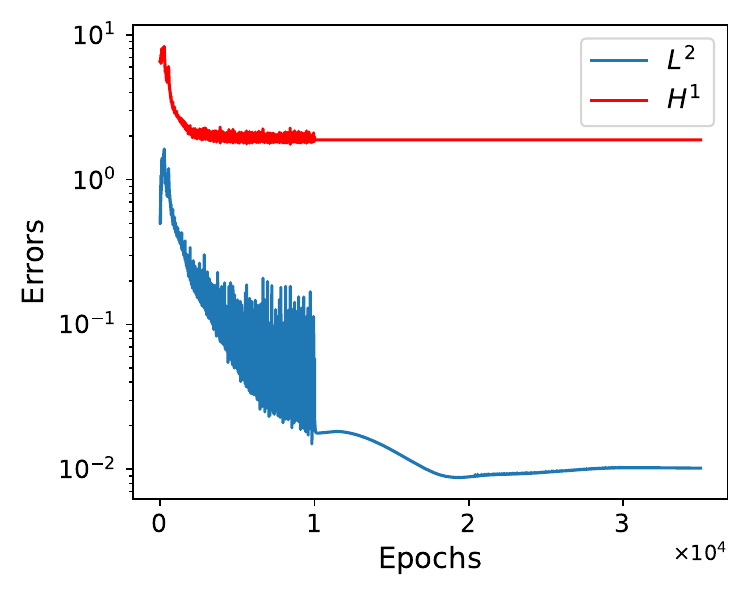}
 }

 \caption{\textbf{(b)(d)(e)(f)(g) are numerical results with the two-scale NN} for Example \ref{exm:1p14} a),  when $k=4, a_1=1, a_2=4$ using $(x^2-1)(y^2-1)N(x, x/\sqrt{\epsilon}, 1/\sqrt{\epsilon}, y, y/\sqrt{\epsilon}, 1/\sqrt{\epsilon})$,  $\alpha=1, \alpha_1=0, N_c=22500,  $ { epochs=35000}. We enforce $e_{rel}(x,y)=0$ on the nodes where the exact solution $u=0$. The relative $l^2$ error:   $e_{l^2}^{rel} =4.1907\times 10^{-4}.$\protect\\
 \textbf{(c) is the numerical result with the one-scale NN}, using $(x^2-1)(y^2-1)N(x, y)$, the other configurations are the same as the two-scale method. 
 }
 \label{fig:exm1p14_k4a14} 
 \end{figure}


\begin{figure}[!ht]
	\centering
\subfigure[absolute error of two-scale NN solution ]
{
\includegraphics[width=0.24\textwidth]{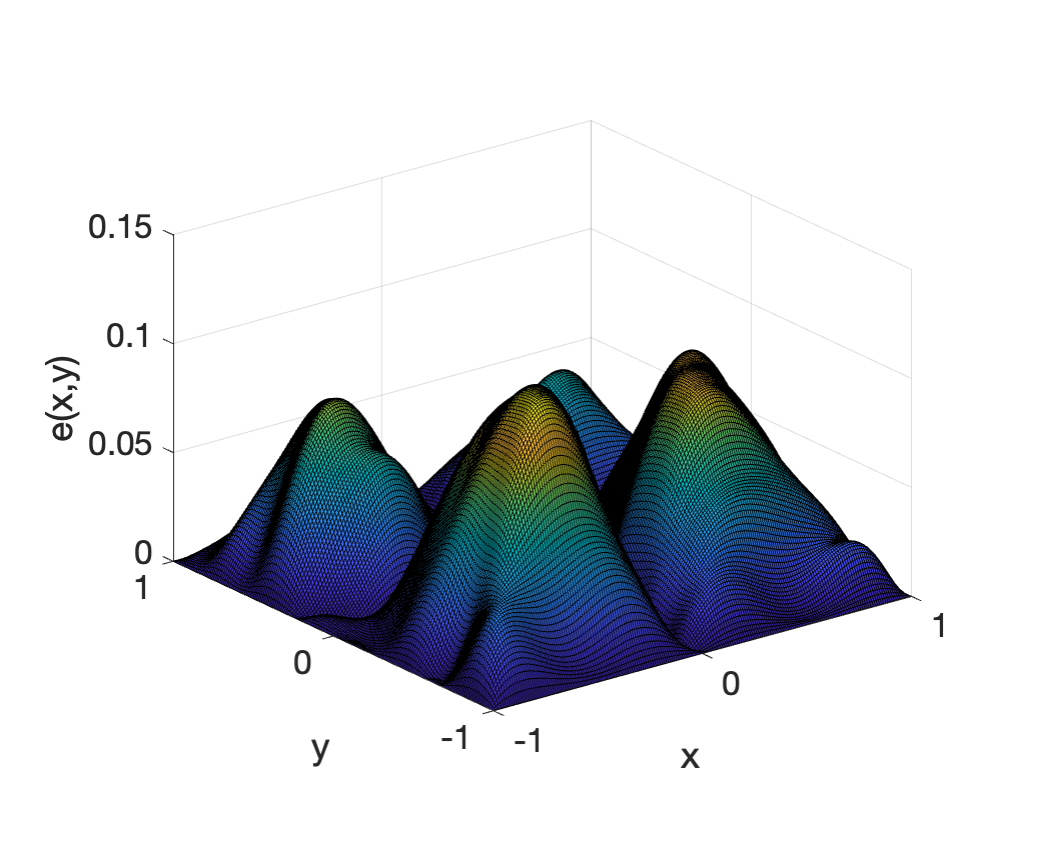}
}
\subfigure[loss history]
 { 
 \includegraphics[width=0.22\textwidth]{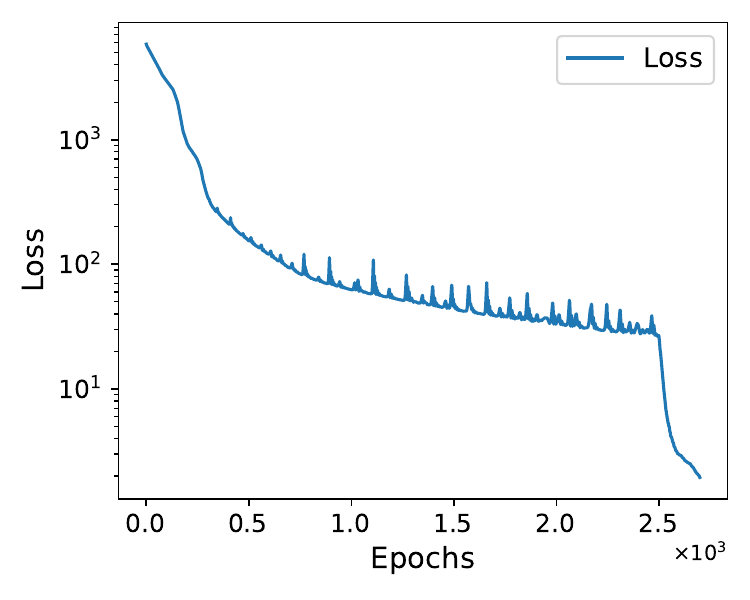}}
\subfigure[errors history]
{
\includegraphics[width=0.22\textwidth]{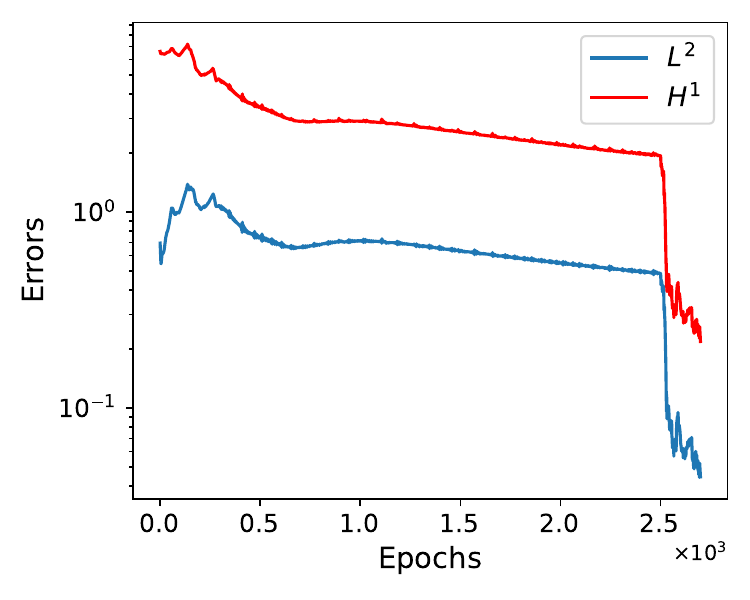}}

\subfigure[$e_0$ with MLNN \newline without Fourier feature]
{
\includegraphics[width=0.24\textwidth]{exm5/error_exm1p14_1slv0.pdf}
}
\subfigure[$e_1$ with MLNN \newline without Fourier feature]
{
\includegraphics[width=0.24\textwidth]{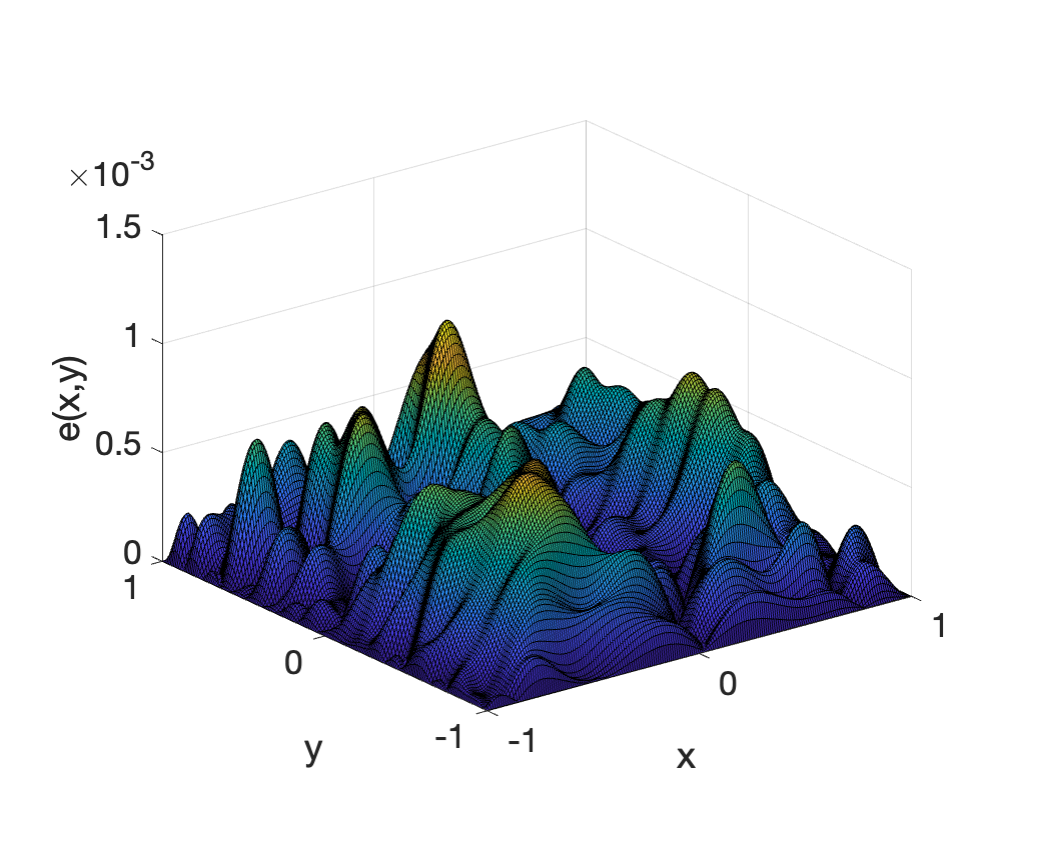}}
\subfigure[loss history]
 { 
 \includegraphics[width=0.22\textwidth]{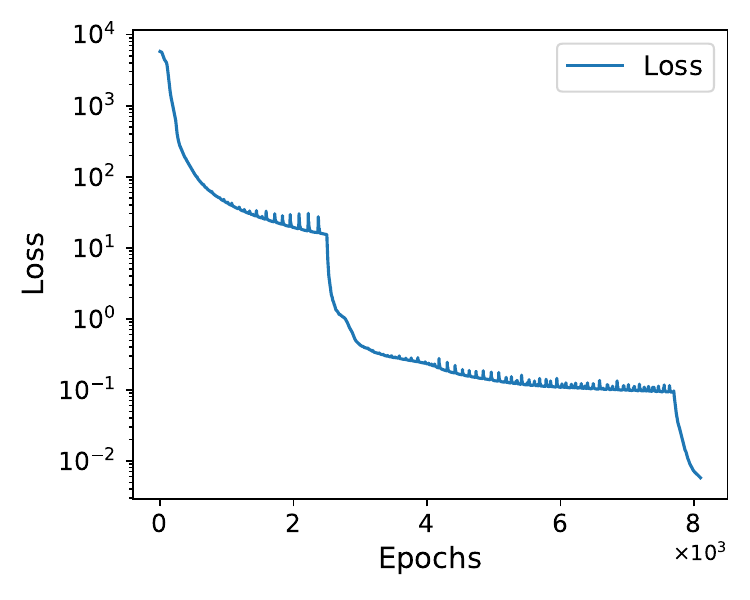}}
\subfigure[errors history]
{
\includegraphics[width=0.22\textwidth]{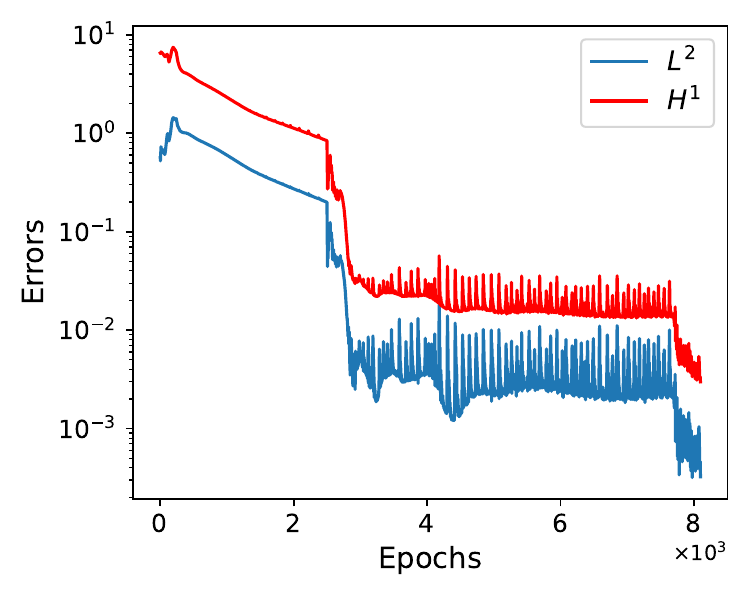}}

\subfigure[$e_0$ (left) and $e_1$ (right) with MLNN, adding Fourier feature]
{
\includegraphics[width=0.25\textwidth]{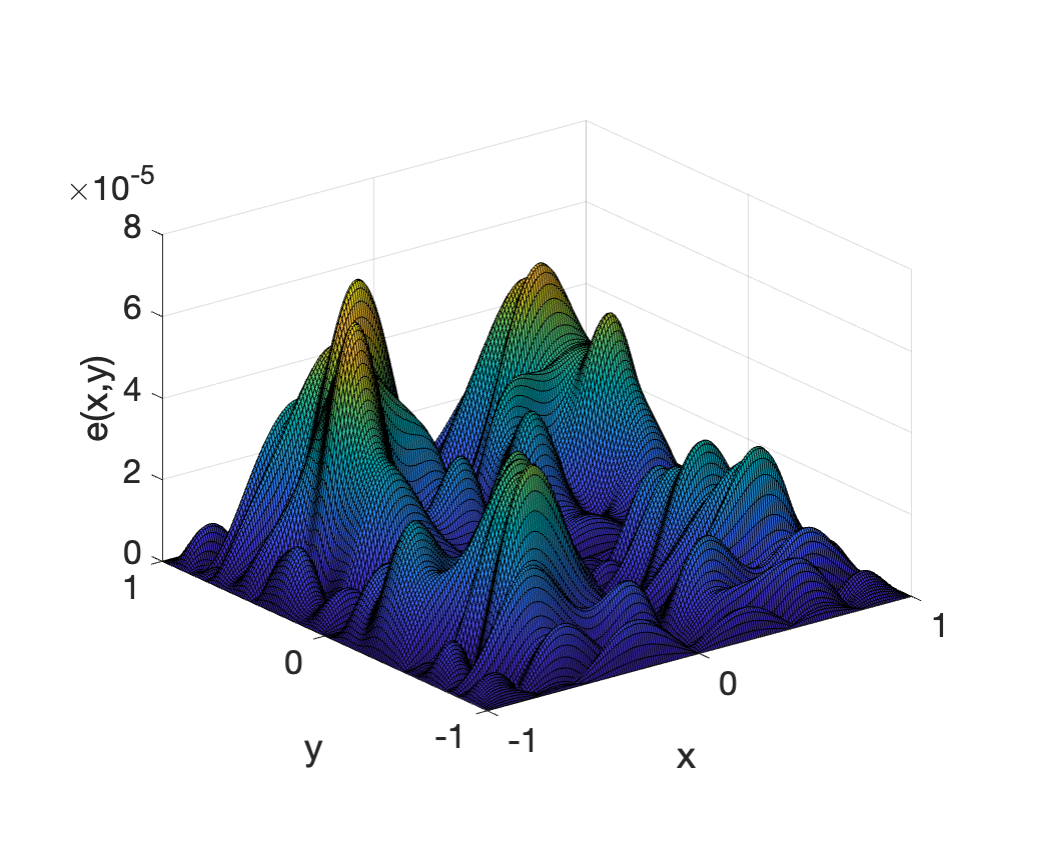}
\includegraphics[width=0.25\textwidth]{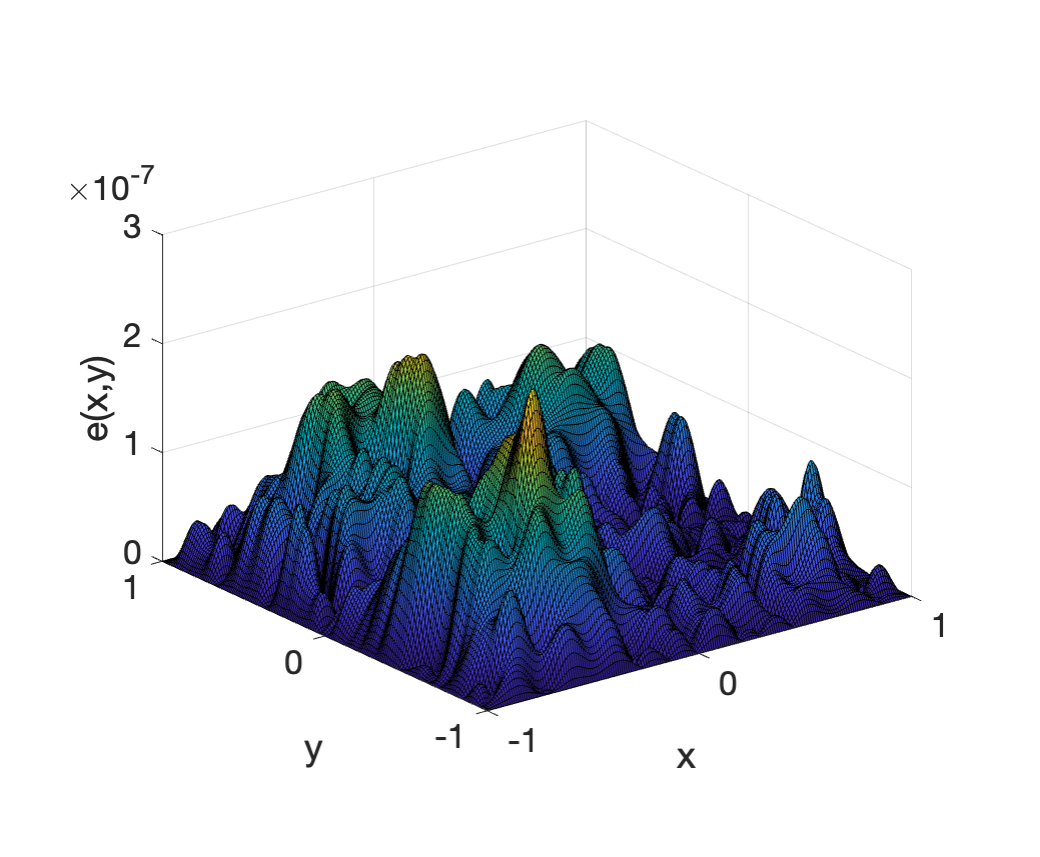}
}

         \caption{ \textbf{(a)-(c)}: Results for Example \ref{exm:1p14} a) with two-scale NN:  $(x^2-1)(y^2-1)N(x,x/\sqrt{\epsilon}, 1/\sqrt{\epsilon},y,y/\sqrt{\epsilon},1\sqrt{\epsilon})$, when $k=4, a_1=1, a_2=4$, with $\alpha=1, \alpha_1=0$, $N_c=10000$. \\
         \textbf{(d)-(g)}: MLNN with NN:  $(x^2-1)(y^2-1)N(x, y)$ ; \\
         \textbf{(h):} MLNN with Fourier feature in \eqref{fourier_2d}. 
         Hyper-parameters are specified in Table \ref{table:multilvNN_hyper}.   }
\label{fig:helm_k14_multilvNN_errors}
\end{figure}

\begin{figure} [!ht]
\centering
\subfigure[exact solution]
{
\includegraphics[width=0.31\textwidth]{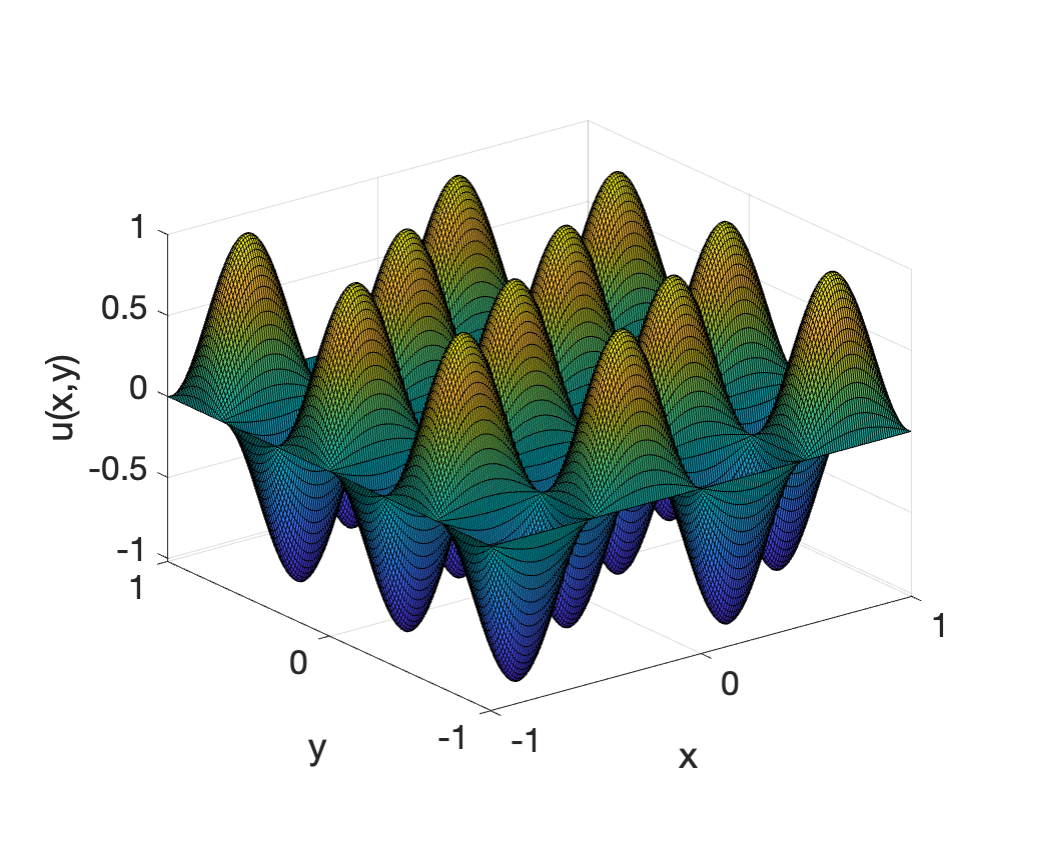}
}
\subfigure[NN solution]
{
\includegraphics[width=0.31\textwidth]{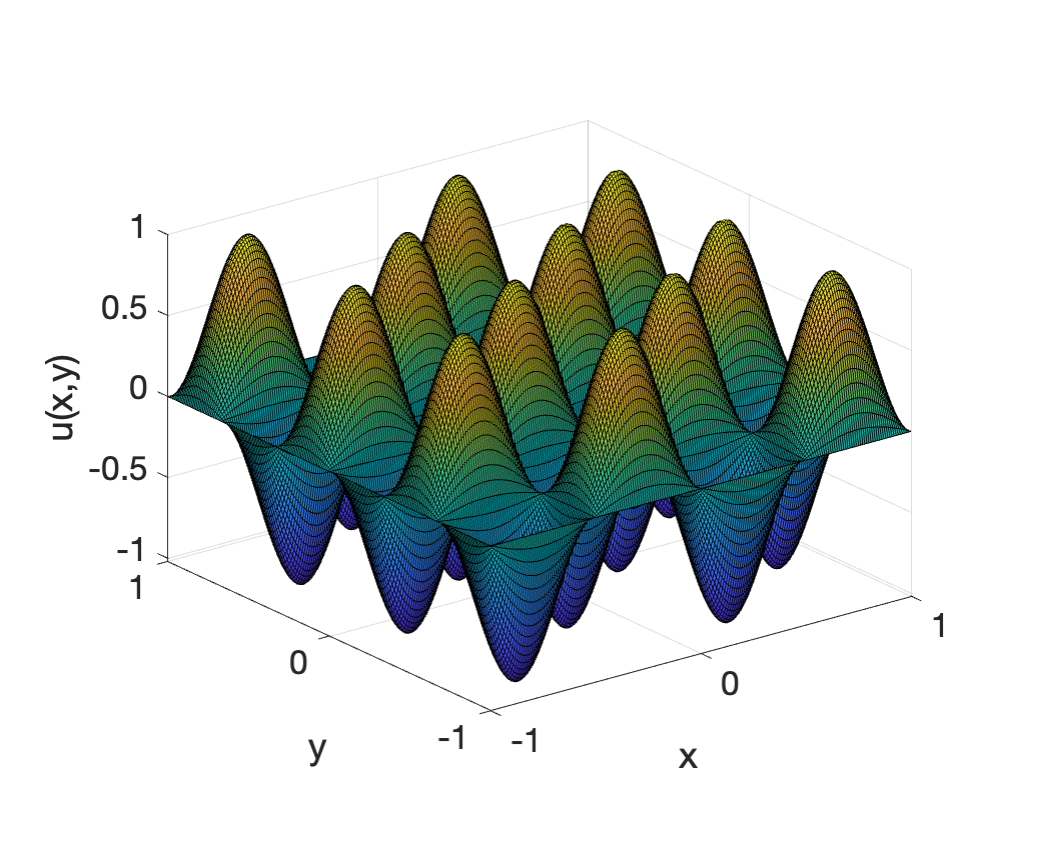}}
\subfigure[absolute error]
{
\includegraphics[width=0.31\textwidth]{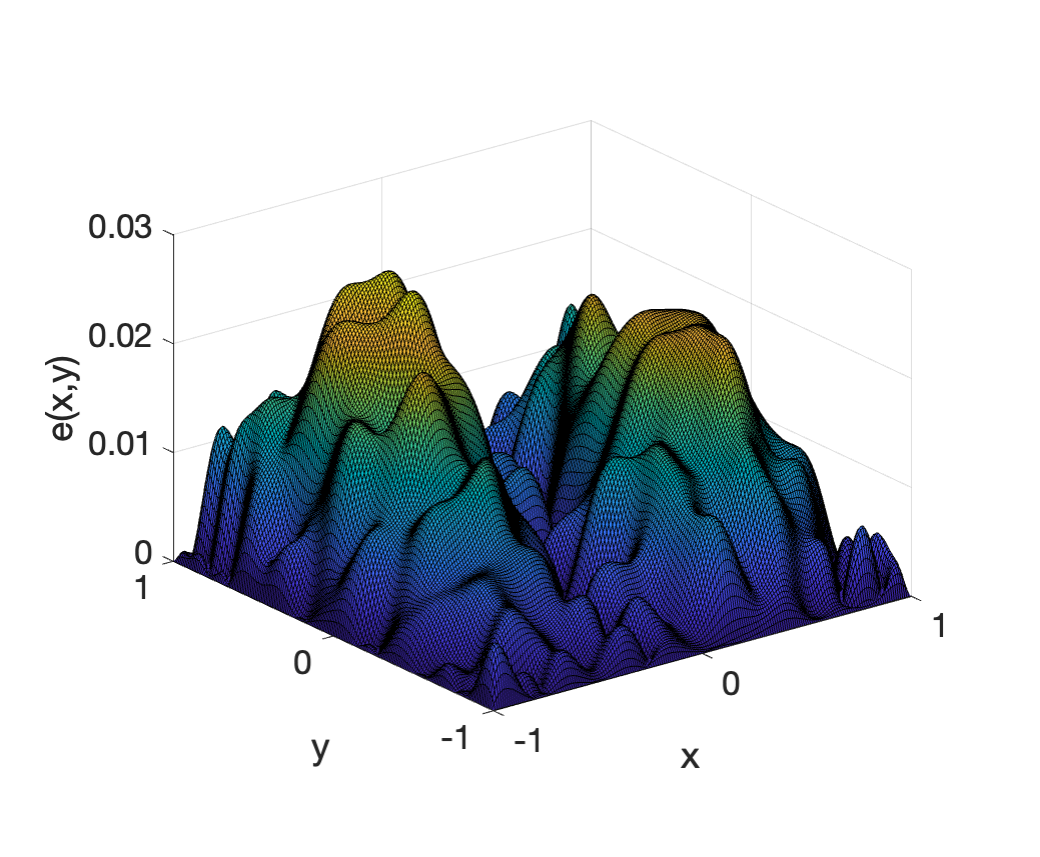}}

\subfigure[relative error]
 { 
 \includegraphics[width=0.31\textwidth]{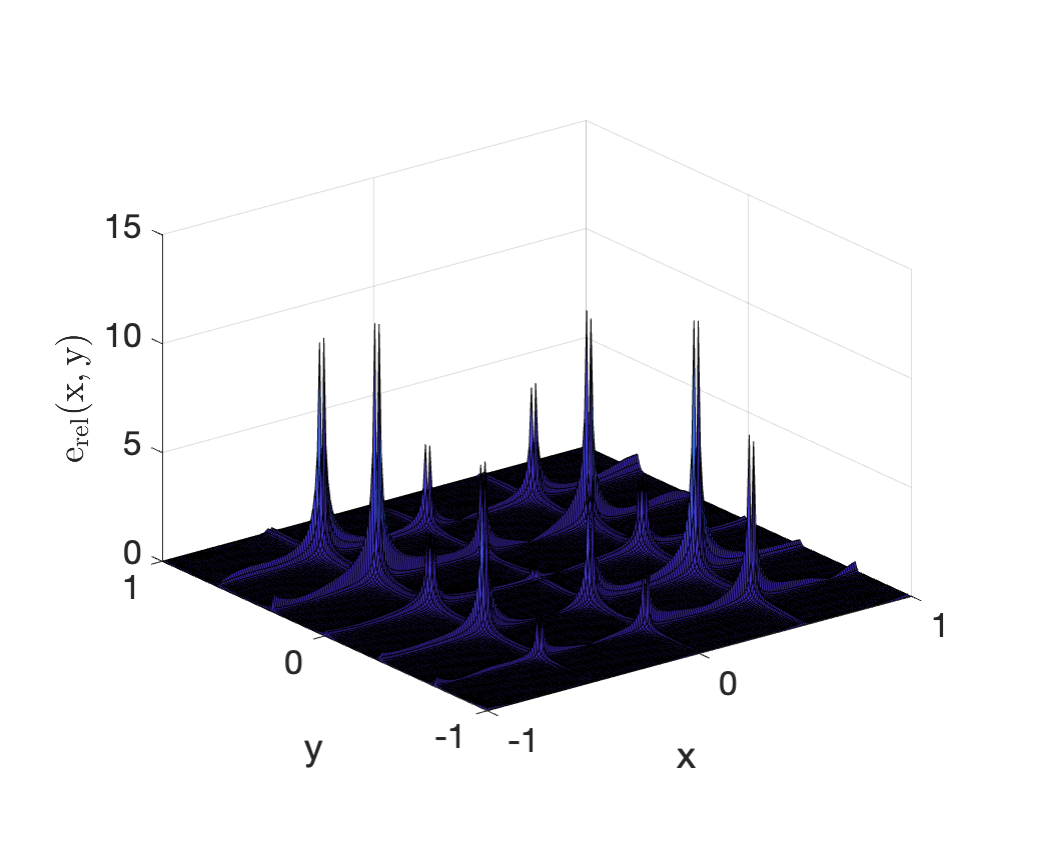}}
\subfigure[loss history]
 { 
 \includegraphics[width=0.26\textwidth]{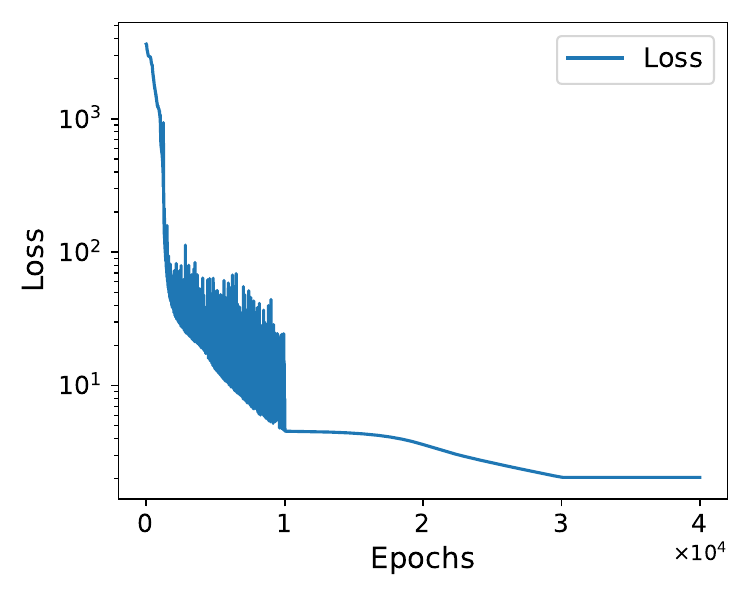}}
\subfigure[errors history]
 { 
 \includegraphics[width=0.26\textwidth]{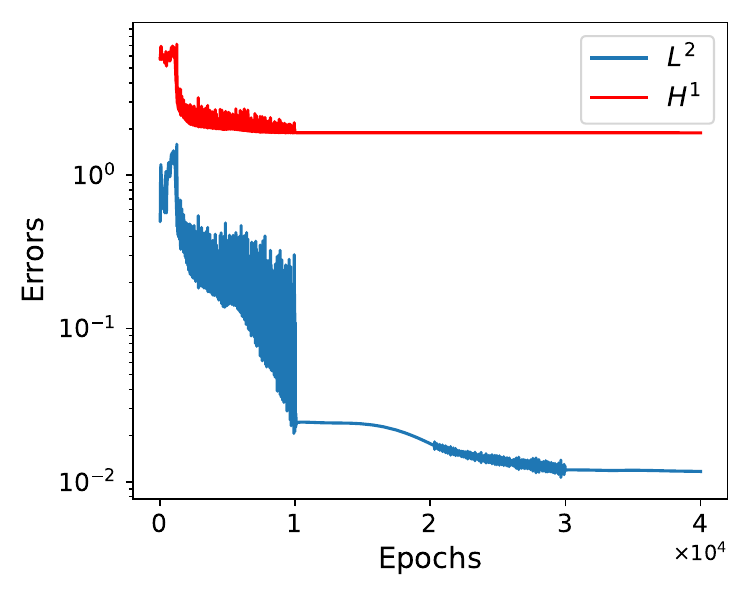}}
 \caption{Numerical results for Example \ref{exm:1p14} a),  when $k=3, a_1=2, a_2=3$ using $(x^2-1)(y^2-1)N(x, x/\sqrt{\epsilon}, 1/\sqrt{\epsilon}, y, y/\sqrt{\epsilon}, 1/\sqrt{\epsilon})$,  $\alpha=1, \alpha_1=0, N_c=22500,  $ epochs=40000. We enforce $e_{rel}(x,y)=0$ on the nodes where the exact solution $u=0$. The relative $l^2$ error: $e_{l^2}^{rel} =5.5318\times 10^{-4}.$
 }
 \label{fig:exm1p14_k3a23} 
 \end{figure}


For \textbf{\textbf{Example \ref{exm:1p14} b)}}, to avoid evaluating singular derivatives of the solution, we exclude the sampling points close to the origin. To utilize NN solutions that can exactly enforce the Dirichlet boundary condition for the problem, we consider the following auxiliary problem:

\begin{align}\label{eq:helm_Feng_auxprob}
-\Delta \tilde{u} - k^2 \tilde{u} &= f \quad \text{in}~\Omega=(-1,1)^{2}, \\
\tilde{u} &= 0 \quad \text{on}~\partial \Omega. \nonumber
\end{align}
Then $u=\tilde{u}+g_b$, where $g_b$ is defined to be
\begin{equation*}
g_b(x,y)=x^2 (g_2(y)-g(1,1))+y^2 (g_1(x)-g(1,1)) +g(1,1),
\end{equation*}
with $g_2(y)=g(1,y)$, $g_1(x)=g(x,1)$. Taking $\epsilon=1/k^2$, we  utilize the neural networks \((x^2-1)(y^2-1)N(x, x/{\sqrt{\epsilon}}, y, 1/{\sqrt{\epsilon}}, y/{\sqrt{\epsilon}}, 1/{\sqrt{\epsilon}})\) to solve the auxiliary problem \eqref{eq:helm_Feng_auxprob} and, consequently, obtain the solution \(u\) to the original Helmholtz problem.
 We collect the numerical results in Figure \ref{fig:helmFeng}. The figure demonstrates that the two-scale NN solution accurately captures the radial oscillation pattern for a modest wave number $k=10$, an achievement the corresponding one-scale approach does not replicate under the same configuration. 

\begin{figure}[!ht]
	\centering
          \subfigure[exact solution]{
	     \includegraphics[width=0.30\textwidth]{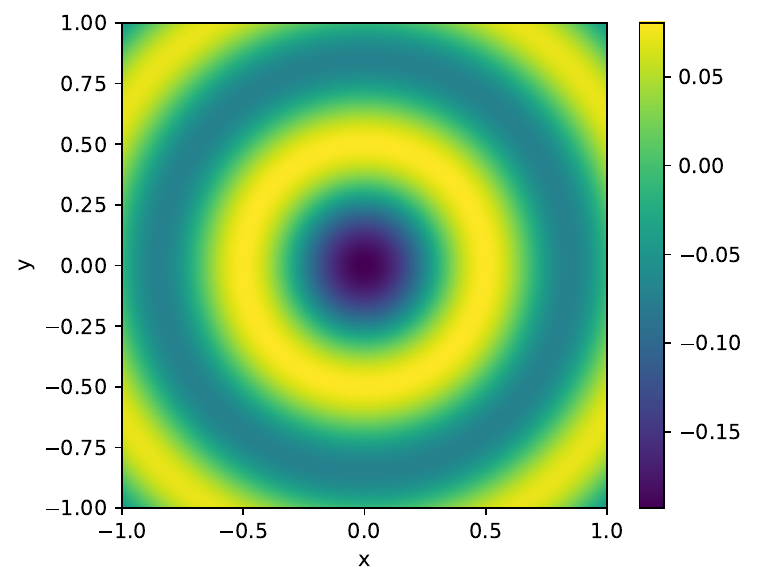}}
        \subfigure[NN solution]{
        \includegraphics[width=0.30\textwidth]{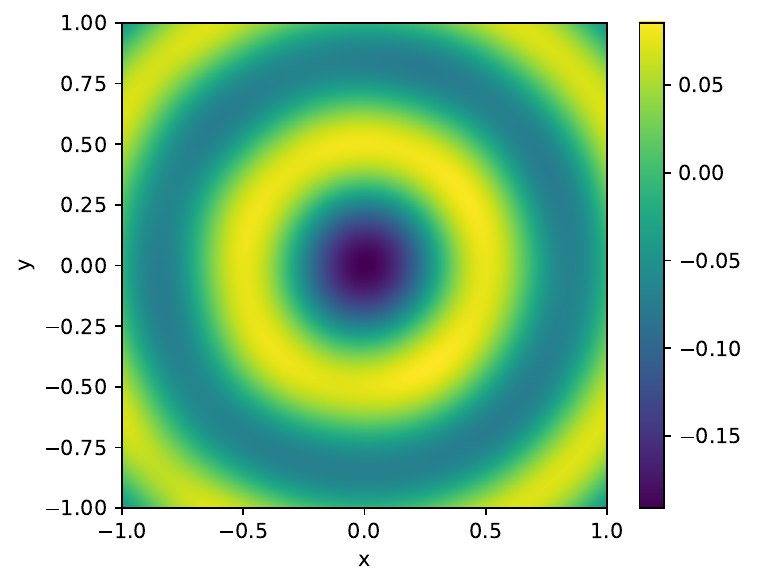}}
        
        \subfigure[absolute error]{
        \includegraphics[width=0.32\textwidth]{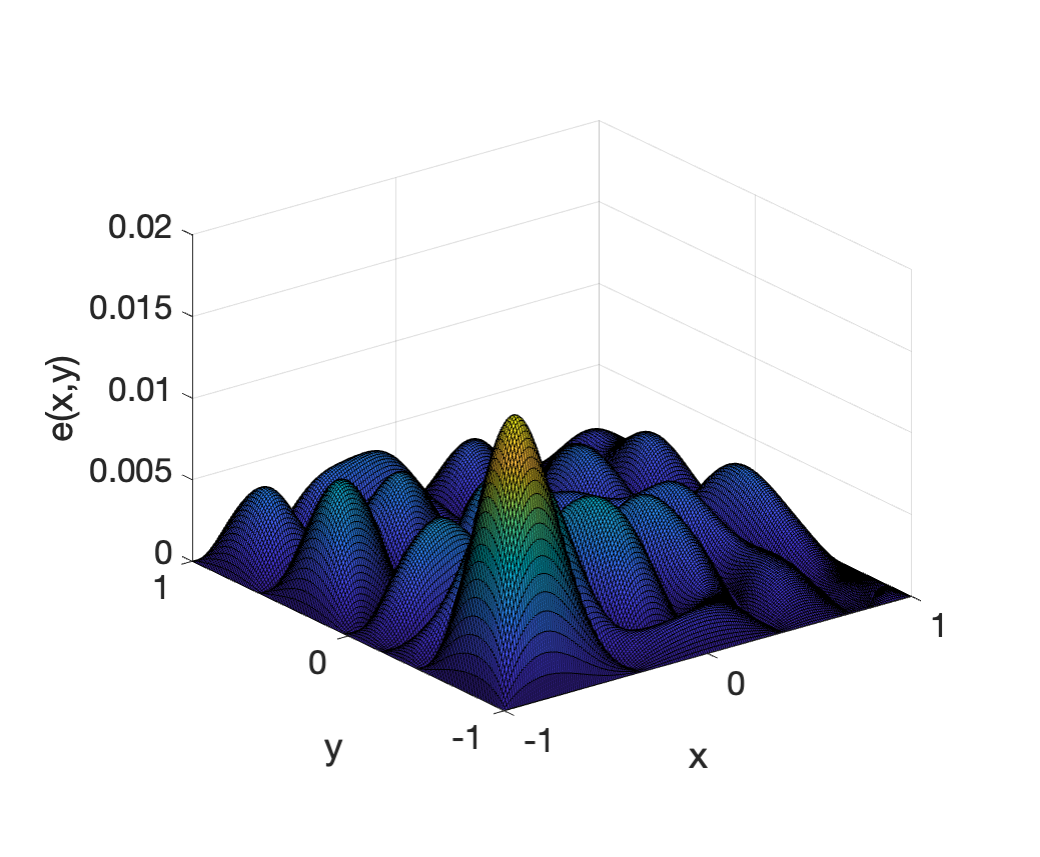}}
        \subfigure[relative error]{
        \includegraphics[width=0.32\textwidth]{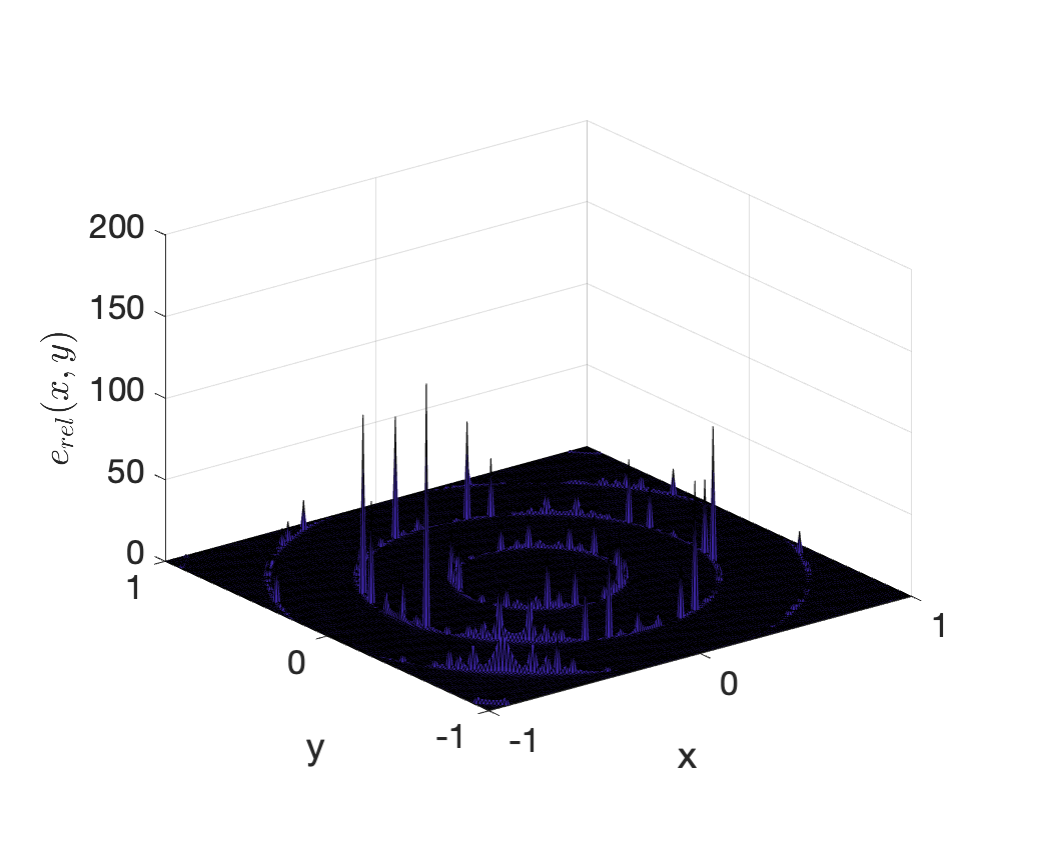}}	
        \subfigure[loss history]{
        \includegraphics[width=0.3\textwidth]{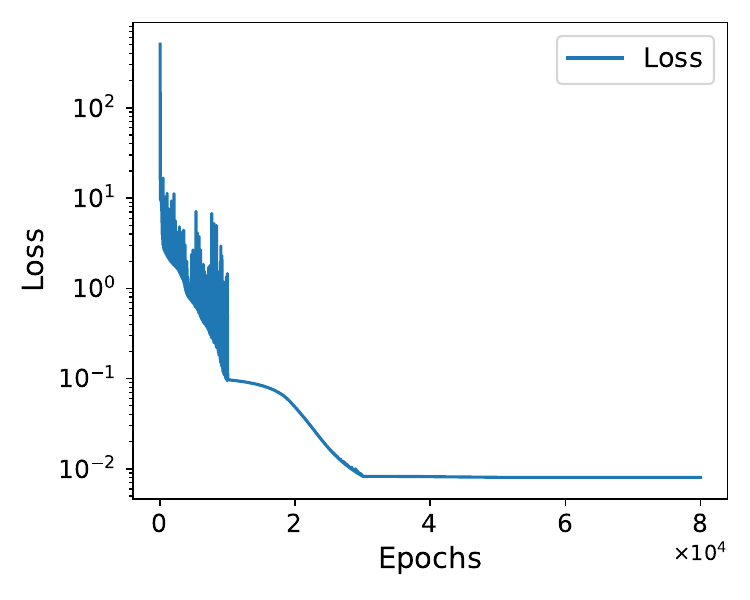}
        }

  \caption{ Numerical results for Example \ref{exm:1p14} b) when $k=10$  using  $(x^2-1)(y^2-1)N(x, x/{\sqrt{\epsilon}}, y, 1/{\sqrt{\epsilon}}, y/{\sqrt{\epsilon}}, 1/{\sqrt{\epsilon}})$  with $k=1/{\sqrt{\epsilon}}$, $\alpha=1$, $\alpha_1=0, N_c=22497,$ epochs=80000. The relative $l^2$ error is: { $e_{l^2}^{rel}=3.41\times 10^{-3}$}.}
\label{fig:helmFeng}
\end{figure}

{
\begin{exm}[1D ODE with one boundary layer, comparison with  multi-level neural networks]\label{exm:burgers1d-steady}
\begin{equation*}
	-\varepsilon u'' + u'  =1, \quad   \quad 
	u(0)=0, \quad u (1)=0. 
	\end{equation*}	
\end{exm}

The example is the same as Example 5.2 in \cite{multilvNN23}. 
We compare the results of solving this problem using the two-scale NN method and MLNN. To generate the MLNN results, we replicate the neural network configuration with the Fourier feature as described in Example 5.2 of \cite{multilvNN23}:
\begin{equation}\label{fourier_1d}
\sin \left(\boldsymbol{\omega_M} \pi (1-x)\right)N(\gamma(x)),
\end{equation}
where \(\gamma(x) = [\cos(\boldsymbol{\omega_M} x), \sin(\boldsymbol{\omega_M} x)]\), ${\boldsymbol{\omega_M}}=(\omega_1, ..., \omega_M) $. We denote the absolute error between the exact solution and the $i$-th level approximated solution as $e_i=\abs{u-\sum_{k=0}^i \left(\tilde{u}_k/\prod_{j=0}^{k} \mu_j \right)}$. 

The MLNN results are obtained with the hyper-parameters listed in Table \ref{table:MLNNexm5p2} and presented in Figures \ref{fig:MLNN_1e_2} and \ref{fig:MLNN_1e_3}. 
For the case $\epsilon=10^{-2}$, as shown in Figure \ref{fig:MLNN_1e_2}, even though we use only one level of correction, the accuracy of the approximation achieves $10^{-7}$ in the absolute error. The loss and error histories demonstrate fast convergence and are consistent with Figure 17 in \cite{multilvNN23}. The accuracy can be further refined to $10^{-9}$ if three levels of corrections are applied. 
However, for the case $\epsilon=10^{-3}$, as shown in Figure \ref{fig:MLNN_1e_3}, the MLNN results have not captured the exact solution, despite the employment of three levels of corrections.

As a comparison, we present the results using the two-scale NN method in Figure \ref{fig:exm5p2_2s}, with the NN size of $(3, 20, 20, 20, 20, 1)$ and hyper-parameters in Table \ref{tab:lr_pt}. 
When $\epsilon=10^{-2}$, as illustrated in Figure \ref{fig:exm5p2_2s} (a) to (d), the absolute error reaches the magnitude of $10^{-4}$ after 35000 iterations. This is less precise compared to the MLNN results, where an accuracy of the order of $10^{-7}$ is achieved with just one level of correction.
When $\epsilon=10^{-3}$, we solve the equation directly without employing the successive training strategy in Algorithm \ref{alg:succesive-training-two-scale}. As depicted in Figure \ref{fig:exm5p2_2s} (e) to (h), the two-scale NN method successfully captures the key features of the exact solution. The overall accuracy is similar to that shown in Figure \ref{fig:results_exm1p1_1e-3} for Example \ref{exm:1p1} (a similar 1D problem with one boundary layer).

In summary, this case study demonstrates that although the two-scale NN method may not necessarily surpass the MLNN method in terms of accuracy, it offers a simple approach to solving PDEs with small parameters. Furthermore, as the parameter in the PDE becomes smaller, the two-scale NN method can still provide reasonable accuracy with only a little or no special treatment.

}
{
\begin{table}[!htb]

\centering
\label{tab:hyperparameters}
\begin{tabular}{l|c|c|c|c}
\hline
Hyper-parameters & $\tilde{u}_0$ & $\tilde{u}_1$ & $\tilde{u}_2$ & $\tilde{u}_3$ \\
\hline
\# Hidden layers  & 1 & 1 & 1 & 1 \\
Width of each hidden layer & 5 & 10 & 20 & 20 \\
\# Adam iterations & 4000 & 4000 & 4000 & 10000 \\
\# L-BFGS iterations & 500 & 1000 & 2000 & 0 \\
\# Wave numbers $M$ & 3 & 5 & 7 & 3 \\
\hline
\end{tabular}

\caption{ Hyper-parameters used in Example \ref{exm:burgers1d-steady}  with multi-level neural networks \cite{multilvNN23}}

\label{table:MLNNexm5p2}

\end{table}

\begin{figure} [!ht]
\centering
\subfigure[absolute error $e_1$]
{
\includegraphics[width=0.3\textwidth]{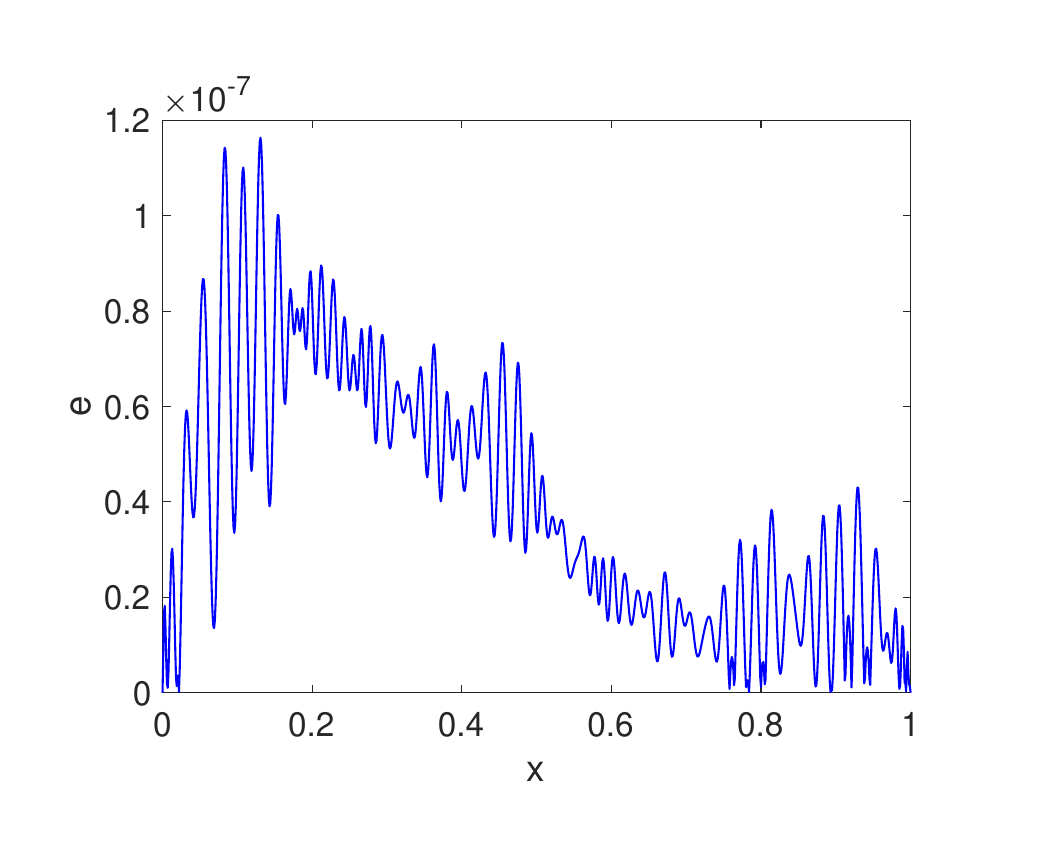}
}
\subfigure[loss history]
{
\includegraphics[width=0.28\textwidth]{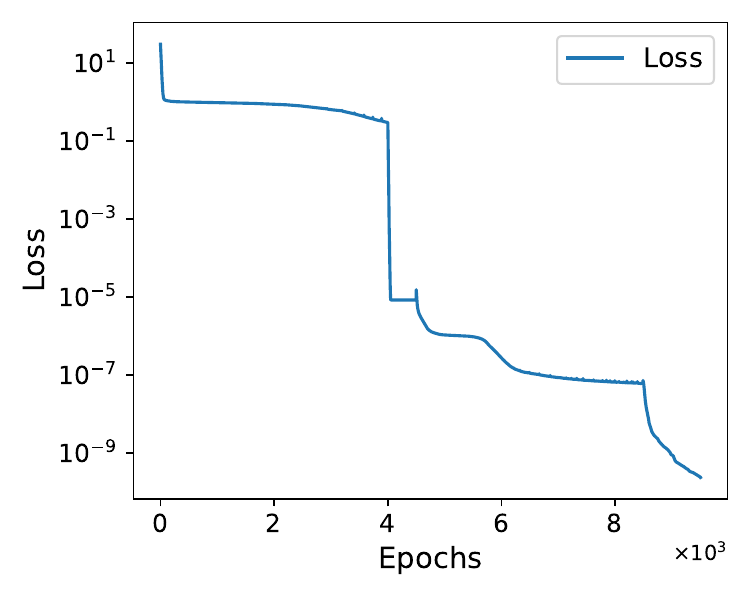}}
\subfigure[errors history]
{
\includegraphics[width=0.28\textwidth]{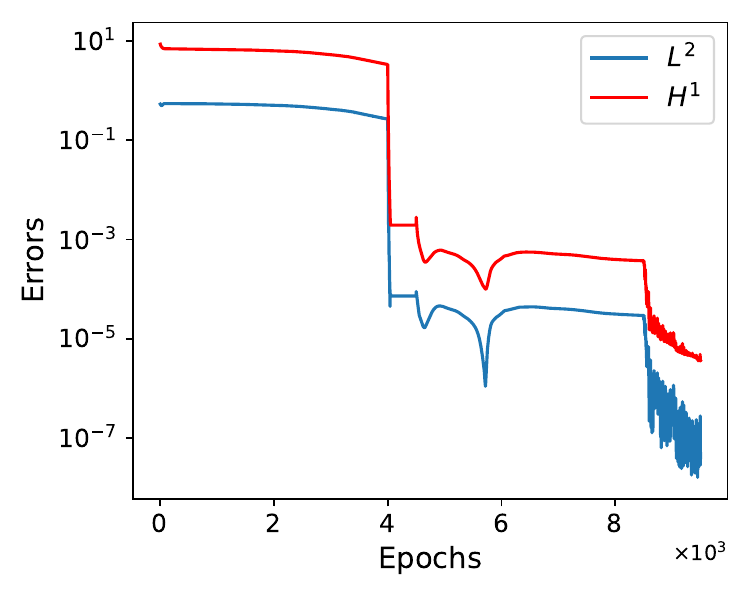}}
 \caption{  Numerical results for Example \ref{exm:burgers1d-steady} when $\epsilon=10^{-2}$, using MLNN adding the Fourier feature in \eqref{fourier_1d}, with one level of correction. }
 
 \label{fig:MLNN_1e_2} 
 \end{figure}

}

\begin{figure} [!ht]
\centering
\subfigure[$\tilde{u}_0$]
{
\includegraphics[width=0.23\textwidth]{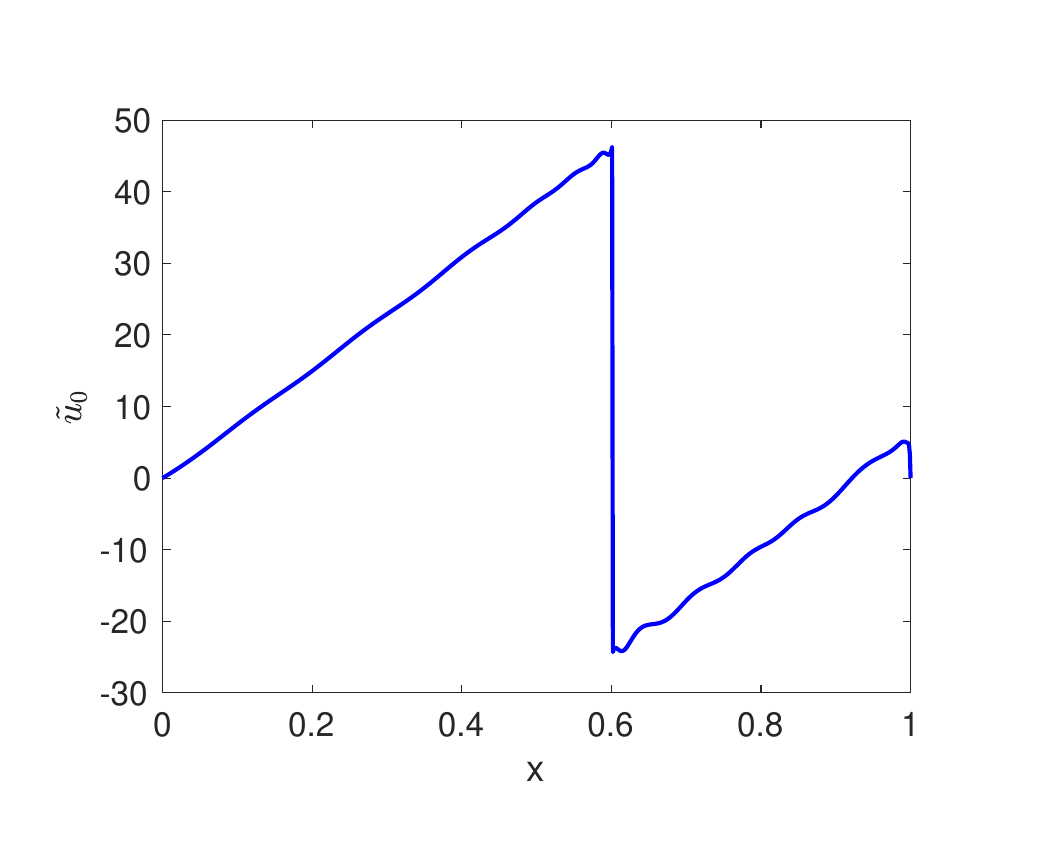}
}
\subfigure[$\tilde{u}_1$]
{
\includegraphics[width=0.23\textwidth]{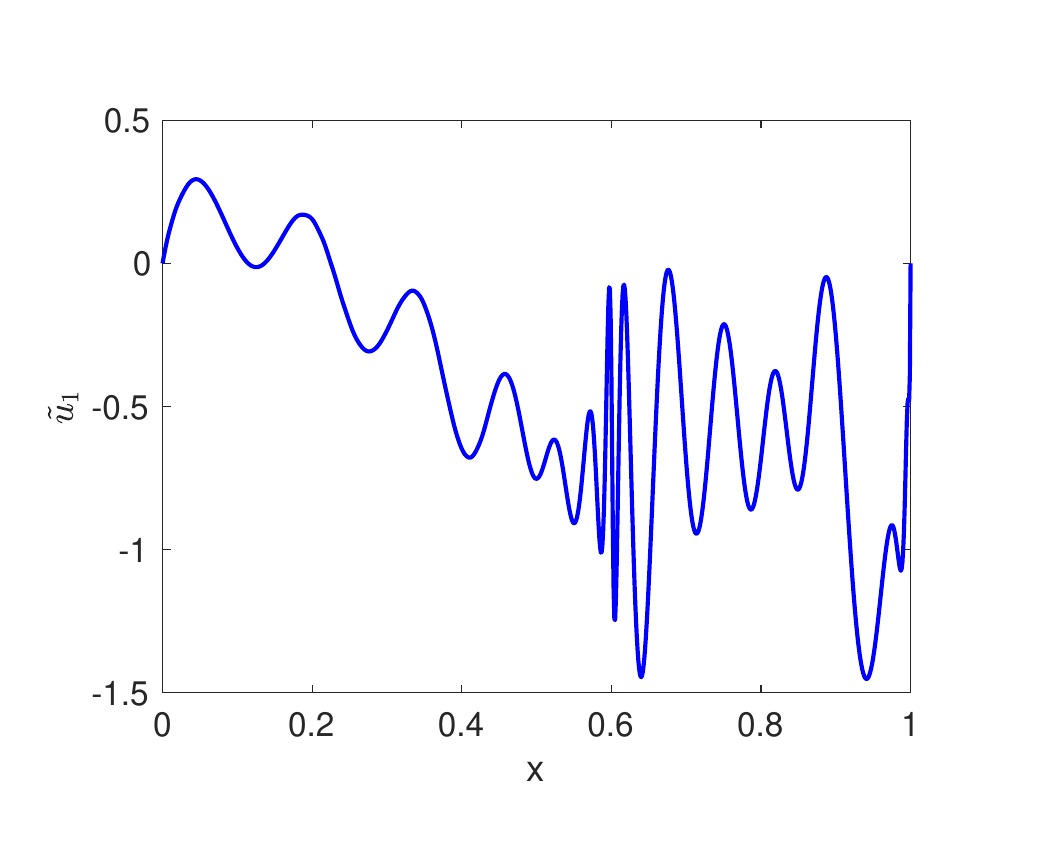}}
\subfigure[$\tilde{u}_2$]
{
\includegraphics[width=0.23\textwidth]{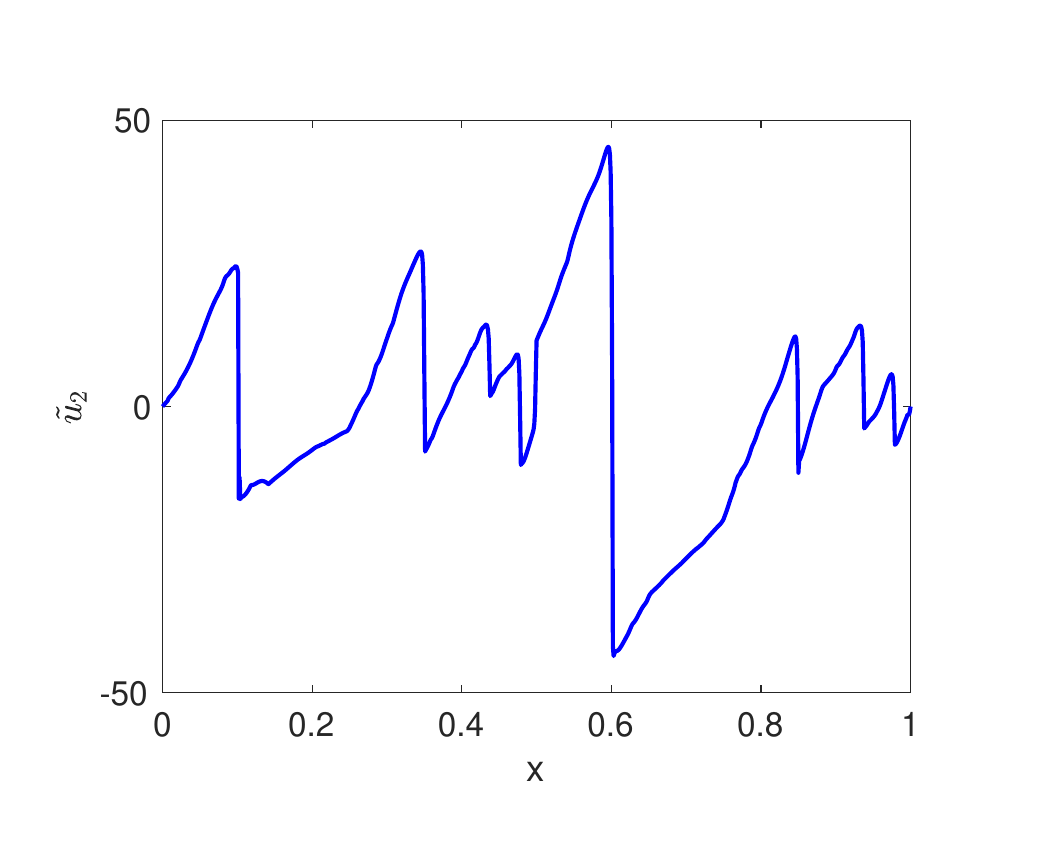}}
\subfigure[$\tilde{u}_3$]
{
\includegraphics[width=0.23\textwidth]{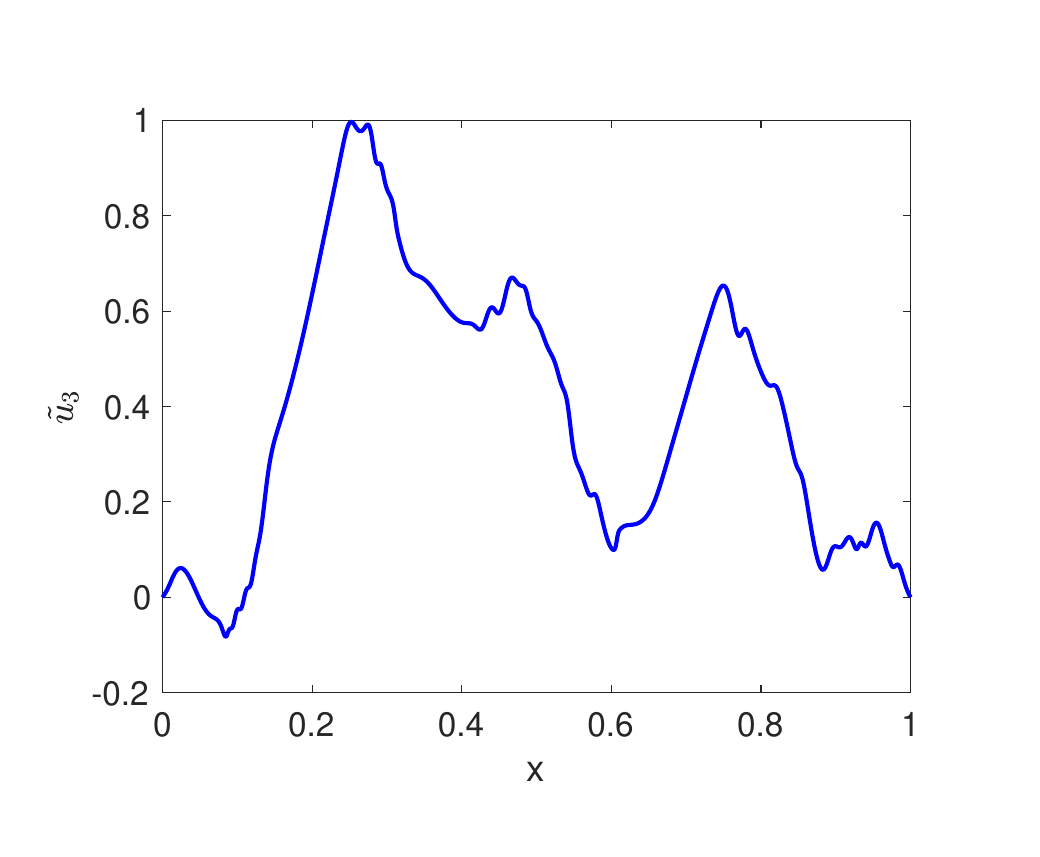}}

\subfigure[absolute error $e_3$]
{
\includegraphics[width=0.23\textwidth]{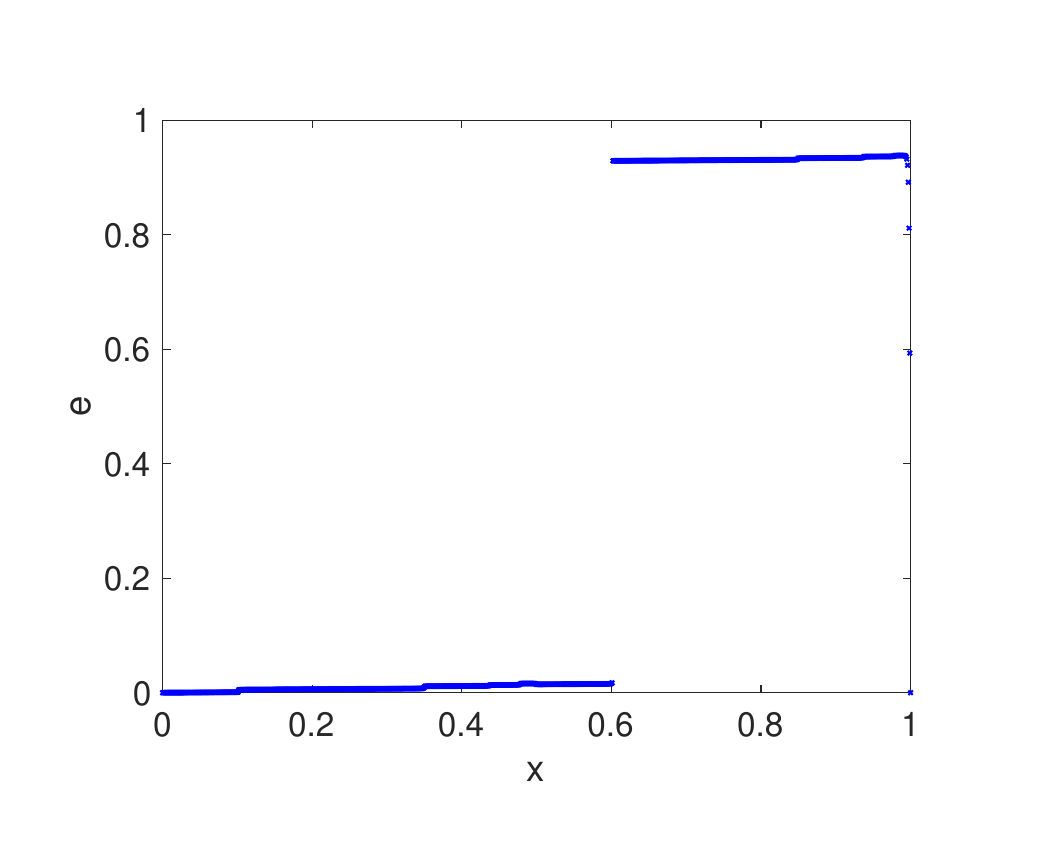}}
\subfigure[loss history]
{
\includegraphics[width=0.23\textwidth]{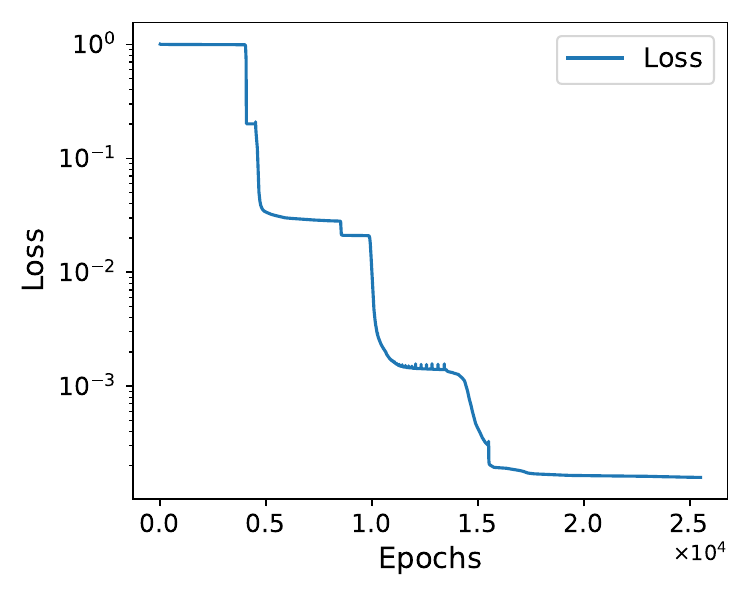}}
\subfigure[errors history]
{
\includegraphics[width=0.23\textwidth]{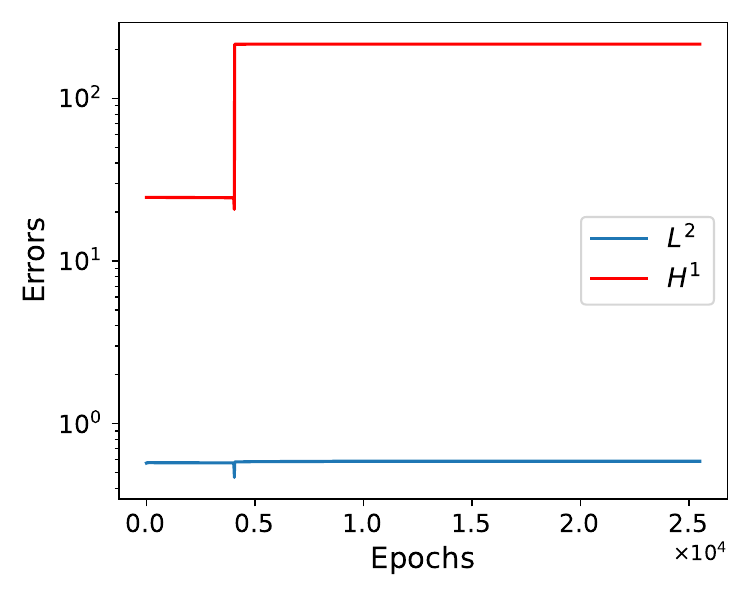}}
 \caption{ Numerical results for Example \ref{exm:burgers1d-steady}  when $\epsilon=10^{-3} $, using MLNN adding the Fourier feature in \eqref{fourier_1d}, with three levels of correction. }
 
 \label{fig:MLNN_1e_3} 
 \end{figure}

\begin{figure} [!ht]
\centering
\subfigure[exact and NN solutions for $\epsilon=10^{-2}$]
{
\includegraphics[width=0.23\textwidth]{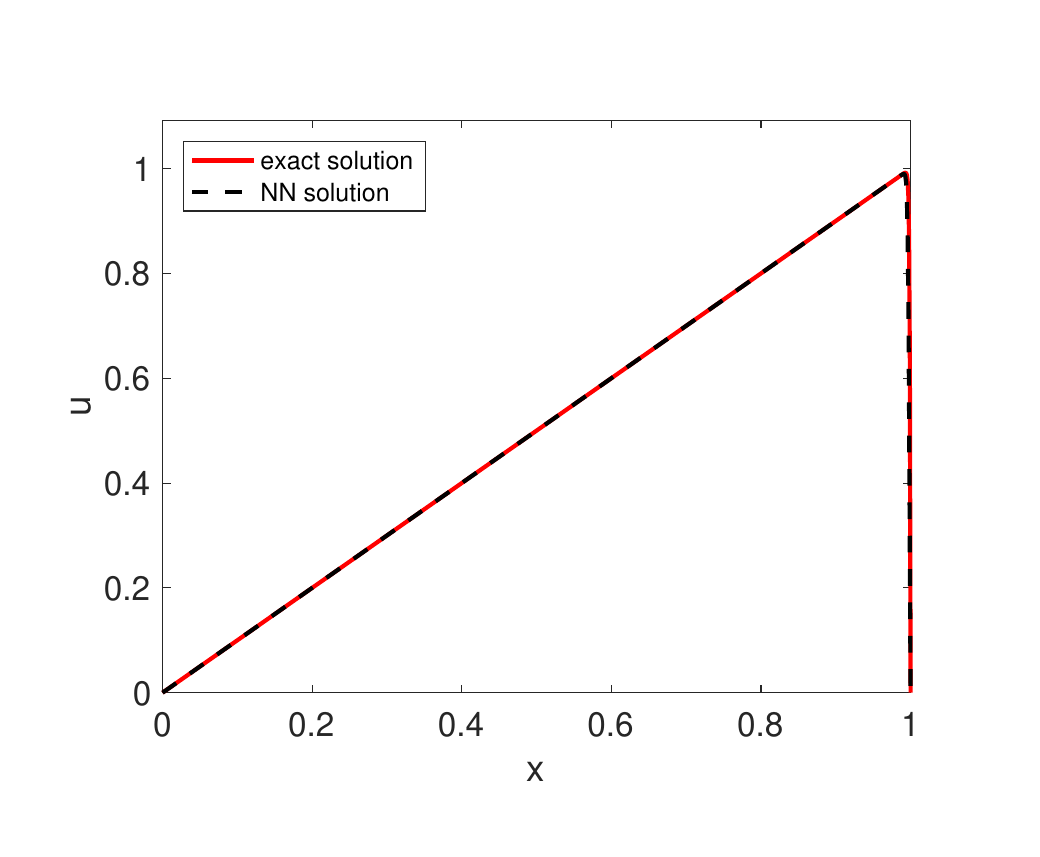}}
\subfigure[absolute error]
 { 
\includegraphics[width=0.23\textwidth]{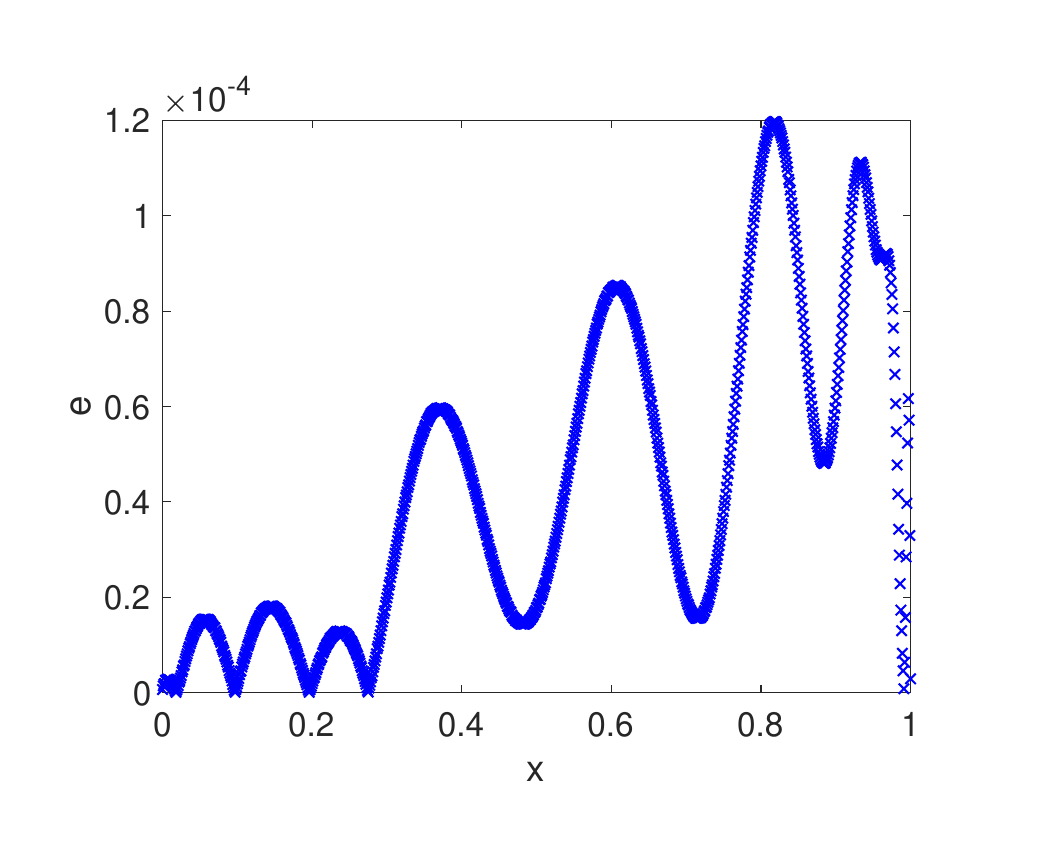}}
 \subfigure[loss history]
{
\includegraphics[width=0.22\textwidth]{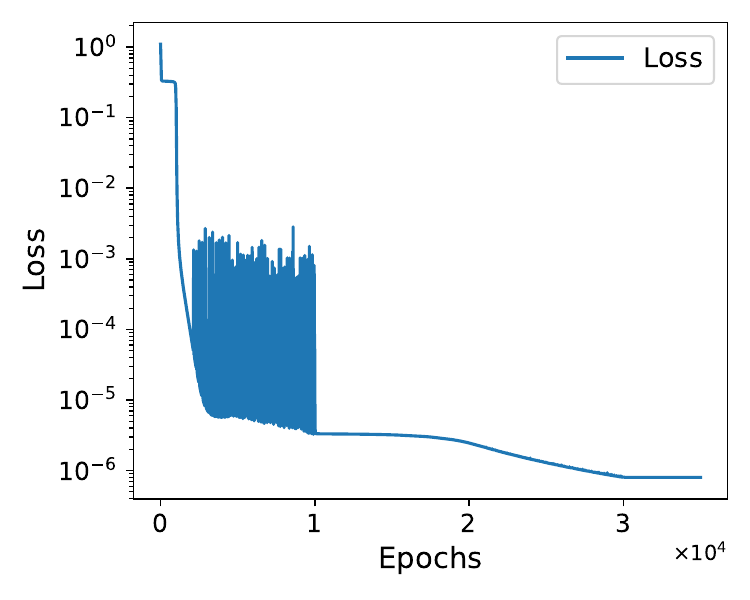}}
\subfigure[errors history]
 { 
\includegraphics[width=0.22\textwidth]{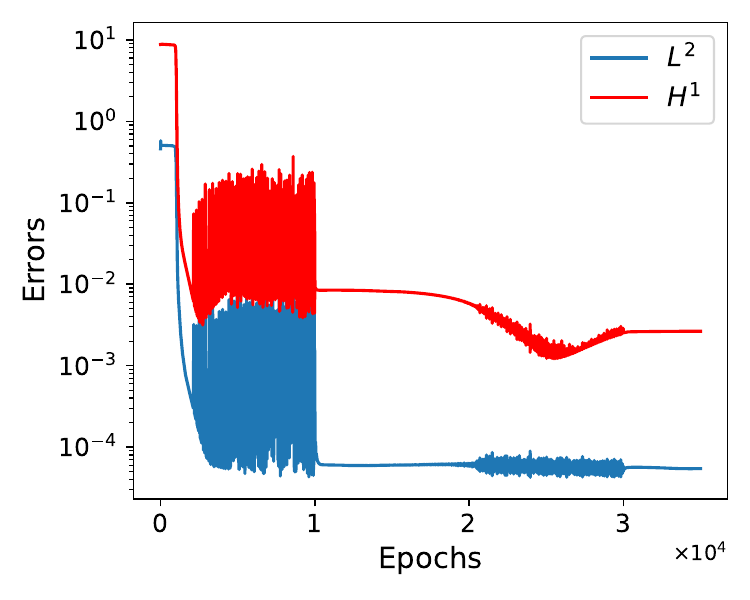}}

\subfigure[exact and NN solutions for $\epsilon=10^{-3}$]
{
\includegraphics[width=0.23\textwidth]{figs/exm7/uexactNN_exm5p2_2sc_1e-3.pdf}}
\subfigure[absolute error]
 { 
\includegraphics[width=0.23\textwidth]{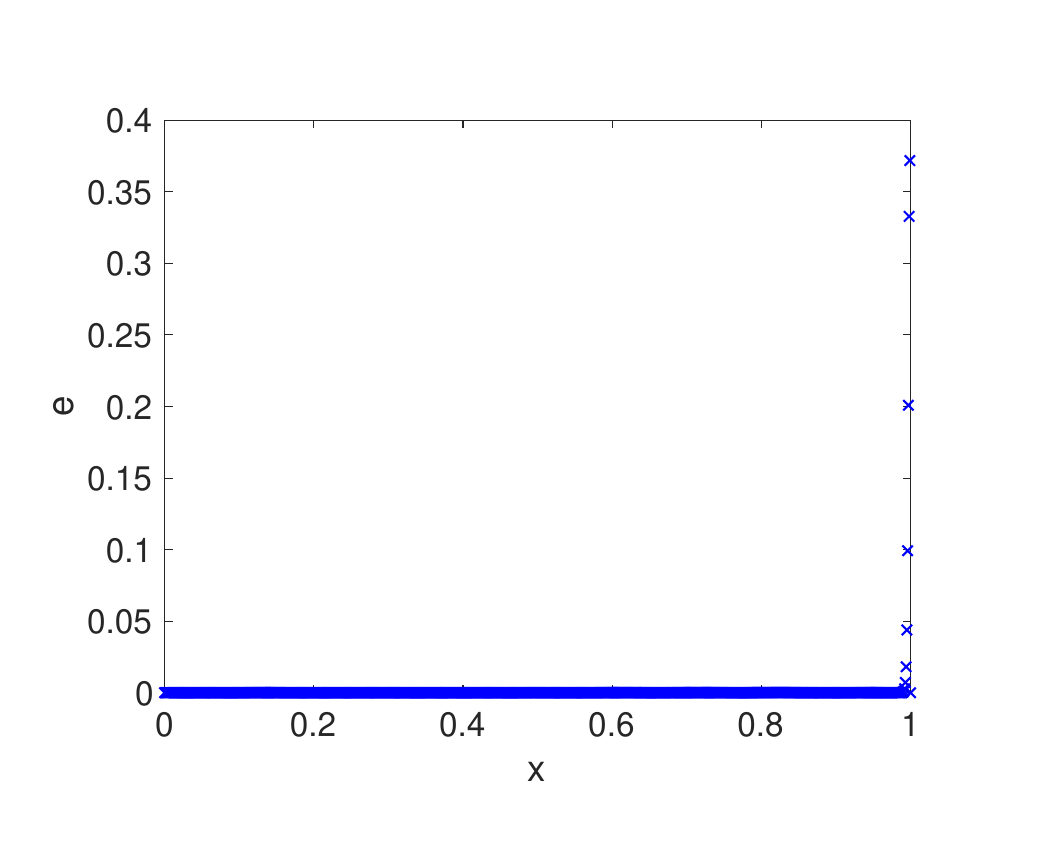}}
 \subfigure[loss history]
{
\includegraphics[width=0.22\textwidth]{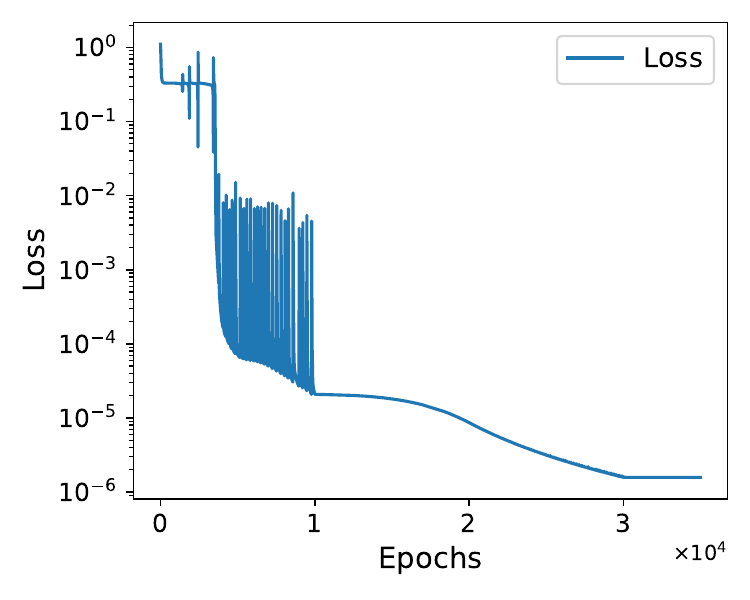}}
\subfigure[errors history]
 { 
\includegraphics[width=0.22\textwidth]{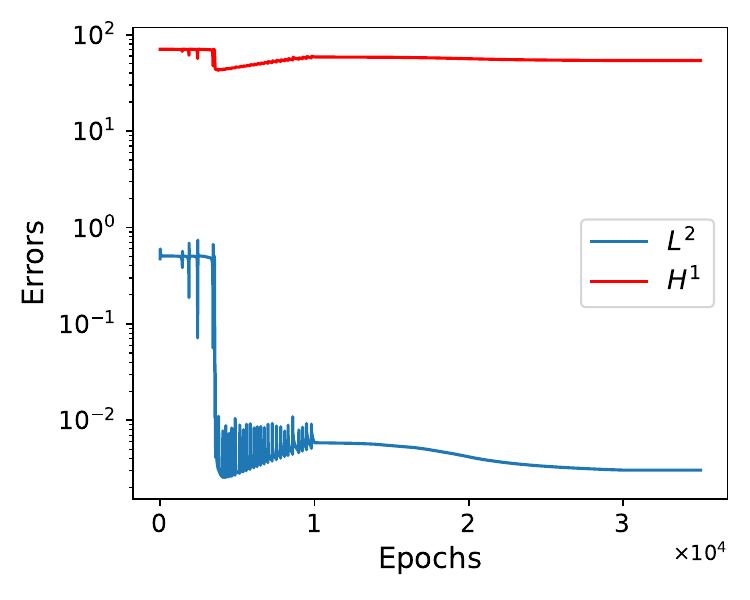}}
 
 \caption{ Numerical results for Example \ref{exm:burgers1d-steady} using {$N(x, (x-0.5)/\sqrt{\epsilon}, 1/\sqrt{\epsilon})$}:  \\
 \textbf{(a)-(d)}:  when $\epsilon=10^{-2}$, with $\alpha=1$, $\alpha_1=0$, $N_c=200$, epochs=35000; \\
 \textbf{(e)-(h)}: when $\epsilon=10^{-3} $, with $\alpha=1$, $\alpha_1=0$, $N_c=450$, epochs=35000.}
 \label{fig:exm5p2_2s} %
 \end{figure}


{
\begin{exm}[1D steady-state Burgers' equation]\label{exm:1d_nonlinear}
\begin{equation*}
	 u'' -ku u'  =0, \quad 0<x<1,  \quad 
	u(0)=-1, \quad u (1)=-2/(k+2). 
	\end{equation*}	
\end{exm}
This example, a nonlinear 1D ODE with a boundary layer, is the same as Example 7 in \cite{multilvNN23}. We solve this problem with the two-scale NN method with the NN size of $(3,20, 20, 20, 20, 1)$ and discuss the results alongside those presented in Example 7 of \cite{multilvNN23} using MLNN. Here, we take $\epsilon=1/k$. 

When $k=8$ ($\epsilon=1/8$, a larger $\epsilon$), as shown in Figure \ref{fig:exm7_2s} (a) to (d), the accuracy of the two-scale NN method reaches $10^{-5}$. In contrast, the MLNN results for this example, as reported in \cite{multilvNN23}, achieve an accuracy of $10^{-13}$ with three levels of correction. Although the two-scale NN results do not reach the accuracy of those provided by MLNN, it is worth mentioning that the two-scale NN results are obtained without adding Fourier features and are optimized with only a first-order optimizer (Adam).

When $\epsilon$ is modestly small, such as $\epsilon=10^{-2}$ ($k=100$), as shown in Figure \ref{fig:exm7_2s} (e) to (h), the two-scale NN method can directly handle it without any special treatment, achieving an accuracy of $10^{-5}$.
As $\epsilon$ decreases further, for example to $\epsilon=10^{-3}$ ($k=1000$), the successive training strategy in Algorithm \ref{alg:succesive-training-two-scale} can be utilized to enhance the training outcomes. 
%
As demonstrated in Figure \ref{fig:exm7_2s_1e_3}, the two-scale NN method achieves an accuracy of $10^{-3}$ for $\epsilon=10^{-3}$ by employing Algorithm \ref{alg:succesive-training-two-scale} with starting $\epsilon_0=10^{-2}$ and $\ell=2$.

Similar to our observations in Example \ref{exm:burgers1d-steady}, this case study shows that while the two-scale NN method may not necessarily exceed the accuracy of the MLNN method, it provides a simple approach to solving PDEs with small parameters. When the parameter in the PDE is modestly small, the two-scale NN method can deliver reasonable accuracy without any special treatment. For significantly smaller values of $\epsilon$, the successive training strategy can be employed to facilitate reasonably accurate outcomes.
}
\begin{figure} [!ht]
\centering
\subfigure[exact and NN solutions for $\epsilon=1/8$]
{
\includegraphics[width=0.23\textwidth]{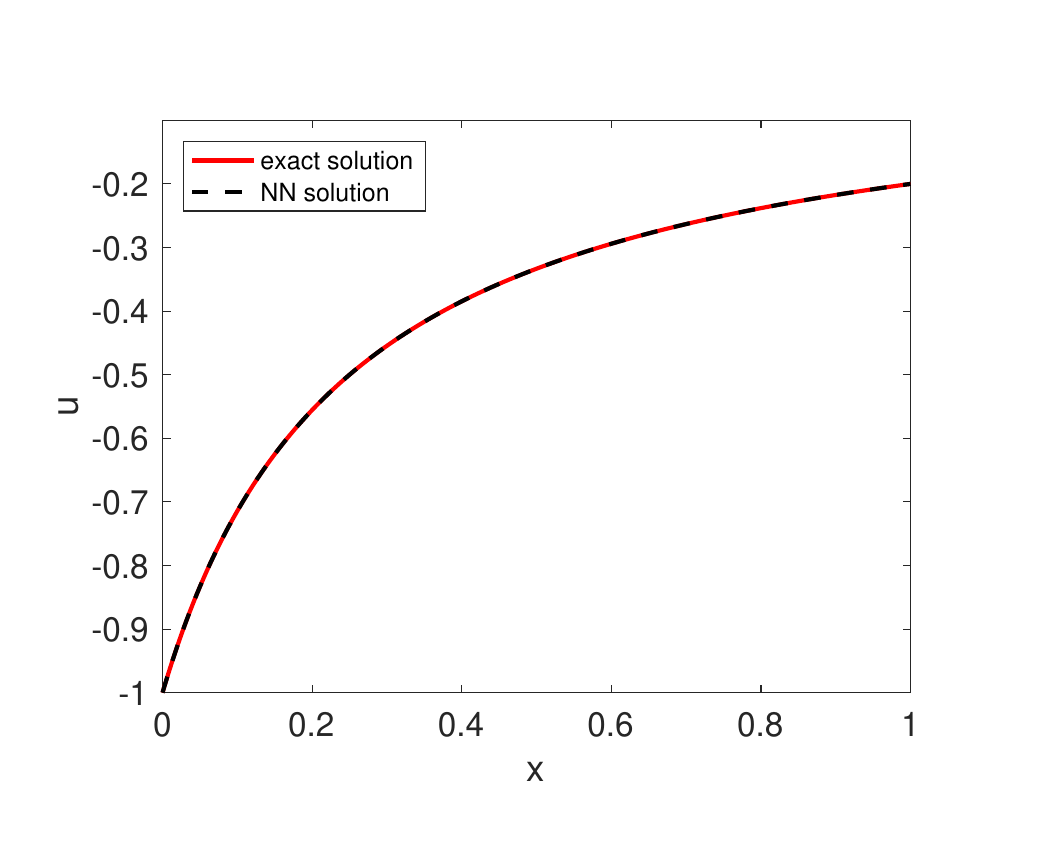}}
\subfigure[absolute error]
 { 
\includegraphics[width=0.23\textwidth]{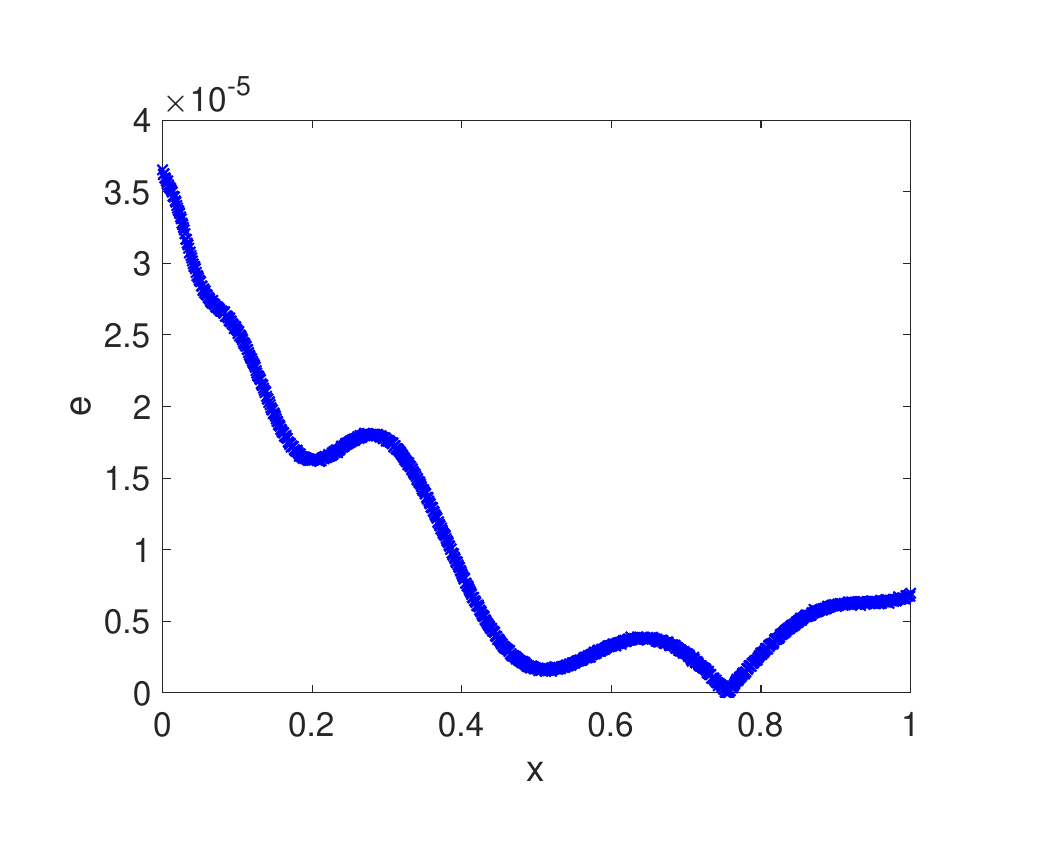}}
 \subfigure[loss history]
{
\includegraphics[width=0.22\textwidth]{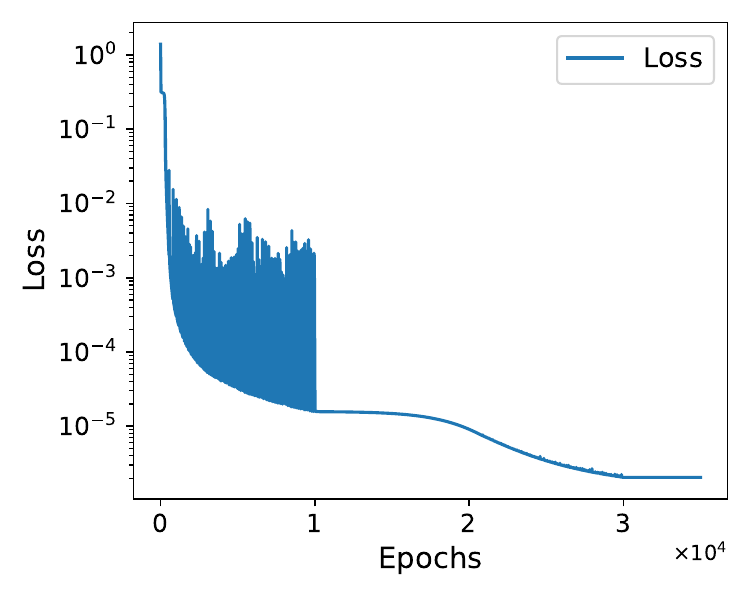}}
\subfigure[errors history]
 { 
\includegraphics[width=0.22\textwidth]{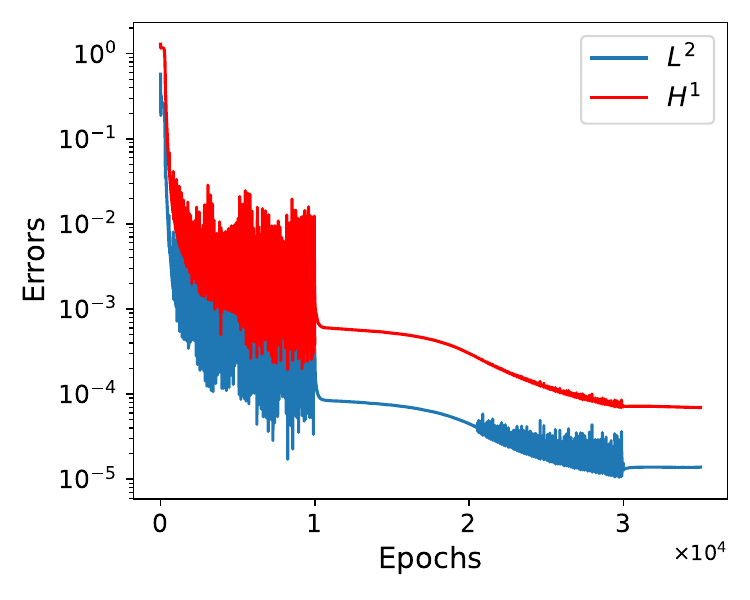}}

\subfigure[exact and NN solutions for $\epsilon=10^{-2}$]
{
\includegraphics[width=0.23\textwidth]{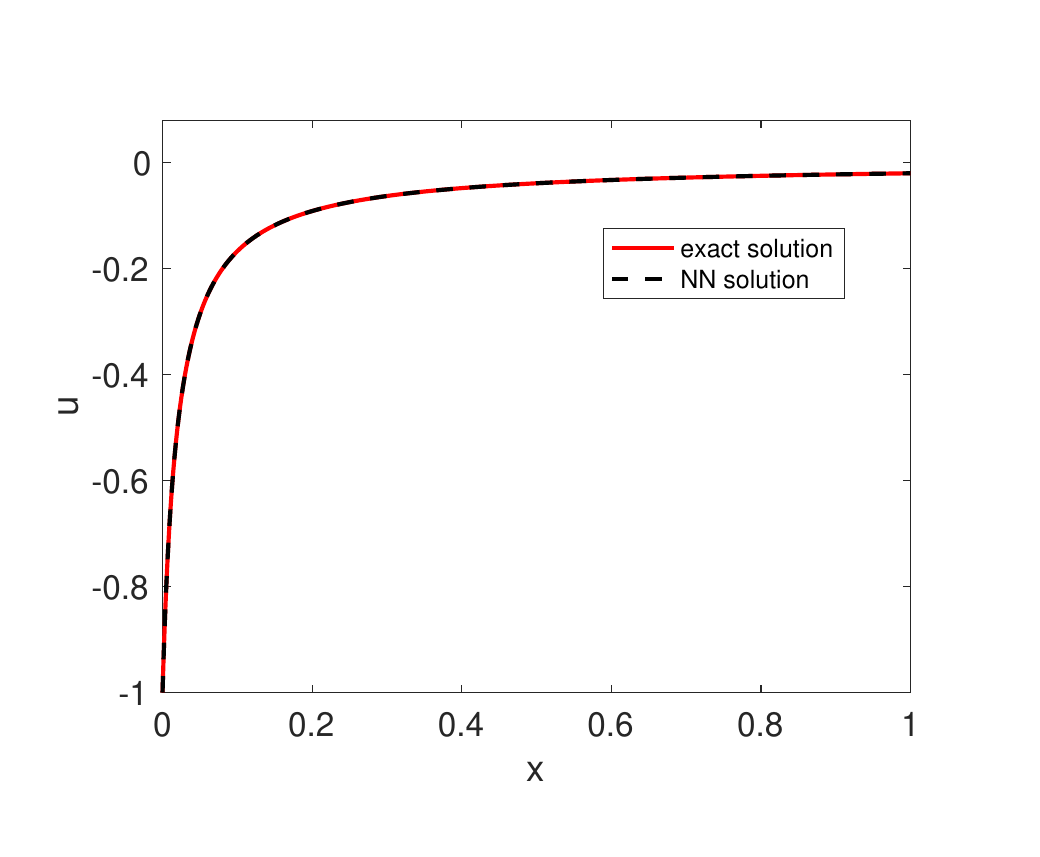}}
\subfigure[absolute error]
 { 
\includegraphics[width=0.23\textwidth]{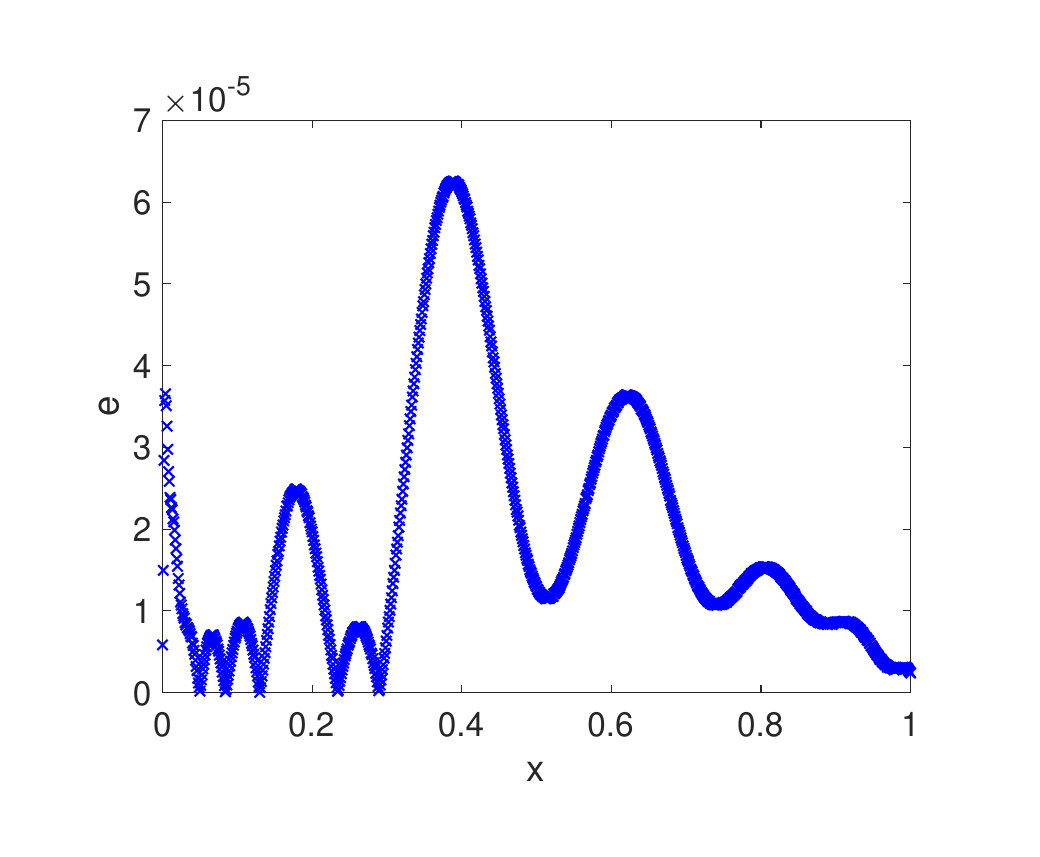}}
 \subfigure[loss history]
{
\includegraphics[width=0.22\textwidth]{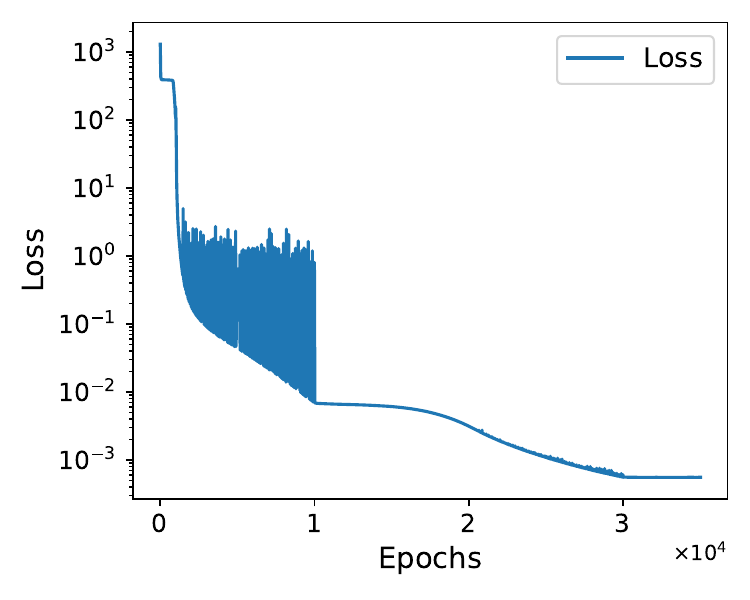}}
\subfigure[errors history]
 { 
\includegraphics[width=0.22\textwidth]{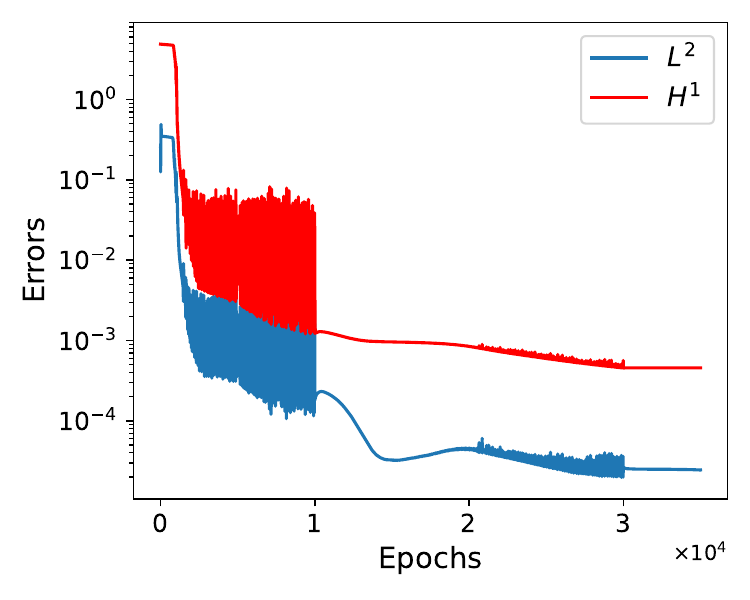}}
 
 \caption{ Numerical results for Example \ref{exm:1d_nonlinear} using {$N(x, (x-0.5)/\sqrt{\epsilon}, 1/\sqrt{\epsilon})$}:  \\
 \textbf{(a)-(d)}:  when $\epsilon=1/8$ $(k=8)$ , with $\alpha=1$, $\alpha_1=0$, $N_c=200$, epochs=35000; \\
 \textbf{(e)-(h)}: when $\epsilon=10^{-2} $ $(k=100)$, with $\alpha=1000$, $\alpha_1=0$, $N_c=450$, epochs=35000.}
 \label{fig:exm7_2s} %
 \end{figure}

\begin{figure} [!ht]
\centering
\subfigure[exact and NN solutions]
{
\includegraphics[width=0.25\textwidth]{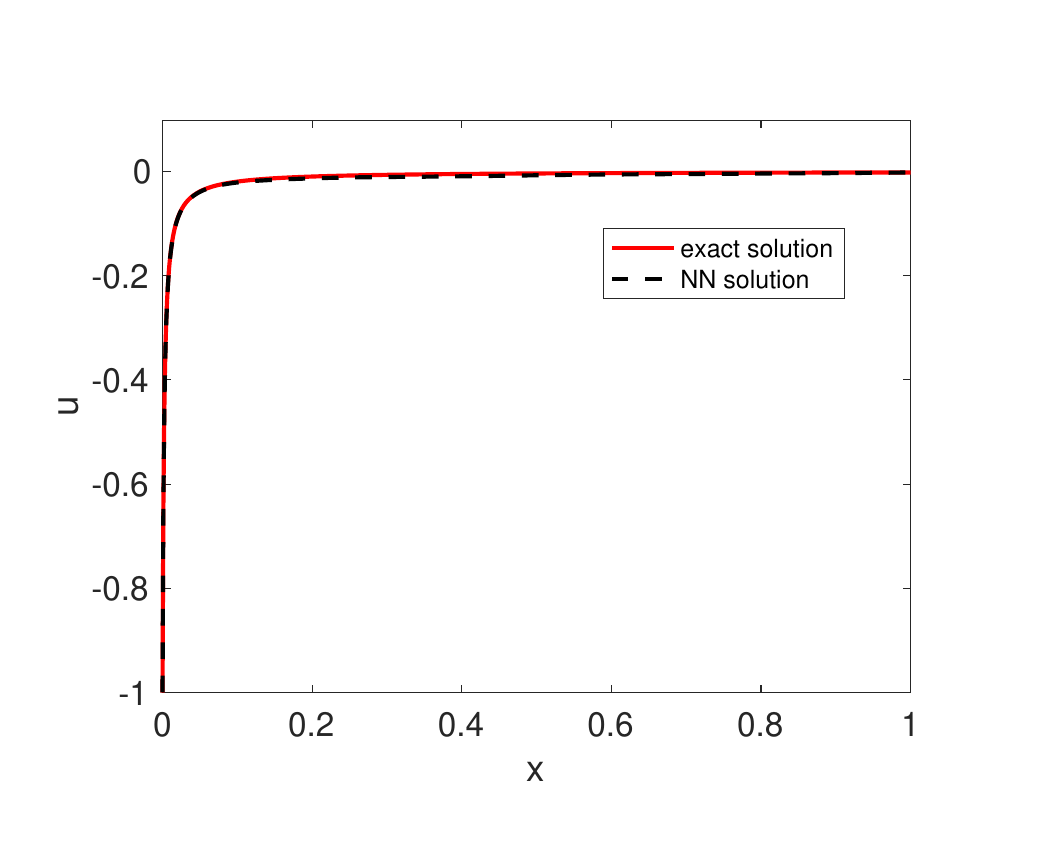}}
\subfigure[absolute error]
 { 
\includegraphics[width=0.25\textwidth]{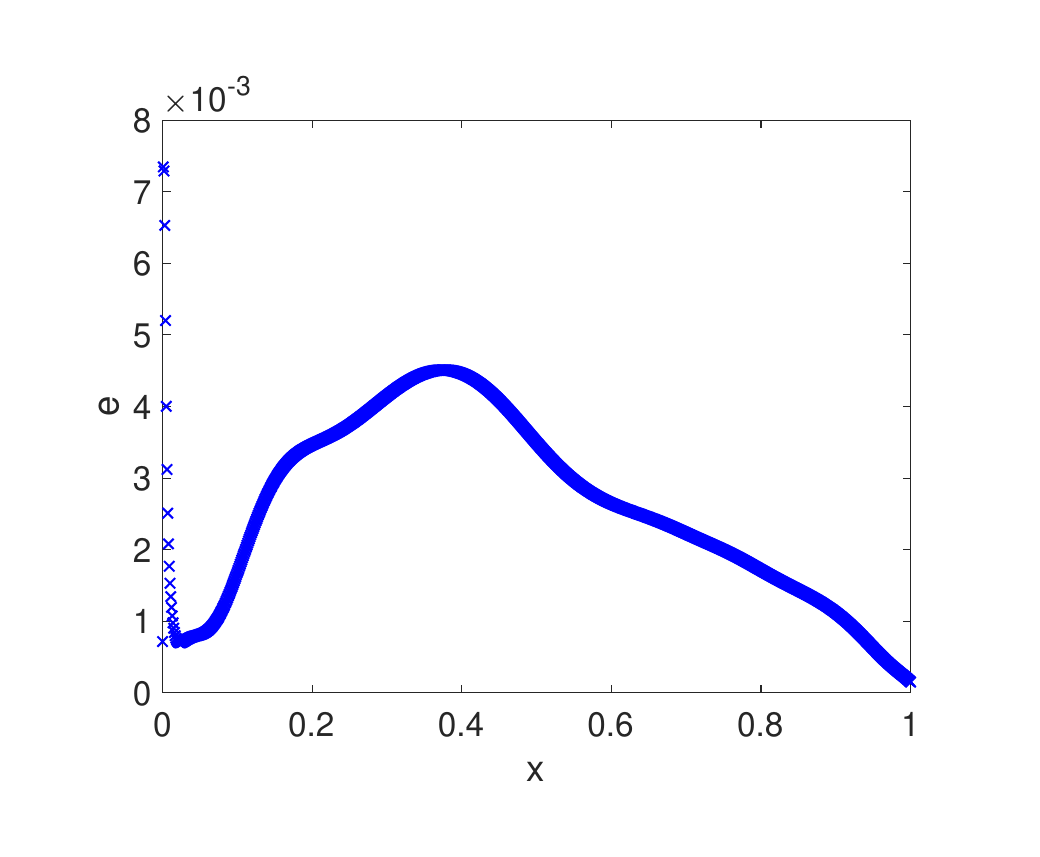}}
 \caption{ Numerical results for Example \ref{exm:1d_nonlinear}  when $\epsilon=10^{-3}$ using $N(x, (x-0.5)/\sqrt{\epsilon}, 1/\sqrt{\epsilon})$ and Algorithm \ref{alg:succesive-training-two-scale}, with parameters specified in Table \ref{tab:exm7_para}.}
 

\label{fig:exm7_2s_1e_3}
 \end{figure}

\begin{table}[ht]
\centering

\setlength{\tabcolsep}{11pt}
\begin{tabular}{c c c c c c}
\hline
$\epsilon$ & $10^{-2}$ (starting $\epsilon_0$)  & $5\times 10^{-3}$ & $2.5\times 10^{-3}$ & $1.25\times10^{-3}$ & $10^{-3}$ \\
\hline
$\alpha$ & 1000 & 1000 & 1000 & 10000 & 10000 \\
$\alpha_1$ & 0 & 0 & 0 & $0$ & $0$ \\
$N_c$ & 450 & 450 & 450 & 450 & 450 \\
LR & PC & $10^{-4}$ & $10^{-4}$ & $10^{-4}$ & $10^{-4}$ \\
epochs & 35000 & 35000 & 35000 & 35000 & 35000 \\
\hline
\end{tabular}
\caption{ Parameters in the loss function \eqref{eq:loss-general} and hyper-parameters for the successive training for Example \ref{exm:1d_nonlinear}, LR is the abbreviation for learning rate, PC is the piecewise constant scheduler in Table \ref{tab:lr_pt}.}
\label{tab:exm7_para}
\end{table}

}


\section{Conclusion and discussion}


In this work, we construct the two-scale neural networks by explicitly incorporating small parameters in { PDEs} into the architecture of feedforward neural networks. The construction enables solving problems with small parameters simply, without modifying the formulations of PINNs or adding  Fourier features in the networks. 
{ To enhance the accuracy of neural network predictions, especially for PDEs with smaller parameters, a successive training strategy is introduced.} 
Extensive numerical tests illustrate that the two-scale neural network method is effective and provides reasonable accuracy in capturing features associated with large derivatives in solutions. These features include boundary layers, inner layers, and oscillations.

Notably, when a parameter in the { PDEs} becomes extremely small, the method struggles to accurately capture large derivatives. 
However, for problems with extremely small parameters, the proposed two-scale neural network method may be employed as a good initial guess for further training. 

\section*{Acknowledgement}
{  We thank Professor Serge Prudhomme and Mr. Ziad Aldirany from Polytechnique Montr\'{e}al for providing the multi-level neural networks code for Examples 1 and 5.2 in their paper \cite{multilvNN23} and Example \ref{exm:1p14} a) in this work.
We also thank Professor Houman Owhadi
from the California Institute of Technology
 for helpful discussions}.
GEK acknowledges support from the DOE SEA-CROGS project (DE-SC0023191) and the MURI-AFOSR project (FA9550-20-1-0358). 


  \bibliographystyle{plain}
  \bibliography{ap-ref,singptb-ref,nn}
\end{document}